\documentclass{amsart}
\usepackage{graphicx} 
\usepackage[utf8]{inputenc}

\usepackage{amssymb}
\usepackage{amsmath}
\usepackage{amsthm}
\usepackage{enumitem}
\usepackage{pgf,tikz}
\usetikzlibrary{arrows}
\usetikzlibrary{decorations.markings}
\usetikzlibrary{backgrounds}
\usepackage{comment}

\numberwithin{equation}{section}

\theoremstyle{plain}
\newtheorem{thm}{Theorem}[section]
\newtheorem{lemma}[thm]{Lemma}
\newtheorem{prop}[thm]{Proposition}
\newtheorem{coroll}[thm]{Corollary}
\newtheorem{claim}[thm]{Claim}
\newtheorem{fact}[thm]{Fact}

\theoremstyle{definition}
\newtheorem{defn}[thm]{Definition}
\newtheorem{notation}[thm]{Notation}

\newtheorem{remark}[thm]{Remark}
\newtheorem{ex}[thm]{Example}
\newtheorem{prob}[thm]{Problem}

\newtheorem{innercustomthm}{Theorem}
\newenvironment{customthm}[1]
{\renewcommand\theinnercustomthm{#1}\innercustomthm}
{\endinnercustomthm}

\usepackage{xspace} 

\definecolor{e}{rgb}{0,0.6,0}
\definecolor{light-gray}{gray}{0.85}

\newcommand{\R}{\mathbb{R}}

\DeclareMathOperator{\Link}{Link}

\DeclareMathOperator{\jaeg}{\mathcal{J}}

\DeclareMathOperator{\blue}{bl}
\DeclareMathOperator{\green}{gr}

\def\jaeg{{\mathcal{J}}}
\def\jaegtree{\mathcal{J}{\rm_{tree}}}
\def\naivejaegtree{{\mathcal{J}{\rm_{tree}^{\rm naive}}}}
\def\quasi-trees{{\rm qt}}

\newcommand{\AuxContr}[1]{{\rm A}_{#1}}

\def\next{{\rm next}}
\newcommand{\Seifert}{\mathcal{S}}

\title{A quasi-tree expansion of the Homfly polynomial}

\author{Tam\'as K\'alm\'an}
\thanks{TK was supported by the Japan Society for the Promotion of Science (JSPS) Grant-in-Aid for Scientific Research C no.\ 23K03108.}
\address{Department of Mathematics,
Institute of Science Tokyo,
H-214, 2-12-1 Ookayama, Meguro-ku, Tokyo 152-8551, Japan}
\email{kalman@math.titech.ac.jp}

\author{Lilla T\'othm\'er\'esz}
\thanks{LT was supported by the Dynasnet European Research Council Synergy project (ERC2018-SYG 810115).}
\address{ELTE Eötvös Loránd University and HUN-REN Alfréd Rényi Institute of Mathematics, Budapest, Hungary}
\email{lilla.tothmeresz@ttk.elte.hu}

\begin{document}

\begin{abstract}
The Homfly polynomial is one of the best-studied invariants of oriented links. Here, building on work of Jaeger and others, we present a new method of computing it. Terms in the expansion are algebraically simple and they are indexed by certain quasi-trees of the Seifert graph of an (almost) arbitrary link diagram. The quasi-trees, in turn, depend on a single basepoint along the link but no more auxiliary data. This allows us to streamline the proofs of some old results and to find new ones as well, such as an upper bound on the total degree in terms of directed joins of the Seifert graph.
\end{abstract}

\maketitle

\section{Introduction}\label{s:intro}

\subsection{Setup and background}
In this paper we work with the standard notion of an oriented link in $\R^3$ and we represent links with planar diagrams in the usual way. See Section \ref{ss:links_and_graphs} for a short review.

The \emph{Homfly polynomial} of such an object, which is in fact an integer-coefficient Laurent polynomial in the indeterminates $v$ and $z$, is uniquely determined by the fact that it is invariant under isotopies of oriented links, that it associates $1$ to the trivial knot (also called the unknot), and, most importantly, its skein relation 
\[v^{-1}\,P_{\includegraphics[totalheight=8pt]{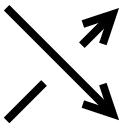}}-v\,P_{\includegraphics[totalheight=8pt]{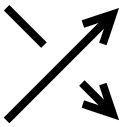}}=z\,P_{\includegraphics[totalheight=8pt]{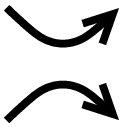}},\]
where, as usual, the subscripts mean link diagrams that contain the indicated parts but are otherwise identical. 

To compute $P$ 
for an oriented link given by a diagram,
one might successively apply the skein relation at positive crossings in the form 
\begin{equation}\label{eq:positive_skein}
    P_{\includegraphics[totalheight=8pt]{skeinpos.eps}}=v^2\,P_{\includegraphics[totalheight=8pt]{skeinneg.eps}}+vz\,P_{\includegraphics[totalheight=8pt]{skeinsmooth.eps}},
\end{equation} or at negative crossings in the form
\begin{equation}\label{eq:negative_skein}
P_{\includegraphics[totalheight=8pt]{skeinneg.eps}}=v^{-2}\,P_{\includegraphics[totalheight=8pt]{skeinpos.eps}}\,-\,v^{-1}z\,P_{\includegraphics[totalheight=8pt]{skeinsmooth.eps}}.
\end{equation}
That is, one reduces the computation of $P$ 
to two other diagrams, one of which is the result of \emph{crossing change} and the other the result of \emph{smoothing}.

Owing to the many available choices of crossing, this offers a multitude of 
ways of determining the Homfly polynomial. Indeed it is rather 
remarkable that the polynomial \emph{exists}, i.e., that all computations reach the same conclusion.

A complete computation of $P$ is achieved if, after branching a finite number of times, all resulting diagrams represent the trivial knot. 
In fact, trivial \emph{links} suffice because their Homfly polynomials can be handled by appealing to connected sums of unknots. This is what one finds in most of the literature but an important aspect of our work is that we fold this step into the rest of the computation. 

The standard strategy to create unknots (or unlinks) 
is (see, e.g., \cite{Jaeger,Kauffman}) to 
make all diagrams \emph{descending} with respect to some basepoint (or a collection of basepoints). This in turn means that whenever the curve reaches a crossing for the first time, it does so on the overcrossing strand; 
then, a``side view'' reveals that a knot with such a diagram is indeed trivial. 
A problem with working with unlinks, one that we will avoid, is that if we allow a descending component to close then a new basepoint has to be chosen for the next component, and so on. The first basepoint may be fixed throughout but the further ones may not, because a priori we do not know which parts of the curve will still be there to process.

The main goal of this paper is to present a streamlined version of such a computation, still using the idea of a descending diagram but further guided by certain graph-theoretical notions. 
Our expression contains fewer summands, and requires 
less auxiliary data, than its predecessors  \cite{Jaeger,Kauffman}. Moreover, by exploiting nontrivial graph theory, it delivers interesting knot-theoretical consequences such as Theorem \ref{thm:max_degree_dijoin} below.

Graphs enter through 
the Seifert graph of the link diagram, which we treat as a signed, directed ribbon graph (see Section \ref{ss:links_and_graphs}). 
Our formula, stated in Theorem \ref{thm:Homfly_formula}, presents the Homfly polynomial as a sum, where the terms correspond to certain quasi-trees of the Seifert graph, and 
are determined by certain activities of these quasi-trees. Beyond the link diagram itself, we use a single auxiliary datum: a vertex and incident edge of the Seifert graph (equivalently, a point on the link). Next, we informally introduce these notions and then state our result. For a more thorough introduction, see Section \ref{s:prelim}.
 
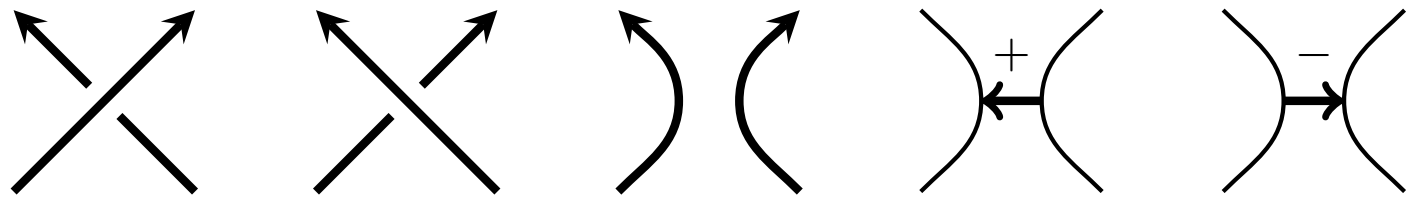
\begin{figure}[h]
    \centering
    \begin{tikzpicture}[scale=.1]
    \draw [ultra thick,->,>=stealth] (0,0)--(12,12);
    \draw [ultra thick] (12,0)--(7,5);
    \draw [ultra thick,->,>=stealth] (5,7)--(0,12);
        \begin{scope}[shift={(20,0)}]
        \draw [ultra thick,->,>=stealth] (12,0)--(0,12);
        \draw [ultra thick] (0,0)--(5,5);
        \draw [ultra thick,->,>=stealth] (7,7)--(12,12);
        \end{scope}
        \begin{scope}[shift={(40,0)}]
        \draw [ultra thick,->,>=stealth] (0,0) to [out=45,in=-90] (4,6) to [out=90,in=-45] (0,12);
        \draw [ultra thick,->,>=stealth] (12,0) to [out=135,in=-90] (8,6) to [out=90,in=-135] (12,12);
        \end{scope}
        \begin{scope}[shift={(60,0)}]
        \draw [thick] (0,0) to [out=45,in=-90] (4,6) to [out=90,in=-45] (0,12);
        \draw [thick] (12,0) to [out=135,in=-90] (8,6) to [out=90,in=-135] (12,12);
        \draw [ultra thick,<-] (4,6)--(8,6);
        \node at (6,9) {\small $+$};
        \end{scope}
        \begin{scope}[shift={(80,0)}]
        \draw [thick] (0,0) to [out=45,in=-90] (4,6) to [out=90,in=-45] (0,12);
        \draw [thick] (12,0) to [out=135,in=-90] (8,6) to [out=90,in=-135] (12,12);
        \draw [ultra thick,->] (4,6)--(8,6);
        \node at (6,9) {\small $-$};
        \end{scope}
    \end{tikzpicture}
    \caption{From left to right: positive crossing, negative crossing, orientation-preserving smoothing, edge of the Seifert graph corresponding to the positive crossing, edge of the Seifert graph corresponding to the negative crossing. 
    }
    \label{fig:posnegsmooth}
\end{figure}

The vertices of the Seifert graph
are the Seifert circles, obtained by smoothing each crossing in the diagram. 
The edges 
correspond to the crossings, 
and that endows them with a \emph{sign} as depicted in Figure \ref{fig:posnegsmooth}. Furthermore the \emph{direction} of each edge is induced by the under-strand at the crossing. One may think of the ribbon structure as an abstraction of the so called Seifert surface induced by the diagram; in particular the boundary of the Seifert graph (as a ribbon surface) is the link itself.

\emph{Quasi-trees} of a ribbon graph are spanning subgraphs with a connected boundary. In many respects quasi-trees, in the ribbon graph setting, are analogous to spanning trees of ordinary graphs. 
The main ingredients of Theorem \ref{thm:Homfly_formula}, Jaeger quasi-trees and their activities,
may be defined in any directed and signed ribbon graph after 
fixing one basepoint along the boundary (equivalently, a vertex and an incident edge). 
The main condition that distinguishes \emph{Jaeger quasi-trees} is that as we go around their boundary, every non-edge is first encountered at its tail. Apart from this we have one additional constraint that is less easy to state, but is computationally just as tractable; see Definition \ref{def:jaegquastree} for the details. In Section \ref{s:Jaeger_and_naive_Jaeger} we will also see that for some classes of link diagrams, the first condition is sufficient to define Jaeger quasi-trees. 

\subsection{Results}

One can sort edges of a Jaeger quasi-tree as either internally active or passive, and similarly, non-edges of a Jaeger quasi-tree as externally active or passive (see Section \ref{ss:first_ex} for informal descriptions and Definitions \ref{def:internal_activity_of_quasi-trees} and \ref{def:external_activity_of_quasi-trees} for the formal ones.) These concepts also involve the basepoint but they do not rely on any other auxiliary data such as an ordering of the edges.
With that, we claim the following formula. 

\begin{thm}\label{thm:Homfly_formula}
Let $L$ be an oriented link, and let $G$ be the 
Seifert graph of a diagram of $L$.  
Let $\jaeg(G,v_0,e_0)$ be the set of Jaeger quasi-trees of $G$ with respect to the basepoint $(v_0,e_0)$. If $\jaeg(G,v_0,e_0)\neq\varnothing$, then
\begin{multline*}
P_{L}(v,z)=\sum_{K\in\jaeg(G,v_0,e_0)} (v^2)^{ip^+(K)}\cdot (v^{-2})^{ip^-(K)}\cdot (vz)^{ep^+(K)}\cdot (-v^{-1}z)^{ep^-(K)} \\ \cdot (vz+vz^{-1}-v^3z^{-1})^{ea^+(K)}\cdot (-v^{-1}z-v^{-1}z^{-1}+v^{-3}z^{-1})^{ea^-(K)},    
\end{multline*}
where $ip^+(K)$ ($ip^-(K)$) denotes the number of internally active positive (negative) edges of $K$, $ep^+(K)$ ($ep^-(K)$) denotes the number of externally passive positive (negative) edges of $K$, and $ea^+(K)$ ($ea^-(K)$) denotes the number of externally active positive (negative) edges of $K$.
\end{thm}

The first four factors of each term come straight out of \eqref{eq:positive_skein} and \eqref{eq:negative_skein}, and the three-term expressions are exactly the Homfly polynomials of the positive and negative Hopf link, respectively. The latter, which we also abbreviate to $\Delta$ and $\bar\Delta$, respectively, arise as a result of grouping certain terms in the calculation and facilitate a more concise end result.

In nice cases such as so called homogeneous diagrams, which include positive and alternating diagrams, one does find Jaeger quasi-trees in the Seifert graph with respect to any basepoint.
Furthermore, any diagram can be modified using Reidemeister II moves so that its Seifert graph has a Jaeger quasi-tree for any basepoint, cf.\ Theorem \ref{thm:each_link_has_diagram_with_Jaeger_quasi-trees}, so requiring $\jaeg(G,v_0,e_0)\neq\varnothing$ is not a serious restriction. 

In Section \ref{ss:first_ex} we compute the formula of Theorem \ref{thm:Homfly_formula} on a small example, and in Section \ref{ss:intuition_behind_main_thm} we explain the main ideas behind it.

Let us say a few words about the computational complexity of the above formula. It is well known that computing the Homfly polynomial is hard \cite{Jaeger_Tutte_and_link}. Thus it may be no surprise that the number of Jaeger quasi-trees 
can be exponential in the number of crossings. On the positive side, for a given quasi-tree of the Seifert graph, deciding whether it is a Jaeger quasi-tree, and computing its activities, can be done in polynomial time in the number of edges of $G$, see Proposition \ref{prop:Jaeger_and_activity_computation_in_polynomial_time}.

Recall that the specialization $\nabla_L(z)=P_L(1,z)$ is the so called \emph{Conway polynomial} of $L$, which is in fact equivalent to the Alexander polynomial. By making that substitution in our main result, we obtain the following.

\begin{coroll}\label{cor:conway}
    If the oriented link $L$ has a connected diagram with only positive crossings, then for the Seifert graph $G$ of any such diagram, and any basepoint, the coefficient of $z^k$ in $\nabla_L(z)$, for any $k$, agrees with the number of Jaeger quasi-trees of $G$ consisting of $|E(G)|-k$ edges.
\end{coroll}

Here $E(G)$ is the set of edges of $G$, that is, the set of crossings of the link diagram.
As a quick check we note that since $G$ itself is a quasi-tree if and only if $L$ is a knot, and in that case it is Jaeger, Corollary \ref{cor:conway} (and indeed Theorem \ref{thm:Homfly_formula}, in the non-positive cases) is consistent with the well known fact that $\nabla(0)$ is $1$ for a knot and $0$ for a multi-component link.

Theorem \ref{thm:Homfly_formula} allows us to effectively estimate the \emph{total} degree of the Homfly polynomial in terms of the Seifert graph of the link diagram. Recall that for a graph a \emph{cut} is the set of edges with differently colored endpoints in some bicoloring of the vertices and it is a \emph{directed cut} if the heads of the edges all have the same color; moreover a \emph{directed join} is a set of edges that intersects each (non-empty) directed cut.

\begin{thm}\label{thm:max_degree_dijoin}
    Let $D$ be a link diagram such that its Seifert graph $G$ has a Jaeger quasi-tree with respect to some basepoint. Then the degree of $P_{D}$ (where by the degree of a monomial we mean the sum of its $v$-degree and $z$-degree) is at most $2(|E(G)|-\nu(G))$, where $\nu(G)$ is the minimum size of a directed join of $G$. If 
    $D$ is a positive diagram, then this estimate is tight. 
\end{thm}

Estimates for the $v$- and $z$-degrees \emph{separately} have been well known \cite{morton1986seifert} and we also address some of those in Section \ref{s:applications}. Furthermore, we construct a new and more transparent proof of Murasugi--Przytyczki's remarkable product formula \cite{MP_szorzodas}.

Our activity generating function from Theorem \ref{thm:Homfly_formula} is reminiscent of Morse's \cite{Morse_interlace} and Moffatt's \cite{moffatt2025constructing} quasi-tree generating function definitions of the ribbon graph polynomial (also called embedded graph Tutte polynomial or two-variable Bollobás--Riordan polynomial). In a special case (see Section \ref{ss:Homfly_well_def_on_subdivided_graphs}), we establish that our formula indeed gives the ribbon graph polynomial of a modified ribbon graph. This is a further example of the rich connection between knot polynomials and ribbon graph polynomials (see \cite{thistlethwaite1987spanning,chmutov2006kauffman,dasbach2008jones,moffatt2011unsigned,moffatt_knot_inv} for a variety of 
results of this type). 

We note that, though our motivation is the study of the Homfly polynomial, the formula of Theorem \ref{thm:Homfly_formula} can be computed for an arbitrary ribbon digraph. This gives rise to a range of interesting open questions concerning ribbon graphs.

\subsection*{Outline.}
In the remainder of this section, we outline the ideas behind Theorem \ref{thm:Homfly_formula} through an example.
In Section \ref{s:prelim}, we formally introduce the necessary definitions.
In Section \ref{s:Jaeger_trees_def_and_formula} we define Jaeger quasi-trees of arbitrary ribbon digraphs and supply a necessary and sufficient condition for their existence.
In Section \ref{s:activities} we define the activities of Jaeger quasi-trees and characterize them graph-theoretically.
In Section \ref{s:formula_proof} we apply the machinery developed in Sections \ref{s:Jaeger_trees_def_and_formula} and \ref{s:activities} to Seifert graphs of links, and prove Theorem \ref{thm:Homfly_formula}.
In Section \ref{s:applications} we illustrate the usefulness of Theorem \ref{thm:Homfly_formula} by giving an upper bound for the maximal total degree of the Homfly polynomial. We also provide simple new proofs for some old results concerning the Homfly polynomial.
In Section \ref{s:abstract_homfly} we 
evaluate the formula of Theorem \ref{thm:Homfly_formula} for signed ribbon digraphs that are not Seifert graphs of a link diagram, and see that in this case the result may depend on the chosen basepoint. Yet we do find another special case where the polynomial is well-defined, again with a topological explanation.
In Section \ref{s:Jaeger_and_naive_Jaeger} we 
re-examine the notion of Jaeger quasi-tree and establish cases, such as the Seifert graph of a special alternating link diagram, when the `naive' condition is sufficient to define them.

\subsection{An example}\label{ss:first_ex}

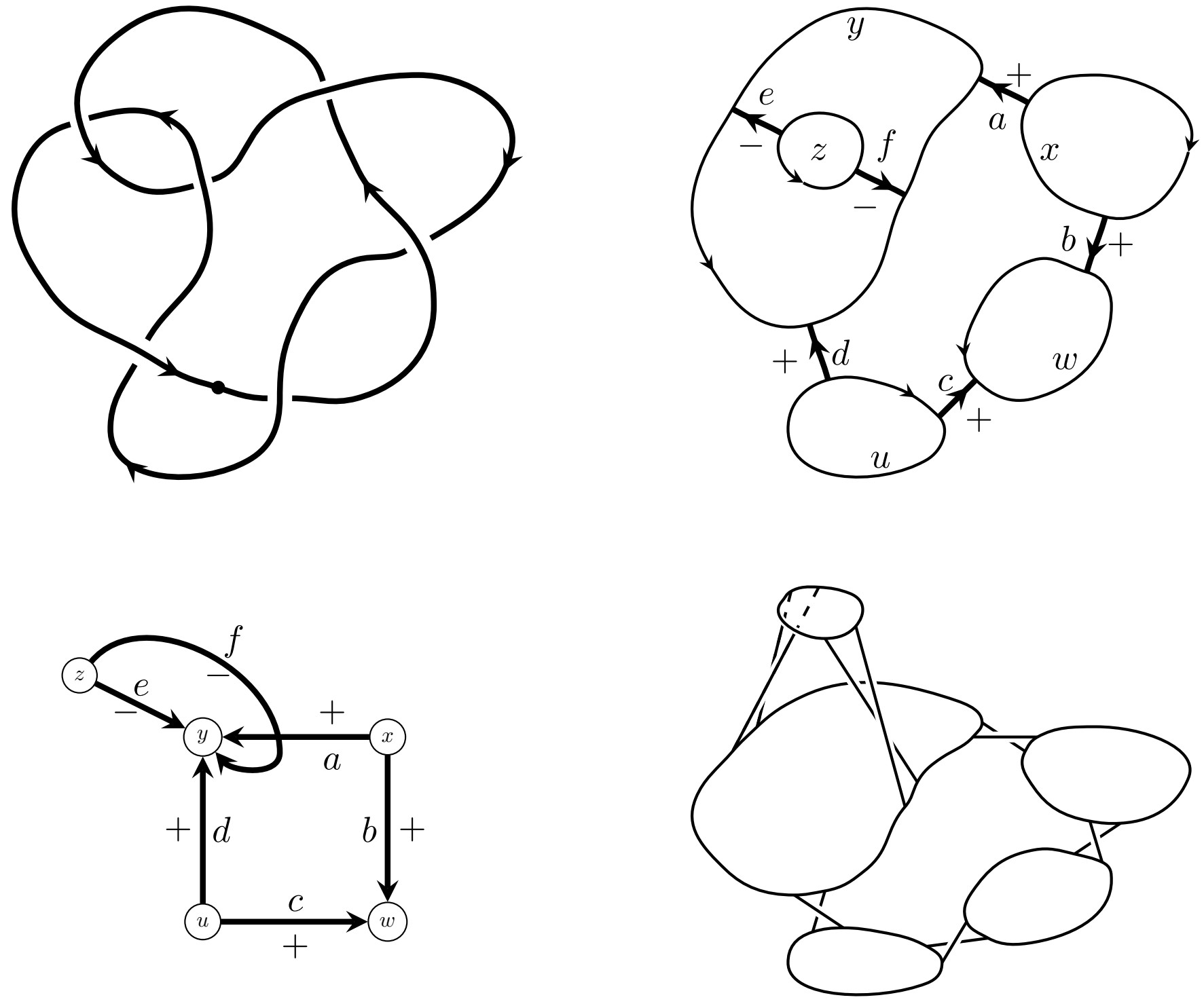
\begin{figure}[h]
    \centering
    \begin{tikzpicture}[decoration={markings,
    mark= at position 0.8 with {\arrow{stealth}}},>=stealth]
    \begin{scope}[scale=.3,shift={(0,0)}]
    \draw [ultra thick,postaction={decorate}] (4.2,4.9) to [out=60,in=-130] (5.5,6.5) to [out=50,in=-75] (6,10) to [out=105,in=-30] (5,12) to [out=150,in=15] (2.3,12.1);
    \draw [ultra thick,postaction={decorate}] (1.7,11.9) to [out=-165,in=75] (0,10) to [out=-105,in=125] (1,6.5) to [out=-55,in=150] (4,4.5) to [out=-30,in=165] (6,3.5) to [out=-15,in=180] (8.1,3);
    \draw [ultra thick,postaction={decorate}] (8.9,3) to [out=0,in=-165] (11,3) to [out=15,in=-90] (13.5,6) to [out=90,in=-60] (13,8) to [out=120,in=-65] (11,10.5) to [out=115,in=-75] (10.1,12.7);
    \draw [ultra thick,postaction={decorate}] (9.9,13.3) to [out=105,in=-25] (8,15) to [out=155,in=35] (3.5,15) to [out=-145,in=105] (2,12) to [out=-75,in=150] (3.5,10) to [out=-30,in=-165] (5.7,9.9);
    \draw [ultra thick,postaction={decorate}] (6.3,10.1) to [out=15,in=-135] (8,12) to [out=45,in=-165] (10,13) to [out=15,in=180] (13,13.5) to [out=0,in=75] (16,11) to [out=-105,in=30] (13.4,8.2);
    \draw [ultra thick,postaction={decorate}] (12.6,7.8) to [out=-150,in=15] (11,7.5) to [out=-165,in=65] (9,5.5) to [out=-115,in=90] (8.5,3) to [out=-90,in=30] (7.5,1) to [out=-150,in=-90] (3,2) to [out=90,in=-120] (3.8,4.1);
    \node [circle,fill,scale=0.4] at (6.5,3.35) {};
    \end{scope}
    \begin{scope}[scale=.3,shift={(22,0)}]
    \draw [->,thick] (3.5,10) to [out=-30,in=-120] (5.2,10.4) to [out=60,in=-30](5,12) to [out=150,in=60] (2.8,11.6) to [out=-120,in=150] (3.5,10); 
    \draw [thick,postaction={decorate}] (5.5,6.5) to [out=50,in=-120] (6.8,9.6) to [out=60,in=-135] (8,12) to [out=45,in=-120] (9.2,13.4) to [out=60,in=-25] (8,15) to [out=155,in=35] (3.5,15) to [out=-145,in=60] (1.2,12.4) to [out=-120,in=75] (0,10) to [out=-105,in=125] (1,6.5) to [out=-55,in=-165] (3.7,5.4) to [out=15,in=-130] (5.5,6.5);
    \draw [thick,postaction={decorate}] (7.5,1) to [out=-150,in=-90] (3,2) to [out=90,in=-165] (4.3,3.6) to [out=15,in=165] (6,3.5) to [out=-15,in=135] (7.9,2.4) to [out=-45,in=30] (7.5,1);
    \draw [thick,postaction={decorate}] (11,3) to [out=15,in=-90] (13.5,6) to [out=90,in=-15] (12.7,7.1) to [out=165,in=15] (11,7.5) to [out=-165,in=65] (9,5.5) to [out=-115,in=135] (9.1,3.6) to [out=-45,in=-165] (11,3);
    \draw [->,thick] (16,11) to [out=-105,in=-15] (13.3,8.9) to [out=165,in=-65] (11,10.5) to [out=115,in=-120] (10.8,12.6) to [out=60,in=180] (13,13.5) to [out=0,in=75] (16,11);
    \draw [ultra thick,postaction={decorate}] (2.8,11.6)  to [out=150,in=-30] (1.2,12.4);
    \draw [ultra thick,postaction={decorate}] (5.2,10.4)  to [out=-30,in=150] (6.8,9.6);
    \draw [ultra thick,postaction={decorate}] (4.3,3.6)  to [out=105,in=-75] (3.7,5.4);
    \draw [ultra thick,postaction={decorate}] (7.9,2.4) to [out=45,in=-135] (9.1,3.6);
    \draw [ultra thick,postaction={decorate}] (13.3,8.9)  to [out=-105,in=75] (12.7,7.1);
    \draw [ultra thick,postaction={decorate}] (10.8,12.6) to [out=150,in=-30] (9.2,13.4);
    \node at (6,1) {$u$};
    \node at (12,4.2) {$w$};
    \node at (11.5,11) {$x$};
    \node at (4,11) {$z$};
    \node at (5.2,15) {$y$};
    \node at (10.5,13.5) {$+$};
    \node at (13.8,8) {$+$};
    \node at (9.2,2.3) {$+$};
    \node at (2.9,4.2) {$+$};
    \node at (1.8,11.2) {$-$};
    \node at (5.5,9.2) {$-$};
    \node at (9.8,12) {$a$};
    \node at (12.1,8.2) {$b$};
    \node at (8.1,3.5) {$c$};
    \node at (4.7,4.5) {$d$};
    \node at (2.3,12.8) {$e$};
    \node at (6.2,11.2) {$f$};
    \end{scope}
    \begin{scope}[scale=.3,shift={(0,-20)}]
    \node [circle,scale=.6,draw] (a) at (6,6) {$u$};
    \node [circle,scale=.6,draw] (b) at (12,6) {$w$};
    \node [circle,scale=.6,draw] (c) at (12,12) {$x$};
    \node [circle,scale=.6,draw] (d) at (6,12) {$y$};
    \node [circle,scale=.6,draw] (e) at (2,14) {$z$};
    \draw [<-,ultra thick] (d) to (e);
    \draw [<-,ultra thick] (d) to [out=-50,in=-90] (8.5,11.5) to [out=90,in=50] (e);
    \draw [->,ultra thick] (a) to (d);
    \draw [->,ultra thick] (a) to (b);
    \draw [->,ultra thick] (c) to (b);
    \draw [->,ultra thick] (c) to (d);
    \node [circle] at (3.5,12.8) {$-$};
    \node [circle] at (6.5,14) {$-$};
    \node [circle] at (10.2,12.8) {$+$};
    \node [circle] at (12.8,9) {$+$};
    \node [circle] at (9,5.2) {$+$};
    \node [circle] at (5.2,9) {$+$};
    \node at (10.2,11.2) {$a$};
    \node at (11.4,9) {$b$};
    \node at (9,6.6) {$c$};
    \node at (6.6,9) {$d$};
    \node at (4,13.6) {$e$};
    \node at (7,15.1) {$f$};
    \end{scope}
    \begin{scope}[xscale=.3,yscale=.2,shift={(22,-25)}]
        \draw [thick] (2.9,18.9) to (2,13.5);
	\draw [dashed,thick] (3.1,20.1) to (2,13.5);
    \draw [line width=3, white] (1.9,14.3) -- (3.3,18.3);
	\draw [dashed,thick] (4,20.3) to (1.2,12.3);
	\draw [thick] (3.3,18.3) to (1.2,12.3);
    \draw [line width=2, white] (2.8,19.6) to [out=-120,in=150] (3.5,18);
	\draw [thick] (5,20) to [out=150,in=-5] (4,20.3) to [out=175,in=40] (3.1,20.1) to [out=-140,in=60] (2.8,19.6) to [out=-120,in=150] (3.5,18) to [out=-30,in=180] (4.2,17.8) to [out=0,in=-120] (5.2,18.4) to [out=60,in=-30] (5,20);
	\draw [thick] (5.5,6.5) to [out=50,in=-120] (6.8,9.6) to [out=60,in=-110] (7.2,10.8) to [out=70,in=-135] (8,12) to [out=45,in=-120] (9,13) to [out=60,in=-90] (9.3,13.7) to [out=90,in=-25] (8,15) to [out=155,in=35] (3.5,15) to [out=-145,in=55] (2,13.5) to [out=-125,in=60] (1.2,12.3) to [out=-120,in=75] (0,10) to [out=-105,in=125] (1,6.5) to [out=-55,in=180] (3.2,5.3) to [out=0,in=-165] (4.2,5.5) to [out=15,in=-130] (5.5,6.5);
	\draw [thick] (7.5,1) to [out=-150,in=-90] (3,2) to [out=90,in=-155] (3.8,3.5) to [out=25,in=180] (4.8,3.7) to [out=0,in=165] (6,3.5) to [out=-15,in=140] (7.5,2.8) to [out=-40,in=90] (8,2) to [out=-90,in=30] (7.5,1);
	\draw [thick] (11,3) to [out=15,in=-90] (13.5,6) to [out=90,in=-35] (13.2,6.9) to [out=145,in=-20] (12.3,7.3) to [out=160,in=15] (11,7.5) to [out=-165,in=65] (9,5.5) to [out=-115, in=105] (8.8,4) to [out=-75,in=140] (9.5,3.2) to [out=-40,in=-165] (11,3);
	\draw [thick] (11,10.5) to [out=115,in=-110] (10.7,12.3) to [out=70,in=-140] (11.1,13) to [out=40,in=180] (13,13.5) to [out=0,in=75] (16,11) to [out=-105,in=5] (13.8,8.8) to [out=-175,in=-15] (12.8,8.9) to [out=165,in=-65] (11,10.5);
    \draw [line width=3, white] (4.8,16.4) -- (5.4,15);
    \draw [thick] (4.2,17.8) to (7.2,10.8);
    \draw [line width=3, white] (5.6,16.2) -- (6.4,11.8);
	\draw [thick] (5.2,18.4) to (6.8,9.6);
    \draw [thick] (9.3,13.7) -- (10.7,12.3);
    \draw [line width=3, white] (9.5,13) -- (10.6,13);
    \draw [thick] (9,13) -- (11.1,13);
    \draw [thick] (13.8,8.8) -- (12.3,7.3);
    \draw [line width=3, white] (12.9,8.4) -- (13.1,7.4);
    \draw [thick] (12.8,8.9) -- (13.2,6.9);
    \draw [thick] (9.5,3.2) -- (7.5,2.8);
    \draw [line width=3, white] (8.6,3.5) -- (8.2,2.5);
    \draw [thick] (8.8,4) -- (8,2);
    \draw [thick] (3.8,3.5) -- (4.2,5.5);
    \draw [line width=3, white] (4.3,4.2) -- (3.7,4.8);
    \draw [thick] (4.8,3.7) -- (3.2,5.3);
    \end{scope}
    \end{tikzpicture}
    \caption{Upper left: a diagram of the twist knot $6_1$. Upper right: The Seifert circles of the diagram and the directed edges of its Seifert graph. Lower left: An abstract presentation of the Seifert graph.
    Lower right: The ribbon surface of the ribbon structure, embedded in three-space as a Seifert surface of the knot.} 
    \label{fig:6_1_link_diagram_and_Seifert}
\end{figure}

Let us compute the formula of Theorem \ref{thm:Homfly_formula} for the standard diagram of the $6_1$ knot, shown in the upper left panel of Figure \ref{fig:6_1_link_diagram_and_Seifert}.
The upper right panel of Figure \ref{fig:6_1_link_diagram_and_Seifert} shows the Seifert graph in our case.
Its edges correspond to the six crossings and each is directed by the undercrossing strand. 
The lower left panel shows the same 
graph, drawn in a more conventional way and so that the positive orientation of the plane induces its ribbon structure. 

For the computation we need to choose a basepoint, which will be the point marked by a dot on the upper left panel of Figure \ref{fig:6_1_link_diagram_and_Seifert}. Equivalently, this corresponds to picking the vertex-edge pair $(u,c)$ in the Seifert graph.

By definition, one may continuously traverse the entire boundary of a quasi-tree starting from the basepoint (more formally, we call this the tour of the quasi-tree). Now a quasi-tree is \emph{naive Jaeger} if, during its tour, each non-edge is first reached at its tail. Jaeger quasi-trees form a subset of naive Jaeger quasi-trees,
defined by a stronger and more complicated condition that we defer to
Definition \ref{def:jaegquastree}. However, the sets of Jaeger and naive Jaeger quasi-trees often 
coincide, and that includes this example.

\begin{figure}
\begin{tikzpicture}[>=stealth,scale=.28]
        \centering
    \begin{scope}[shift={(-15,0)}]
    \node at (9,9) {\small $K_1$};
    \node [circle,scale=.8,draw] (a) at (6,6) {$u$};
    \node [circle,scale=.8,draw] (b) at (12,6) {$w$};
    \node [circle,scale=.8,draw] (c) at (12,12) {$x$};
    \node [circle,scale=.8,draw] (d) at (6,12) {$y$};
    \node [circle,scale=.8,draw] (e) at (2,14) {$z$};
    \draw [color=red,<-,ultra thick,>=stealth,dotted] (d) to (e);
    \draw [color=red,<-,ultra thick] (d) to [out=-48,in=-90] (8.5,11.5) to [out=90,in=50] (e);
    \draw [color=red,->,ultra thick] (a) to (d);
    \draw [color=blue,->,ultra thick,dotted] (a) to (b);
    \draw [color=red,->,ultra thick] (c) to (b);
    \draw [color=red,->,ultra thick] (c) to (d);
    \node [circle] at (3.5,12.8) {\small $-$};
    \node [circle] at (6.5,14) {\small $-$};
    \node [circle] at (10.2,12.8) {\small $+$};
    \node [circle] at (12.8,9) {\small $+$};
    \node [circle] at (9,5.2) {\small $+$};
    \node [circle] at (5.2,9) {\small $+$};
    \node at (10.2,11.2) {\small $a$};
    \node at (11.4,9) {\small $b$};
    \node at (9,6.6) {\small $c$};
    \node at (6.6,9) {\small $d$};
    \node at (4,13.6) {\small $e$};
    \node at (7,15.1) {\small $f$};
    \node [circle,fill,scale=.3,color=gray] at (6,4.8) {};
    \draw [->,color=gray] (6,4.8) to (7,4.8);
    \end{scope}
    \begin{scope}[shift={(0,0)}]
    \node at (9,9) {\small $K_2$};
    \node [circle,scale=.8,draw] (a) at (6,6) {$u$};
    \node [circle,scale=.8,draw] (b) at (12,6) {$w$};
    \node [circle,scale=.8,draw] (c) at (12,12) {$x$};
    \node [circle,scale=.8,draw] (d) at (6,12) {$y$};
    \node [circle,scale=.8,draw] (e) at (2,14) {$z$};
    \draw [color=red,<-,ultra thick,>=stealth,dotted] (d) to (e);
    \draw [color=red,<-,ultra thick] (d) to [out=-48,in=-90] (8.5,11.5) to [out=90,in=50] (e);
    \draw [color=red,->,ultra thick] (a) to (d);
    \draw [color=blue,->,ultra thick] (a) to (b);
    \draw [color=red,->,ultra thick] (c) to (b);
    \draw [color=blue,->,ultra thick,dotted] (c) to (d);
    \node [circle] at (3.5,12.8) {\small $-$};
    \node [circle] at (6.5,14) {\small $-$};
    \node [circle] at (10.2,12.8) {\small $+$};
    \node [circle] at (12.8,9) {\small $+$};
    \node [circle] at (9,5.2) {\small $+$};
    \node [circle] at (5.2,9) {\small $+$};
    \node at (10.2,11.2) {\small $a$};
    \node at (11.4,9) {\small $b$};
    \node at (9,6.6) {\small $c$};
    \node at (6.6,9) {\small $d$};
    \node at (4,13.6) {\small $e$};
    \node at (7,15.1) {\small $f$};
    \node [circle,fill,scale=.3,color=gray] at (6,4.8) {};
    \draw [->,color=gray] (6,4.8) to (7,4.8);
    \end{scope}
    \begin{scope}[shift={(15,0)}]
    \node at (9,9) {\small $K_3$};
    \node [circle,scale=.8,draw] (a) at (6,6) {$u$};
    \node [circle,scale=.8,draw] (b) at (12,6) {$w$};
    \node [circle,scale=.8,draw] (c) at (12,12) {$x$};
    \node [circle,scale=.8,draw] (d) at (6,12) {$y$};
    \node [circle,scale=.8,draw] (e) at (2,14) {$z$};
    \draw [color=red,<-,ultra thick] (d) to (e);
    \draw [color=red,<-,ultra thick] (d) to [out=-48,in=-90] (8.5,11.5) to [out=90,in=50] (e);
    \draw [color=red,->,ultra thick] (a) to (d);
    \draw [color=blue,->,ultra thick] (a) to (b);
    \draw [color=red,->,ultra thick] (c) to (b);
    \draw [color=blue,->,ultra thick] (c) to (d);
    \node [circle] at (3.5,12.8) {\small $-$};
    \node [circle] at (6.5,14) {\small $-$};
    \node [circle] at (10.2,12.8) {\small $+$};
    \node [circle] at (12.8,9) {\small $+$};
    \node [circle] at (9,5.2) {\small $+$};
    \node [circle] at (5.2,9) {\small $+$};
    \node at (10.2,11.2) {\small $a$};
    \node at (11.4,9) {\small $b$};
    \node at (9,6.6) {\small $c$};
    \node at (6.6,9) {\small $d$};
    \node at (4,13.6) {\small $e$};
    \node at (7,15.1) {\small $f$};
    \node [circle,fill,scale=.3,color=gray] at (6,4.8) {};
    \draw [->,color=gray] (6,4.8) to (7,4.8);
    \end{scope}
    \end{tikzpicture}
    \caption{The Jaeger quasi-trees (solid edges) of the Seifert graph of the standard diagram of the $6_1$ knot, with basepoint $(u,c)$. Red edges are active, blue edges are passive.}
    \label{fig:61_Jaeger_kvazifai}
\end{figure}
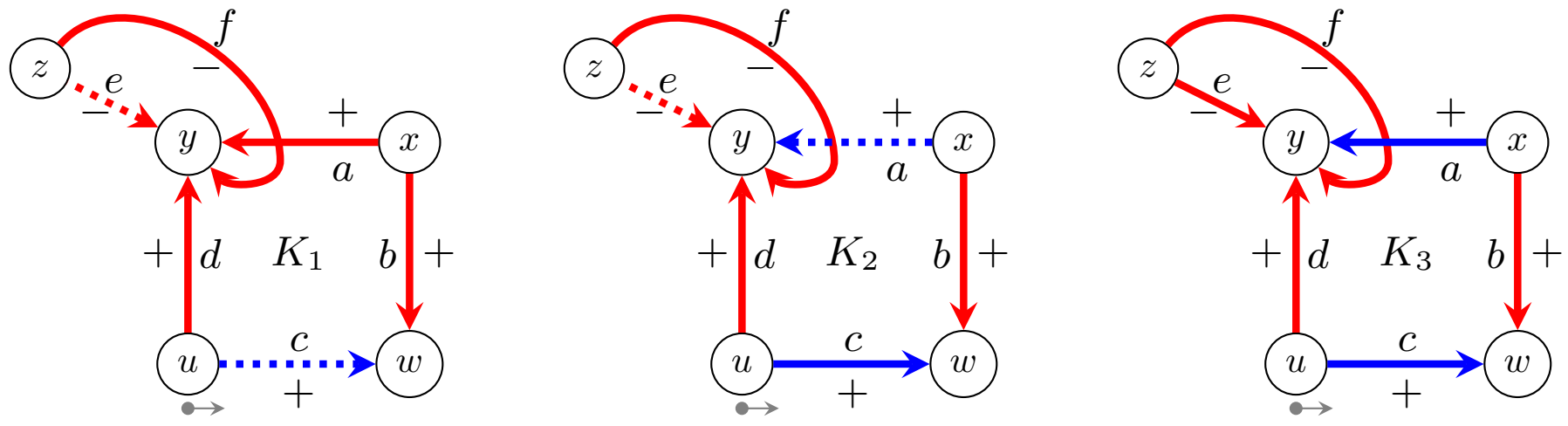

With $(u,c)$ as basepoint, our Seifert graph has three Jaeger quasi-trees. 
These can be seen in Figure \ref{fig:61_Jaeger_kvazifai}. Two of them are spanning trees, and the third one is the entire graph. 
Note that in this presentation of the Seifert graph and its subgraphs, a tour can be visualized by tracing the boundary of a slightly thickened version of the planar image, starting from the point indicated at the lower left. (See Section \ref{ss:quasi-trees_and_tours} for a more formal introduction and examples.)

For any two quasi-trees, their tours agree up to some point, when there is a first edge that is in their symmetric difference. 
We call this the
\emph{first difference}; 
see the more formal definition in Section \ref{ss:quasi-trees_and_tours}. We say that an edge/nonedge of a Jaeger quasi-tree $K$ is \emph{internally/externally active} if it is \emph{not} a first difference between $K$ and any other Jaeger quasi-tree. An internal/external edge is \emph{passive} if it is not active. 

For example, the tours of $K_2$ and $K_3$ in Figure \ref{fig:61_Jaeger_kvazifai}, starting from $(u,c)$, both travel along the edge $c$ first, then go on to $b$, then reach $a$. Now $a\in K_3-K_2$, whence this is the first difference between the tours of $K_2$ and $K_3$. Consequently, $a$ is externally passive in $K_2$ and internally passive in $K_3$. In Figure \ref{fig:61_Jaeger_kvazifai}, active edges are drawn in red (solid if the edge is internal, dotted if external) and passive edges in blue.

The following table shows the activities of our Jaeger quasi-trees: 
\[\begin{array}{c|c|c|c|c|c|c|c|c}
     & ia^+ & ia^- & ip^+ & ip^- & ea^+ & ea^- & ep^+ & ep^- \\
     \hline
     K_1 & 3 & 1 & 0 & 0 & 0 & 1 & 1 & 0 \\
     \hline
     K_2 & 2 & 1 & 1 & 0 & 0 & 1 & 1 & 0 \\
     \hline
     K_3 & 2 & 2 & 2 & 0 & 0 & 0 & 0 & 0      
\end{array}\]
According to Theorem \ref{thm:Homfly_formula}, the contribution of $K_1$ to the Homfly polynomial of $6_1$ is $vz \cdot (-v^{-1}z-v^{-1}z^{-1}+v^{-3}z^{-1})$. The contribution of $K_2$ is $v^2 \cdot vz \cdot (-v^{-1}z-v^{-1}z^{-1}+v^{-3}z^{-1})$, while the contribution of $K_3$ is $v^4$. Summing these up, the constant term cancels and we get that the Homfly polynomial of $6_1$ is 
\begin{equation*}
\begin{array}{cccc}
& -z^2 & -v^2z^2 & \\
    +v^{-2} & & -v^2 & +v^4. 
\end{array}
\end{equation*}

\subsection{Sketch of proof}\label{ss:intuition_behind_main_thm} 
Theorem \ref{thm:Homfly_formula} is established via a computation tree, an example of which is shown
in Figure \ref{fig:computationtree}.
Our strategy is to trace the link, starting from the basepoint, and at certain crossings apply the skein relation. With that we successively reduce the computation of the Homfly polynomial $P$ to computing 
$P$ of two ``simpler'' links. 
At the level of Seifert graphs, applying 
\eqref{eq:positive_skein} or \eqref{eq:negative_skein} 
results in two new 
graphs, one in which the edge corresponding to the crossing in question is reversed, and another one, where the edge is deleted.

It all boils down to deciding when to apply the skein relation. 
A first idea could be to look for descending, thus trivial, links. 
To achieve this we could set the rule, as Jaeger \cite{Jaeger} and Kauffman \cite{Kauffman} do, that 
when we arrive at a crossing for the first time on the understrand, then we always apply the skein relation. This, in both resulting diagrams, 
would eliminate the ``bad'' crossing. 
However, it might happen that we return to the basepoint without traversing the entire diagram.
Hence, after finishing one component, one would need to choose another basepoint on some other component in order to address the rest of the link. This 
would require additional auxiliary data, cf.\ \cite{Jaeger}. 

We solve this problem, and at the same time reduce the overall amount of computation, by sometimes ignoring a crossing even if we reached it on the understrand.
By doing this judiciously we are able to get, instead of descending unlinks, non-descending unknots in the form of connected sums of descending unknots.

Our decisions of when to perform the skein relation are guided by Jaeger quasi-trees. We start with a diagram and choose a basepoint so that at least one Jaeger quasi-tree exists.
This can always be guaranteed by some small modifications of the link diagram; see Section \ref{ss:treeexistence} for the details.

We perform a walk along the link 
and when we reach a crossing (say, corresponding to the edge $e$ of the current Seifert graph $G$) 
for the first time, 
we ask if the removal of $e$ from $G$ and/or the inclusion of $e$ in the 
already visited subgraph of $G$ 
results in a graph with an initial walk that can be completed to the tour of a Jaeger quasi-tree.
\begin{itemize}
\item If the answer is yes in both cases, then we apply the skein relation \eqref{eq:positive_skein} or \eqref{eq:negative_skein} to the crossing. This results in a Seifert graph where $e$ is deleted, and another one where $e$ is reversed. 
The Jaeger quasi-trees 
whose tours agree with our walk up to (and including) $e$ are split up between the two computation branches. The ones on the branch where $e$ was deleted have $e$ as an \textbf{externally passive} edge, and the ones on the branch where $e$ was kept have $e$ as an \textbf{internally passive} edge.
\item If only inclusion works then we leave the edge/crossing as is, and walk on to the next one. (In particular, this case occurs every time an edge is first reached at its head, i.e., on the overstrand. But it also applies in other instances, which is where we deviate from 
the descending unlink approach.) 
In this case, clearly, $e$ is \textbf{internally active} in each Jaeger quasi-tree 
whose tour agrees with our walk up to $e$.
\item Finally, in the case when only the removal of the edge $e$ is feasible, we do that by smoothing the corresponding crossing. We \emph{do not} consider crossing change, instead
we apply some special weights. 
Namely, 
$P_{L_G}=\Delta\cdot P_{L_{G-e}}$ if $e$ is positive, and $P_{L_G}=\bar\Delta\cdot P_{L_{G-e}}$ if $e$ is negative. 
We explain this through an example in Section \ref{ss:ex_cont}, cf.\ Figure \ref{fig:shorthand}, and then in Lemma \ref{l:handling_edges_that_need_to_be_removed} in general.  
This idea is a rather simple computational shortcut, independent of the previous point, yet it represents a further significant reduction in the number of diagrams to be considered.
In this case $e$ is obviously \textbf{externally active} in each Jaeger quasi-tree 
whose tour agrees with our walk up to $e$.
\end{itemize}

In all cases and resulting diagrams/Seifert graphs, the walk resumes continuously, and it stops when it returns to the basepoint. 
Note that at the end of each branch of this process, there is exactly one Jaeger quasi-tree---more precisely, a subgraph of the original unoriented Seifert graph, whose orientation is obtained by reversing the internally passive edges of a Jaeger quasi-tree. Lemma \ref{l:Jaeger_quasi-trees_link_trivial} says that the links corresponding to these Seifert graphs are all trivial knots. Therefore the method achieves a full computation of the Homfly polynomial of the original input, in the form given in Theorem \ref{thm:Homfly_formula}.

\subsection{Example, continued}\label{ss:ex_cont}

\begin{figure}[h]
    \centering
    \begin{tikzpicture}[scale=.15,decoration={markings,
    mark= at position 0.6 with {\arrow{stealth}}},>=stealth]
    \tikzstyle{o}=[circle,scale=.6,draw];

    \draw (-3,-16) rectangle (18,17);
    \draw (27,-16) rectangle (48,17);
    \draw[thick] (57,-16) rectangle (78,17);
    \draw (-3,-56) rectangle (18,-23);
    \draw (27,-56) rectangle (48,-23);
    \draw[thick] (57,-56) rectangle (78,-23);
    \draw[thick] (-3,-96) rectangle (18,-63);
    \draw [->] (20,0) -- (26,0);
    \node at (23,2) {\small $vz$};
    \draw [->] (50,0) -- (56,0);
    \node at (53,2) {\small $\bar\Delta$};
    \draw [->] (20,-35) -- (26,-35);
    \node at (23,-33) {\small $vz$};
    \draw [->] (50,-35) -- (56,-35);
    \node at (53,-33) {\small $\bar\Delta$};
    \draw [->] (6,-18) -- (6,-22);
    \node at (8,-20) {\small $v^2$};
    \draw [->] (6,-58) -- (6,-62);
    \node at (8,-60) {\small $v^2$};

\node at (50,-70) {$\left(\bar\Delta=\begin{aligned}
    &-v^{-1}z\\+v^{-3}z^{-1}&-v^{-1}z^{-1}
\end{aligned}\right)$};

\begin{scope}[shift={(0,0)}]
\node at (14,2) {\small $A$};
    \draw [red!30,line width=2,->] (4,11) to (6,10) to (8,9);
    \draw [red!30,line width=2,->] (4,11) to (2,12) to (0,13);
    \draw [red!30,line width=2,->] (5.5,2) to [out=115,in=-70] (4,4.5) to [out=110,in=-80] (3.5,6);
    \draw [red!30,line width=2,->] (5.5,2) to [out=20,in=-155] (8.5,3) to [out=25,in=-140] (11,5);
    \draw [red!30,line width=2,->] (13.5,11) to [out=-95,in=70] (13,8) to [out=-110,in=40] (11,5);
    \draw [red!30,line width=2,->] (13.5,11) to [out=150,in=-20] (10,13) to [out=160,in=-10] (7.5,13.5);
    \draw [green!30,fill] (5.5,2) circle [radius=.6];
    \draw [blue!30,fill] (11,5) circle [radius=.6];
    \draw [green!30,fill] (13.5,11) circle [radius=.6];
    \draw [green!30,fill] (4,11) circle [radius=.6];
    \draw [blue!30,line width=3] (5.5,6.5) to [out=25,in=-110] (8,9) to [out=70,in=-75] (8,12) to [out=105,in=-60] (7.5,13.5) to [out=120,in=15] (3.5,15) to [out=-165,in=65] (0,13) to [out=-115,in=100] (0,10) to [out=-80,in=180] (3.5,6) to [out=0,in=-155] (5.5,6.5);
    \draw [black!40,ultra thick,
    ] (4.2,4.9) to [out=60,in=-130] (5.5,6.5) to [out=50,in=-75] (6,10) to [out=105,in=-30] (5,12) to [out=150,in=15] (2.3,12.1);
    \draw [black!40,ultra thick,postaction={decorate}] (1.7,11.9) to [out=-165,in=75] (0,10) to [out=-105,in=125] (1,6.5) to [out=-55,in=150] (4,4.5) to [out=-30,in=165] (6,3.5);
    \draw [ultra thick] (6,3.5) to [out=-15,in=180] (8.1,3);
    \draw [black!40,ultra thick,postaction={decorate}] (8.9,3) to [out=0,in=-165] (11,3) to [out=15,in=-90] (13.5,6) to [out=90,in=-60] (13,8) to [out=120,in=-65] (11,10.5) to [out=115,in=-75] (10.1,12.7);
    \draw [black!40,ultra thick,postaction={decorate}] (9.9,13.3) to [out=105,in=-25] (8,15) to [out=155,in=35] (3.5,15) to [out=-145,in=105] (2,12) to [out=-75,in=150] (3.5,10) to [out=-30,in=-165] (5.7,9.9);
    \draw [black!40,ultra thick,postaction={decorate}] (6.3,10.1) to [out=15,in=-135] (8,12) to [out=45,in=-165] (10,13) to [out=15,in=180] (13,13.5) to [out=0,in=75] (16,11) to [out=-105,in=30] (13.4,8.2);
    \draw [black!40,ultra thick,
    ] (12.6,7.8) to [out=-150,in=15] (11,7.5) to [out=-165,in=65] (9,5.5) to [out=-115,in=90] (8.5,3) to [out=-90,in=30] (7.5,1) to [out=-150,in=-90] (3,2) to [out=90,in=-120] (3.8,4.1);
    \draw [fill] (6,3.5) circle [radius=.3];
    \end{scope}

    \begin{scope}[>=stealth,shift={(-2,-20)},scale=1.2]
    \node at (9,9) {\small $A$};
    \node [o,green] (u) at (6,6) {\color{black} $u$};
    \node [o,blue] (w) at (12,6) {\color{black} $v$};
    \node [o,green] (x) at (12,12) {\color{black} $x$};
    \node [o,blue] (y) at (6,12) {\color{black} $y$};
    \node [o,green] (z) at (3,14) {\color{black} $z$};
    \draw [red!30,<-] (y) to (z);
    \draw [red!30,<-] (y) to [out=-48,in=-90] (9,12) to [out=90,in=50] (z);
    \draw [red!30,->] (u) to (y);
    \draw [red,->,ultra thick] (u) to (w);
    \draw [red!30,->] (x) to (w);
    \draw [red!30,->] (x) to (y);
    \node [circle,fill,scale=.3,color=gray] at (6,4.5) {};
    \draw [->,color=gray] (6,4.5) to (7,4.5);
    \node at (10.2,12.8) {\small $a$};
    \node at (12.6,9) {\small $b$};
    \node at (9,5.2) {\small $c$};
    \node at (5.2,8.8) {\small $d$};
    \node at (3.8,12.6) {\small $e$};
    \node at (7,15) {\small $f$};
    \end{scope}
    
    \begin{scope}[shift={(30,0)}]
    \node at (14,2) {\small $B$};
    \draw [red!30,line width=2,->] (4,11) to (6,10) to (8,9);
    \draw [red!30,line width=2,->] (4,11) to (2,12) to (0,13);
    \draw [red!30,line width=2,->] (5.5,2) to [out=115,in=-70] (4,4.5) to [out=110,in=-80] (3.5,6);
    \draw [red!30,line width=2,->] (13.5,11) to [out=-95,in=70] (13,8) to [out=-110,in=40] (11,5);
    \draw [red!30,line width=2,->] (13.5,11) to [out=150,in=-20] (10,13) to [out=160,in=-10] (7.5,13.5);
    \draw [green!30,fill] (5.5,2) circle [radius=.6];
    \draw [blue!30,fill] (11,5) circle [radius=.6];
    \draw [green!30,fill] (13.5,11) circle [radius=.6];
    \draw [green!30,fill] (4,11) circle [radius=.6];
    \draw [blue!30,line width=3] (5.5,6.5) to [out=25,in=-110] (8,9) to [out=70,in=-75] (8,12) to [out=105,in=-60] (7.5,13.5) to [out=120,in=15] (3.5,15) to [out=-165,in=65] (0,13) to [out=-115,in=100] (0,10) to [out=-80,in=180] (3.5,6) to [out=0,in=-155] (5.5,6.5);
    \draw [ultra thick,
    ] (4.2,4.9) to [out=60,in=-130] (5.5,6.5) to [out=50,in=-75] (6,10) to [out=105,in=-30] (5,12) to [out=150,in=15] (2.3,12.1);
    \draw [black!40,ultra thick,postaction={decorate}] (1.7,11.9) to [out=-165,in=75] (0,10) to [out=-105,in=125] (1,6.5) to [out=-55,in=150] (4,4.5) to [out=-30,in=165] (6,3.5); 
    \draw [ultra thick,postaction={decorate}] (6,3.5) to [out=-15,in=135] (7.9,2.4) to [out=-45,in=30] (7.5,1) to [out=-150,in=-90] (3,2) to [out=90,in=-120] (3.8,4.1);
    \draw [black!40,ultra thick,postaction={decorate}] (9.9,13.3) to [out=105,in=-25] (8,15) to [out=155,in=35] (3.5,15) to [out=-145,in=105] (2,12) to [out=-75,in=150] (3.5,10) to [out=-30,in=-165] (5.7,9.9);
    \draw [black!40,ultra thick,postaction={decorate}] (6.3,10.1) to [out=15,in=-135] (8,12) to [out=45,in=-165] (10,13) to [out=15,in=180] (13,13.5) to [out=0,in=75] (16,11) to [out=-105,in=30] (13.4,8.2);
    \draw [black!40,ultra thick,postaction={decorate}] (12.6,7.8) to [out=-150,in=15] (11,7.5) to [out=-165,in=65] (9,5.5) to [out=-115,in=135] (9.1,3.6) to [out=-45,in=-165] (11,3)to [out=15,in=-90] (13.5,6) to [out=90,in=-60] (13,8) to [out=120,in=-65] (11,10.5) to [out=115,in=-75] (10.1,12.7);
    \draw [fill] (6,3.5) circle [radius=.3];
    \end{scope}

    \begin{scope}[>=stealth,shift={(28,-20)},scale=1.2]
    \node at (9,9) {\small $B$};
    \node [o,green] (u) at (6,6) {\color{black} $u$};
    \node [o,blue] (w) at (12,6) {\color{black} $v$};
    \node [o,green] (x) at (12,12) {\color{black} $x$};
    \node [o,blue] (y) at (6,12) {\color{black} $y$};
    \node [o,green] (z) at (3,14) {\color{black} $z$};
    \draw [red,<-,ultra thick] (y) to (z);
    \draw [red,<-,thick] (y) to [out=-48,in=-90] (9,12) to [out=90,in=50] (z);
    \draw [red,->,thick] (u) to (y);
    \draw [color=cyan,thick,dotted] (u) to (w);
    \draw [color=red!30,->] (x) to (w);
    \draw [color=red!30,->] (x) to (y);
    \node [circle,fill,scale=.3,color=gray] at (6,4.5) {};
    \draw [->,color=gray] (6,4.5) to (7,4.5);
    \node at (10.2,12.8) {\small $a$};
    \node at (12.6,9) {\small $b$};
    \node at (9,5.2) {\small $c$};
    \node at (5.2,8.8) {\small $d$};
    \node at (3.6,12.4) {\small $e$};
    \node at (7,15) {\small $f$};
    \end{scope}
    
    \begin{scope}[shift={(60,0)}]
    \node at (14,2) {\small $D$};
    \draw [red!30,line width=2,->] (4,11) to (6,10) to (8,9);
    \draw [red!30,line width=2,->] (5.5,2) to [out=115,in=-70] (4,4.5) to [out=110,in=-80] (3.5,6);
    \draw [red!30,line width=2,->] (13.5,11) to [out=-95,in=70] (13,8) to [out=-110,in=40] (11,5);
    \draw [red!30,line width=2,->] (13.5,11) to [out=150,in=-20] (10,13) to [out=160,in=-10] (7.5,13.5);
    \draw [green!30,fill] (5.5,2) circle [radius=.6];
    \draw [blue!30,fill] (11,5) circle [radius=.6];
    \draw [green!30,fill] (13.5,11) circle [radius=.6];
    \draw [green!30,fill] (4,11) circle [radius=.6];
    \draw [fill] (6,3.5) circle [radius=.3];
    \draw [blue!30,line width=3] (5.5,6.5) to [out=25,in=-110] (8,9) to [out=70,in=-75] (8,12) to [out=105,in=-60] (7.5,13.5) to [out=120,in=15] (3.5,15) to [out=-165,in=65] (0,13) to [out=-115,in=100] (0,10) to [out=-80,in=180] (3.5,6) to [out=0,in=-155] (5.5,6.5);
    \draw [ultra thick,
    ] (4.2,4.9) to [out=60,in=-130] (5.5,6.5) to [out=50,in=-75] (6,10) to [out=105,in=-30] (5,12) to [out=150,in=60] (2.8,11.6) to [out=-120,in=150] (3.5,10) to [out=-30,in=-165] (5.7,9.9);
    \draw [ultra thick,postaction={decorate}](9.9,13.3) to [out=105,in=-25] (8,15) to [out=155,in=55] (-1,14) to [out=-125,in=125] (-.5,7.5) to [out=-55,in=150] (4,4.5) to [out=-30,in=165] (6,3.5);
    \draw [ultra thick,postaction={decorate}] (6,3.5) to [out=-15,in=135] (7.9,2.4) to [out=-45,in=30] (7.5,1) to [out=-150,in=-90] (3,2) to [out=90,in=-120] (3.8,4.1);
    \draw [ultra thick,postaction={decorate}] (6.3,10.1) to [out=15,in=-135] (8,12) to [out=45,in=-165] (10,13) to [out=15,in=180] (13,13.5) to [out=0,in=75] (16,11) to [out=-105,in=30] (13.4,8.2);
    \draw [ultra thick,postaction={decorate}] (12.6,7.8) to [out=-150,in=15] (11,7.5) to [out=-165,in=65] (9,5.5) to [out=-115,in=135] (9.1,3.6) to [out=-45,in=-165] (11,3)to [out=15,in=-90] (13.5,6) to [out=90,in=-60] (13,8) to [out=120,in=-65] (11,10.5) to [out=115,in=-75] (10.1,12.7);
    \end{scope}

    \begin{scope}[>=stealth,shift={(58,-20)},scale=1.2]
    \node at (9,9) {\small $D$};
    \node [o,green] (u) at (6,6) {\color{black} $u$};
    \node [o,blue] (w) at (12,6) {\color{black} $v$};
    \node [o,green] (x) at (12,12) {\color{black} $x$};
    \node [o,blue] (y) at (6,12) {\color{black} $y$};
    \node [o,green] (z) at (3,14) {\color{black} $z$};
    \draw [red,thick,dotted] (y) to (z);
    \draw [red,<-,thick] (y) to [out=-48,in=-90] (9,12) to [out=90,in=50] (z);
    \draw [red,->,thick] (u) to (y);
    \draw [color=cyan,thick,dotted] (u) to (w);
    \draw [red,->,thick] (x) to (w);
    \draw [red,->,thick] (x) to (y);
    \node [circle,fill,scale=.3,color=gray] at (6,4.5) {};
    \draw [->,color=gray] (6,4.5) to (7,4.5);
    \node at (10.2,12.8) {\small $a$};
    \node at (12.6,9) {\small $b$};
    \node at (9,5.2) {\small $c$};
    \node at (5.2,8.8) {\small $d$};
    \node at (3.8,12.6) {\small $e$};
    \node at (7,15) {\small $f$};
    \end{scope}
    
\begin{scope}[shift={(0,-40)}]
\node at (14,2) {\small $C$};
    \draw [red!30,line width=2,->] (4,11) to (6,10) to (8,9);
    \draw [red!30,line width=2,->] (4,11) to (2,12) to (0,13);
    \draw [red!30,line width=2,->] (5.5,2) to [out=115,in=-70] (4,4.5) to [out=110,in=-80] (3.5,6);
    \draw [red!30,line width=2,->] (11,5) to [out=-140,in=25] (8.5,3) to [out=-155,in=20] (5.5,2);
    \draw [red!30,line width=2,->] (13.5,11) to [out=-95,in=70] (13,8) to [out=-110,in=40] (11,5);
    \draw [red!30,line width=2,->] (13.5,11) to [out=150,in=-20] (10,13) to [out=160,in=-10] (7.5,13.5);
    \draw [green!30,fill] (5.5,2) circle [radius=.6];
    \draw [blue!30,fill] (11,5) circle [radius=.6];
    \draw [green!30,fill] (13.5,11) circle [radius=.6];
    \draw [green!30,fill] (4,11) circle [radius=.6];
    \draw [blue!30,line width=3] (5.5,6.5) to [out=25,in=-110] (8,9) to [out=70,in=-75] (8,12) to [out=105,in=-60] (7.5,13.5) to [out=120,in=15] (3.5,15) to [out=-165,in=65] (0,13) to [out=-115,in=100] (0,10) to [out=-80,in=180] (3.5,6) to [out=0,in=-155] (5.5,6.5);
    \draw [black!40,ultra thick,
    ] (4.2,4.9) to [out=60,in=-130] (5.5,6.5) to [out=50,in=-75] (6,10) to [out=105,in=-30] (5,12) to [out=150,in=15] (2.3,12.1);
    \draw [black!40,ultra thick,postaction={decorate}] (1.7,11.9) to [out=-165,in=75] (0,10) to [out=-105,in=125] (1,6.5) to [out=-55,in=150] (4,4.5) to [out=-30,in=165] (6,3.5);
    \draw [ultra thick,postaction={decorate}] (6,3.5) to [out=-15,in=180] (8.5,3) to [out=0,in=-165] (11,3) to [out=15,in=-90] (13.5,6) to [out=90,in=-60] (13,8) to [out=120,in=-65] (11,10.5) to [out=115,in=-75] (10.1,12.7);
    \draw [black!40,ultra thick,postaction={decorate}] (9.9,13.3) to [out=105,in=-25] (8,15) to [out=155,in=35] (3.5,15) to [out=-145,in=105] (2,12) to [out=-75,in=150] (3.5,10) to [out=-30,in=-165] (5.7,9.9);
    \draw [black!40,ultra thick,postaction={decorate}] (6.3,10.1) to [out=15,in=-135] (8,12) to [out=45,in=-165] (10,13) to [out=15,in=180] (13,13.5) to [out=0,in=75] (16,11) to [out=-105,in=30] (13.4,8.2);
    \draw [black!40,ultra thick,postaction={decorate}] (12.6,7.8) to [out=-150,in=15] (11,7.5) to [out=-165,in=65] (9,5.5) to [out=-115,in=90] (8.5,3.4); 
    \draw [black!40,ultra thick,
    ] (8.5,2.6) to [out=-90,in=30] (7.5,1) to [out=-150,in=-90] (3,2) to [out=90,in=-120] (3.8,4.1);
    \draw [fill] (6,3.5) circle [radius=.3];
    \end{scope}

    \begin{scope}[>=stealth,shift={(-3,-60)},scale=1.2]
    \node at (9,9) {\small $C$};
    \node [o,green] (u) at (6,6) {\color{black} $u$};
    \node [o,blue] (w) at (12,6) {\color{black} $v$};
    \node [o,green] (x) at (12,12) {\color{black} $x$};
    \node [o,blue] (y) at (6,12) {\color{black} $y$};
    \node [o,green] (z) at (3,14) {\color{black} $z$};
    \draw [color=red!30,<-] (y) to (z);
    \draw [color=red!30,<-] (y) to [out=-48,in=-90] (9,12) to [out=90,in=50] (z);
    \draw [color=red!30,->] (u) to (y);
    \draw [color=cyan,<-,thick] (u) to (w);
    \draw [red,->,thick] (x) to (w);
    \draw [red,->,ultra thick] (x) to (y);
    \node [circle,fill,scale=.3,color=gray] at (6,4.5) {};
    \draw [->,color=gray] (6,4.5) to (7,4.5);
    \node at (10.2,12.8) {\small $a$};
    \node at (12.6,9) {\small $b$};
    \node at (9,5.2) {\small $c$};
    \node at (5.2,8.8) {\small $d$};
    \node at (3.8,12.6) {\small $e$};
    \node at (7,15) {\small $f$};
    \end{scope}

\begin{scope}[shift={(30,-40)}]
\node at (14,2) {\small $E$};
    \draw [red!30,line width=2,->] (4,11) to (6,10) to (8,9);
    \draw [red!30,line width=2,->] (4,11) to (2,12) to (0,13);
    \draw [red!30,line width=2,->] (5.5,2) to [out=115,in=-70] (4,4.5) to [out=110,in=-80] (3.5,6);
    \draw [red!30,line width=2,->] (11,5) to [out=-140,in=25] (8.5,3) to [out=-155,in=20] (5.5,2);
    \draw [red!30,line width=2,->] (13.5,11) to [out=-95,in=70] (13,8) to [out=-110,in=40] (11,5);
    \draw [green!30,fill] (5.5,2) circle [radius=.6];
    \draw [blue!30,fill] (11,5) circle [radius=.6];
    \draw [green!30,fill] (13.5,11) circle [radius=.6];
    \draw [green!30,fill] (4,11) circle [radius=.6];
    \draw [blue!30,line width=3] (5.5,6.5) to [out=25,in=-110] (8,9) to [out=70,in=-75] (8,12) to [out=105,in=-60] (7.5,13.5) to [out=120,in=15] (3.5,15) to [out=-165,in=65] (0,13) to [out=-115,in=100] (0,10) to [out=-80,in=180] (3.5,6) to [out=0,in=-155] (5.5,6.5);
    \draw [ultra thick,
    ] (4.2,4.9) to [out=60,in=-130] (5.5,6.5) to [out=50,in=-75] (6,10) to [out=105,in=-30] (5,12) to [out=150,in=15] (2.3,12.1);
    \draw [black!40,ultra thick,postaction={decorate}] (1.7,11.9) to [out=-165,in=75] (0,10) to [out=-105,in=125] (1,6.5) to [out=-55,in=150] (4,4.5) to [out=-30,in=165] (6,3.5);
    \draw [ultra thick,postaction={decorate}] (6,3.5) to [out=-15,in=180] (8.5,3) to [out=0,in=-165] (11,3) to [out=15,in=-90] (13.5,6) to [out=90,in=-60] (13,8) to [out=120,in=-65] (11,10.5) to [out=115,in=-120] (10.8,12.6) to [out=60,in=180] (13,13.5) to [out=0,in=75] (16,11) to [out=-105,in=30] (13.4,8.2);
    \draw [black!40,ultra thick,postaction={decorate}] (6.3,10.1) to [out=15,
    in=-20] (5.5,13.5) to [out=160,in=105] (2,12) to [out=-75,in=150] (3.5,10) to [out=-30,in=-165] (5.7,9.9);
    \draw [ultra thick,postaction={decorate}] (12.6,7.8) to [out=-150,in=15] (11,7.5) to [out=-165,in=65] (9,5.5) to [out=-115,in=90] (8.5,3.4); 
    \draw [ultra thick,postaction={decorate}
    ] (8.5,2.6) to [out=-90,in=30] (7.5,1) to [out=-150,in=-90] (3,2) to [out=90,in=-120] (3.8,4.1);
    \draw [fill] (6,3.5) circle [radius=.3];
    \end{scope}

    \begin{scope}[>=stealth,shift={(27,-60)},scale=1.2]
    \node at (9,9) {\small $E$};
    \node [o,green] (u) at (6,6) {\color{black} $u$};
    \node [o,blue] (w) at (12,6) {\color{black} $v$};
    \node [o,green] (x) at (12,12) {\color{black} $x$};
    \node [o,blue] (y) at (6,12) {\color{black} $y$};
    \node [o,green] (z) at (3,14) {\color{black} $z$};
    \draw [red,<-,ultra thick] (y) to (z);
    \draw [red,<-,thick] (y) to [out=-48,in=-90] (9,12) to [out=90,in=50] (z);
    \draw [red,->,thick] (u) to (y);
    \draw [color=cyan,<-,thick] (u) to (w);
    \draw [red,->,thick] (x) to (w);
    \draw [color=cyan,thick,dotted] (x) to (y);
    \node [circle,fill,scale=.3,color=gray] at (6,4.5) {};
    \draw [->,color=gray] (6,4.5) to (7,4.5);
    \node at (10.2,12.8) {\small $a$};
    \node at (12.6,9) {\small $b$};
    \node at (9,5.2) {\small $c$};
    \node at (5.2,8.8) {\small $d$};
    \node at (3.6,12.4) {\small $e$};
    \node at (7,15) {\small $f$};
    \end{scope}

\begin{scope}[shift={(60,-40)}]
\node at (14,2) {\small $F$};
    \draw [red!30,line width=2,->] (4,11) to (6,10) to (8,9);
    \draw [red!30,line width=2,->] (5.5,2) to [out=115,in=-70] (4,4.5) to [out=110,in=-80] (3.5,6);
    \draw [red!30,line width=2,->] (11,5) to [out=-140,in=25] (8.5,3) to [out=-155,in=20] (5.5,2);
    \draw [red!30,line width=2,->] (13.5,11) to [out=-95,in=70] (13,8) to [out=-110,in=40] (11,5);
    \draw [green!30,fill] (5.5,2) circle [radius=.6];
    \draw [blue!30,fill] (11,5) circle [radius=.6];
    \draw [green!30,fill] (13.5,11) circle [radius=.6];
    \draw [green!30,fill] (4,11) circle [radius=.6];
    \draw [fill] (6,3.5) circle [radius=.3];
    \draw [blue!30,line width=3] (5.5,6.5) to [out=25,in=-110] (8,9) to [out=70,in=-75] (8,12) to [out=105,in=-60] (7.5,13.5) to [out=120,in=15] (3.5,15) to [out=-165,in=65] (0,13) to [out=-115,in=100] (0,10) to [out=-80,in=180] (3.5,6) to [out=0,in=-155] (5.5,6.5);
    \draw [ultra thick,postaction={decorate}
    ](6.3,10.1) to [out=15,in=-135] (8,12) to [out=45,in=-40] (7.5,15) to [out=140,in=55] (-1,14) to [out=-125,in=125] (-.5,7.5) to [out=-55,in=150] (4,4.5) to [out=-30,in=165] (6,3.5);
    \draw [ultra thick,postaction={decorate}] (6,3.5) to [out=-15,in=180] (8.5,3) to [out=0,in=-165] (11,3) to [out=15,in=-90] (13.5,6) to [out=90,in=-60] (13,8) to [out=120,in=-65] (11,10.5) to [out=115,in=-120] (10.8,12.6) to [out=60,in=180] (13,13.5) to [out=0,in=75] (16,11) to [out=-105,in=30] (13.4,8.2);
\draw [ultra thick,
    ] (4.2,4.9) to [out=60,in=-130] (5.5,6.5) to [out=50,in=-75] (6,10) to [out=105,in=-30] (5,12) to [out=150,in=60] (2.8,11.6) to [out=-120,in=150] (3.5,10) to [out=-30,in=-165] (5.7,9.9);
     \draw [ultra thick,postaction={decorate}] (12.6,7.8) to [out=-150,in=15] (11,7.5) to [out=-165,in=65] (9,5.5) to [out=-115,in=90] (8.5,3.4); 
    \draw [ultra thick,
    ] (8.5,2.6) to [out=-90,in=30] (7.5,1) to [out=-150,in=-90] (3,2) to [out=90,in=-120] (3.8,4.1);
    \end{scope}
    
    \begin{scope}[>=stealth,shift={(57,-60)},scale=1.2]
    \node at (9,9) {\small $F$};
    \node [o,green] (u) at (6,6) {\color{black} $u$};
    \node [o,blue] (w) at (12,6) {\color{black} $v$};
    \node [o,green] (x) at (12,12) {\color{black} $x$};
    \node [o,blue] (y) at (6,12) {\color{black} $y$};
    \node [o,green] (z) at (3,14) {\color{black} $z$};
    \draw [red,thick,dotted] (y) to (z);
    \draw [red,<-,thick] (y) to [out=-48,in=-90] (9,12) to [out=90,in=50] (z);
    \draw [red,->,thick] (u) to (y);
    \draw [color=cyan,<-,thick] (u) to (w);
    \draw [red,->,thick] (x) to (w);
    \draw [color=cyan,thick,dotted] (x) to (y);
    \node [circle,fill,scale=.3,color=gray] at (6,4.5) {};
    \draw [->,color=gray] (6,4.5) to (7,4.5);
    \node at (10.2,12.8) {\small $a$};
    \node at (12.6,9) {\small $b$};
    \node at (9,5.2) {\small $c$};
    \node at (5.2,8.8) {\small $d$};
    \node at (3.8,12.6) {\small $e$};
    \node at (7,15) {\small $f$};
    \end{scope}

\begin{scope}[shift={(0,-80)}]
\node at (14,2) {\small $G$};
    \draw [red!30,line width=2,->] (4,11) to (6,10) to (8,9);
    \draw [red!30,line width=2,->] (4,11) to (2,12) to (0,13);
    \draw [red!30,line width=2,->] (5.5,2) to [out=115,in=-70] (4,4.5) to [out=110,in=-80] (3.5,6);
    \draw [red!30,line width=2,->] (11,5) to [out=-140,in=25] (8.5,3) to [out=-155,in=20] (5.5,2);
    \draw [red!30,line width=2,->] (13.5,11) to [out=-95,in=70] (13,8) to [out=-110,in=40] (11,5);
    \draw [red!30,line width=2,->] (7.5,13.5) to [out=-10,in=160] (10,13) to [out=-20,in=150] (13.5,11);
    \draw [green!30,fill] (5.5,2) circle [radius=.6];
    \draw [blue!30,fill] (11,5) circle [radius=.6];
    \draw [green!30,fill] (13.5,11) circle [radius=.6];
    \draw [green!30,fill] (4,11) circle [radius=.6];
    \draw [fill] (6,3.5) circle [radius=.3];
    \draw [blue!30,line width=3] (5.5,6.5) to [out=25,in=-110] (8,9) to [out=70,in=-75] (8,12) to [out=105,in=-60] (7.5,13.5) to [out=120,in=15] (3.5,15) to [out=-165,in=65] (0,13) to [out=-115,in=100] (0,10) to [out=-80,in=180] (3.5,6) to [out=0,in=-155] (5.5,6.5);
    \draw [ultra thick,
    ] (4.2,4.9) to [out=60,in=-130] (5.5,6.5) to [out=50,in=-75] (6,10) to [out=105,in=-30] (5,12) to [out=150,in=15] (2.3,12.1);
    \draw [ultra thick,postaction={decorate}] (1.7,11.9) to [out=-165,in=75] (0,10) to [out=-105,in=125] (1,6.5) to [out=-55,in=150] (4,4.5) to [out=-30,in=165] (6,3.5) to [out=-15,in=180] (8.5,3) to [out=0,in=-165] (11,3) to [out=15,in=-90] (13.5,6) to [out=90,in=-60] (13,8) to [out=120,in=-65] (11,10.5) to [out=115,in=-75] (10,13) to [out=105,in=-25] (8,15) to [out=155,in=35] (3.5,15) to [out=-145,in=105] (2,12) to [out=-75,in=150] (3.5,10) to [out=-30,in=-165] (5.7,9.9);
    \draw [ultra thick,postaction={decorate}] (6.3,10.1) to [out=15,in=-135] (8,12) to [out=45,in=-165] (9.7,12.9);
     \draw [ultra thick,postaction={decorate}] (10.3,13.1) to [out=15,in=180] (13,13.5) to [out=0,in=75] (16,11) to [out=-105,in=30] (13.4,8.2);
    \draw [ultra thick,postaction={decorate}] (12.6,7.8) to [out=-150,in=15] (11,7.5) to [out=-165,in=65] (9,5.5) to [out=-115,in=90] (8.5,3.4); 
    \draw [ultra thick,
    ] (8.5,2.6) to [out=-90,in=30] (7.5,1) to [out=-150,in=-90] (3,2) to [out=90,in=-120] (3.8,4.1);
    \end{scope}

    \begin{scope}[>=stealth,shift={(-3,-100)},scale=1.2]
    \node at (9,9) {\small $G$};
    \node [o,green] (u) at (6,6) {\color{black} $u$};
    \node [o,blue] (w) at (12,6) {\color{black} $v$};
    \node [o,green] (x) at (12,12) {\color{black} $x$};
    \node [o,blue] (y) at (6,12) {\color{black} $y$};
    \node [o,green] (z) at (3,14) {\color{black} $z$};
    \draw [red,<-,thick] (y) to (z);
    \draw [red,<-,thick] (y) to [out=-48,in=-90] (9,12) to [out=90,in=50] (z);
    \draw [red,->,thick] (u) to (y);
    \draw [color=cyan,<-,thick] (u) to (w);
    \draw [red,->,thick] (x) to (w);
    \draw [color=cyan,<-,thick] (x) to (y);
    \node [circle,fill,scale=.3,color=gray] at (6,4.5) {};
    \draw [->,color=gray] (6,4.5) to (7,4.5);
    \node at (10.2,12.8) {\small $a$};
    \node at (12.6,9) {\small $b$};
    \node at (9,5.2) {\small $c$};
    \node at (5.2,8.8) {\small $d$};
    \node at (3.8,12.6) {\small $e$};
    \node at (7,15) {\small $f$};
    \end{scope}
\end{tikzpicture}
    \caption{A Homfly computation tree for the standard diagram of $6_1$.}
    \label{fig:computationtree}
\end{figure}
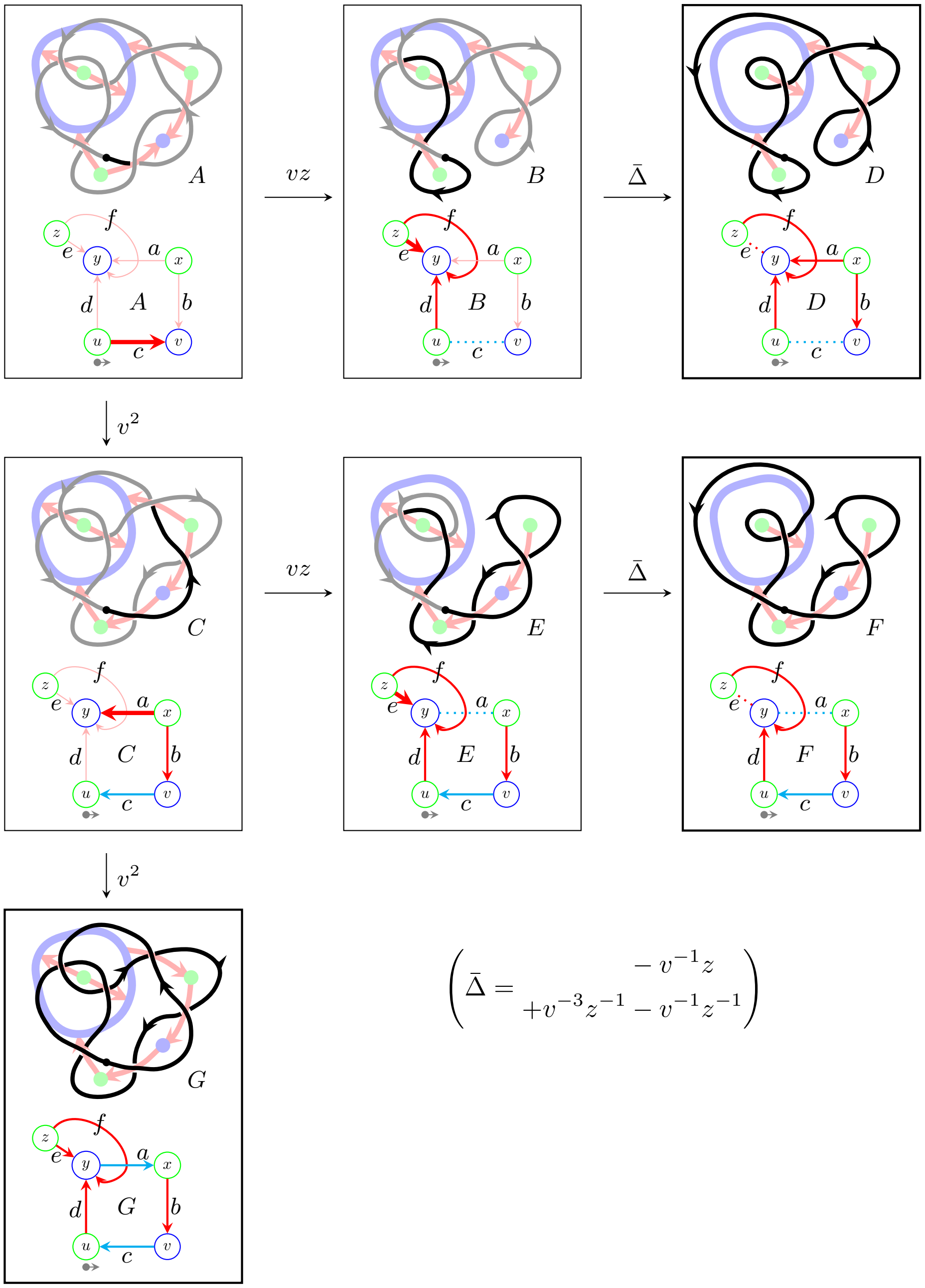

We give a step-by-step breakdown of the computation for the standard diagram of the $6_1$ knot.
It is summarized in Figure \ref{fig:computationtree}. We 
denote the initial
Seifert graph with $A$ and draw it both together with the diagram and separately beneath it. The basepoint 
is as specified in subsection \ref{ss:first_ex}.

The first crossing after the basepoint, or rather the corresponding edge $c$ (in the notation of Figure \ref{fig:6_1_link_diagram_and_Seifert}), can both be included in and excluded from a Jaeger quasi-tree of $A$. 
This can be checked directly but it will also be verified a posteriori by the rest of the computation.
In any case we apply the skein relation \eqref{eq:positive_skein} here, resulting in diagrams with Seifert graphs $B$ and $C$, respectively. 
Continuing with $B$, the next edge is $d$, whose removal would result in a disconnected graph. Hence $d$ cannot be deleted and we simply ignore it, moving on to
the next edge $f$. It is reached at its head, whence it cannot be removed, either. The next edge $e$, on the other hand, needs to be removed, since that is our last chance to connect the two boundary components of $B$. Thus we pass to the diagram $D$.

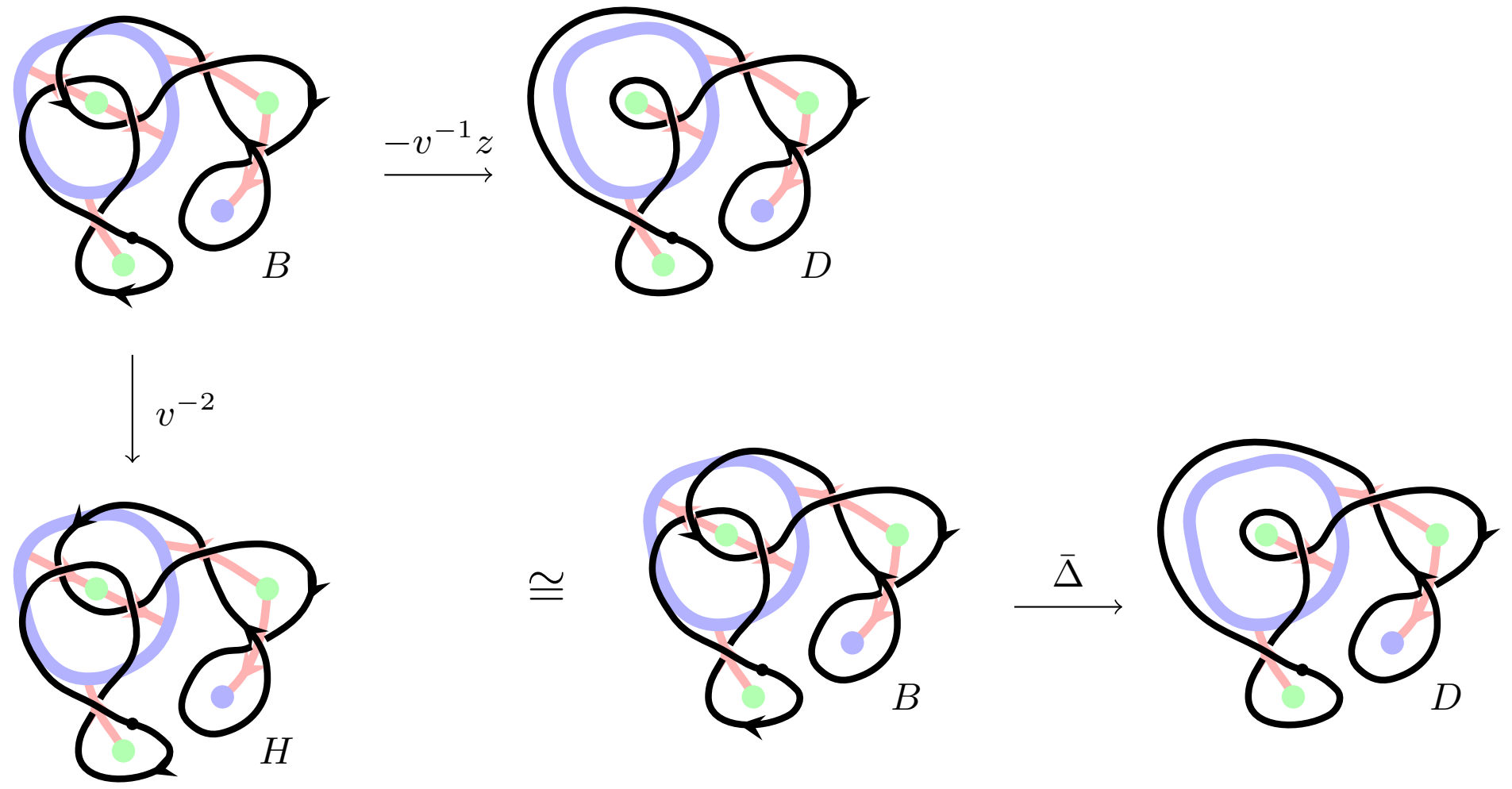
\begin{figure}[h]
    \centering
    \begin{tikzpicture}[scale=.15,decoration={markings,
    mark= at position 0.8 with {\arrow{stealth}}}]

    \draw [->] (20,7) -- (26,7);
    \node at (23,9) {\small $-v^{-1}z$};
    \draw [->] (55,-17) -- (61,-17);
    \node at (58,-15) {\small $\bar\Delta$};
    \draw [->] (6,-3) -- (6,-9);
    \node at (9,-6) {\small $v^{-2}$};
    \node at (29,-16) {\Large $\cong$};
    
    \begin{scope}[shift={(0,0)}]
    \node at (14,2) {\small $B$};
    \draw [red!30,line width=2,postaction={decorate}] (4,11) to (6,10) to (8,9);
    \draw [red!30,line width=2,postaction={decorate}] (4,11) to (2,12) to (0,13);
    \draw [red!30,line width=2,postaction={decorate}] (5.5,2) to [out=115,in=-70] (4,4.5) to [out=110,in=-80] (3.5,6);
    \draw [red!30,line width=2,postaction={decorate}] (13.5,11) to [out=-95,in=70] (13,8) to [out=-110,in=40] (11,5);
    \draw [red!30,line width=2,postaction={decorate}] (13.5,11) to [out=150,in=-20] (10,13) to [out=160,in=-10] (7.5,13.5);
    \draw [green!30,fill] (5.5,2) circle [radius=.6];
    \draw [blue!30,fill] (11,5) circle [radius=.6];
    \draw [green!30,fill] (13.5,11) circle [radius=.6];
    \draw [green!30,fill] (4,11) circle [radius=.6];
    \draw [fill] (6,3.5) circle [radius=.3];
    \draw [blue!30,line width=3] (5.5,6.5) to [out=25,in=-110] (8,9) to [out=70,in=-75] (8,12) to [out=105,in=-60] (7.5,13.5) to [out=120,in=15] (3.5,15) to [out=-165,in=65] (0,13) to [out=-115,in=100] (0,10) to [out=-80,in=180] (3.5,6) to [out=0,in=-155] (5.5,6.5);
    \draw [ultra thick,
    ] (4.2,4.9) to [out=60,in=-130] (5.5,6.5) to [out=50,in=-75] (6,10) to [out=105,in=-30] (5,12) to [out=150,in=15] (2.3,12.1);
    \draw [ultra thick,postaction={decorate}] (1.7,11.9) to [out=-165,in=75] (0,10) to [out=-105,in=125] (1,6.5) to [out=-55,in=150] (4,4.5) to [out=-30,in=165] (6,3.5) to [out=-15,in=135] (7.9,2.4) to [out=-45,in=30] (7.5,1) to [out=-150,in=-90] (3,2) to [out=90,in=-120] (3.8,4.1);
    \draw [ultra thick,postaction={decorate}] (9.9,13.3) to [out=105,in=-25] (8,15) to [out=155,in=35] (3.5,15) to [out=-145,in=105] (2,12) to [out=-75,in=150] (3.5,10) to [out=-30,in=-165] (5.7,9.9);
    \draw [ultra thick,postaction={decorate}] (6.3,10.1) to [out=15,in=-135] (8,12) to [out=45,in=-165] (10,13) to [out=15,in=180] (13,13.5) to [out=0,in=75] (16,11) to [out=-105,in=30] (13.4,8.2);
    \draw [ultra thick,postaction={decorate}] (12.6,7.8) to [out=-150,in=15] (11,7.5) to [out=-165,in=65] (9,5.5) to [out=-115,in=135] (9.1,3.6) to [out=-45,in=-165] (11,3)to [out=15,in=-90] (13.5,6) to [out=90,in=-60] (13,8) to [out=120,in=-65] (11,10.5) to [out=115,in=-75] (10.1,12.7);
    \end{scope}
    
    \begin{scope}[shift={(30,0)}]
    \node at (14,2) {\small $D$};
    \draw [red!30,line width=2,postaction={decorate}] (4,11) to (6,10) to (8,9);
    \draw [red!30,line width=2,postaction={decorate}] (5.5,2) to [out=115,in=-70] (4,4.5) to [out=110,in=-80] (3.5,6);
    \draw [red!30,line width=2,postaction={decorate}] (13.5,11) to [out=-95,in=70] (13,8) to [out=-110,in=40] (11,5);
    \draw [red!30,line width=2,postaction={decorate}] (13.5,11) to [out=150,in=-20] (10,13) to [out=160,in=-10] (7.5,13.5);
    \draw [green!30,fill] (5.5,2) circle [radius=.6];
    \draw [blue!30,fill] (11,5) circle [radius=.6];
    \draw [green!30,fill] (13.5,11) circle [radius=.6];
    \draw [green!30,fill] (4,11) circle [radius=.6];
    \draw [fill] (6,3.5) circle [radius=.3];
    \draw [blue!30,line width=3] (5.5,6.5) to [out=25,in=-110] (8,9) to [out=70,in=-75] (8,12) to [out=105,in=-60] (7.5,13.5) to [out=120,in=15] (3.5,15) to [out=-165,in=65] (0,13) to [out=-115,in=100] (0,10) to [out=-80,in=180] (3.5,6) to [out=0,in=-155] (5.5,6.5);
    \draw [ultra thick,
    ] (4.2,4.9) to [out=60,in=-130] (5.5,6.5) to [out=50,in=-75] (6,10) to [out=105,in=-30] (5,12) to [out=150,in=60] (2.8,11.6) to [out=-120,in=150] (3.5,10) to [out=-30,in=-165] (5.7,9.9);
        \draw [ultra thick,
](9.9,13.3) to [out=105,in=-25] (8,15) to [out=155,in=55] (-1,14) to [out=-125,in=125] (-.5,7.5) to [out=-55,in=150] (4,4.5) to [out=-30,in=165] (6,3.5) to [out=-15,in=135] (7.9,2.4) to [out=-45,in=30] (7.5,1) to [out=-150,in=-90] (3,2) to [out=90,in=-120] (3.8,4.1);
    \draw [ultra thick,postaction={decorate}] (6.3,10.1) to [out=15,in=-135] (8,12) to [out=45,in=-165] (10,13) to [out=15,in=180] (13,13.5) to [out=0,in=75] (16,11) to [out=-105,in=30] (13.4,8.2);
    \draw [ultra thick,postaction={decorate}] (12.6,7.8) to [out=-150,in=15] (11,7.5) to [out=-165,in=65] (9,5.5) to [out=-115,in=135] (9.1,3.6) to [out=-45,in=-165] (11,3)to [out=15,in=-90] (13.5,6) to [out=90,in=-60] (13,8) to [out=120,in=-65] (11,10.5) to [out=115,in=-75] (10.1,12.7);
    \end{scope}
    
    \begin{scope}[shift={(0,-27)}]
    \node at (14,2) {\small $H$};
    \draw [red!30,line width=2,postaction={decorate}] (4,11) to (6,10) to (8,9);
    \draw [red!30,line width=2,postaction={decorate}] (0,13) to (2,12) to (4,11);
    \draw [red!30,line width=2,postaction={decorate}] (5.5,2) to [out=115,in=-70] (4,4.5) to [out=110,in=-80] (3.5,6);
    \draw [red!30,line width=2,postaction={decorate}] (13.5,11) to [out=-95,in=70] (13,8) to [out=-110,in=40] (11,5);
    \draw [red!30,line width=2,postaction={decorate}] (13.5,11) to [out=150,in=-20] (10,13) to [out=160,in=-10] (7.5,13.5);
    \draw [green!30,fill] (5.5,2) circle [radius=.6];
    \draw [blue!30,fill] (11,5) circle [radius=.6];
    \draw [green!30,fill] (13.5,11) circle [radius=.6];
    \draw [green!30,fill] (4,11) circle [radius=.6];
    \draw [fill] (6,3.5) circle [radius=.3];
    \draw [blue!30,line width=3] (5.5,6.5) to [out=25,in=-110] (8,9) to [out=70,in=-75] (8,12) to [out=105,in=-60] (7.5,13.5) to [out=120,in=15] (3.5,15) to [out=-165,in=65] (0,13) to [out=-115,in=100] (0,10) to [out=-80,in=180] (3.5,6) to [out=0,in=-155] (5.5,6.5);
    \draw [ultra thick,postaction={decorate}
    ] (4.2,4.9) to [out=60,in=-130] (5.5,6.5) to [out=50,in=-75] (6,10) to [out=105,in=-30] (5,12) to [out=150,in=15] (2,12) to [out=-165,in=75] (0,10) to [out=-105,in=125] (1,6.5) to [out=-55,in=150] (4,4.5) to [out=-30,in=165] (6,3.5) to [out=-15,in=135] (7.9,2.4) to [out=-45,in=30] (7.5,1) to [out=-150,in=-90] (3,2) to [out=90,in=-120] (3.8,4.1);
    \draw [ultra thick,postaction={decorate}] (9.9,13.3) to [out=105,in=-25] (8,15) to [out=155,in=35] (3.5,15) to [out=-145,in=105] (1.9,12.3);
    \draw [ultra thick,
    ](2.1,11.7) to [out=-75,in=150] (3.5,10) to [out=-30,in=-165] (5.7,9.9);
    \draw [ultra thick,postaction={decorate}] (6.3,10.1) to [out=15,in=-135] (8,12) to [out=45,in=-165] (10,13) to [out=15,in=180] (13,13.5) to [out=0,in=75] (16,11) to [out=-105,in=30] (13.4,8.2);
    \draw [ultra thick,postaction={decorate}] (12.6,7.8) to [out=-150,in=15] (11,7.5) to [out=-165,in=65] (9,5.5) to [out=-115,in=135] (9.1,3.6) to [out=-45,in=-165] (11,3)to [out=15,in=-90] (13.5,6) to [out=90,in=-60] (13,8) to [out=120,in=-65] (11,10.5) to [out=115,in=-75] (10.1,12.7);
    \end{scope}

    \begin{scope}[shift={(35,-24)}]
    \node at (14,2) {\small $B$};
    \draw [red!30,line width=2,postaction={decorate}] (4,11) to (6,10) to (8,9);
    \draw [red!30,line width=2,postaction={decorate}] (4,11) to (2,12) to (0,13);
    \draw [red!30,line width=2,postaction={decorate}] (5.5,2) to [out=115,in=-70] (4,4.5) to [out=110,in=-80] (3.5,6);
    \draw [red!30,line width=2,postaction={decorate}] (13.5,11) to [out=-95,in=70] (13,8) to [out=-110,in=40] (11,5);
    \draw [red!30,line width=2,postaction={decorate}] (13.5,11) to [out=150,in=-20] (10,13) to [out=160,in=-10] (7.5,13.5);
    \draw [green!30,fill] (5.5,2) circle [radius=.6];
    \draw [blue!30,fill] (11,5) circle [radius=.6];
    \draw [green!30,fill] (13.5,11) circle [radius=.6];
    \draw [green!30,fill] (4,11) circle [radius=.6];
    \draw [fill] (6,3.5) circle [radius=.3];
    \draw [blue!30,line width=3] (5.5,6.5) to [out=25,in=-110] (8,9) to [out=70,in=-75] (8,12) to [out=105,in=-60] (7.5,13.5) to [out=120,in=15] (3.5,15) to [out=-165,in=65] (0,13) to [out=-115,in=100] (0,10) to [out=-80,in=180] (3.5,6) to [out=0,in=-155] (5.5,6.5);
    \draw [ultra thick,
    ] (4.2,4.9) to [out=60,in=-130] (5.5,6.5) to [out=50,in=-75] (6,10) to [out=105,in=-30] (5,12) to [out=150,in=15] (2.3,12.1);
    \draw [ultra thick,postaction={decorate}] (1.7,11.9) to [out=-165,in=75] (0,10) to [out=-105,in=125] (1,6.5) to [out=-55,in=150] (4,4.5) to [out=-30,in=165] (6,3.5) to [out=-15,in=135] (7.9,2.4) to [out=-45,in=30] (7.5,1) to [out=-150,in=-90] (3,2) to [out=90,in=-120] (3.8,4.1);
    \draw [ultra thick,postaction={decorate}] (9.9,13.3) to [out=105,in=-25] (8,15) to [out=155,in=35] (3.5,15) to [out=-145,in=105] (2,12) to [out=-75,in=150] (3.5,10) to [out=-30,in=-165] (5.7,9.9);
    \draw [ultra thick,postaction={decorate}] (6.3,10.1) to [out=15,in=-135] (8,12) to [out=45,in=-165] (10,13) to [out=15,in=180] (13,13.5) to [out=0,in=75] (16,11) to [out=-105,in=30] (13.4,8.2);
    \draw [ultra thick,postaction={decorate}] (12.6,7.8) to [out=-150,in=15] (11,7.5) to [out=-165,in=65] (9,5.5) to [out=-115,in=135] (9.1,3.6) to [out=-45,in=-165] (11,3)to [out=15,in=-90] (13.5,6) to [out=90,in=-60] (13,8) to [out=120,in=-65] (11,10.5) to [out=115,in=-75] (10.1,12.7);
    \end{scope}
    
    \begin{scope}[shift={(65,-24)}]
    \node at (14,2) {\small $D$};
    \draw [red!30,line width=2,postaction={decorate}] (4,11) to (6,10) to (8,9);
    \draw [red!30,line width=2,postaction={decorate}] (5.5,2) to [out=115,in=-70] (4,4.5) to [out=110,in=-80] (3.5,6);
    \draw [red!30,line width=2,postaction={decorate}] (13.5,11) to [out=-95,in=70] (13,8) to [out=-110,in=40] (11,5);
    \draw [red!30,line width=2,postaction={decorate}] (13.5,11) to [out=150,in=-20] (10,13) to [out=160,in=-10] (7.5,13.5);
    \draw [green!30,fill] (5.5,2) circle [radius=.6];
    \draw [blue!30,fill] (11,5) circle [radius=.6];
    \draw [green!30,fill] (13.5,11) circle [radius=.6];
    \draw [green!30,fill] (4,11) circle [radius=.6];
    \draw [fill] (6,3.5) circle [radius=.3];
    \draw [blue!30,line width=3] (5.5,6.5) to [out=25,in=-110] (8,9) to [out=70,in=-75] (8,12) to [out=105,in=-60] (7.5,13.5) to [out=120,in=15] (3.5,15) to [out=-165,in=65] (0,13) to [out=-115,in=100] (0,10) to [out=-80,in=180] (3.5,6) to [out=0,in=-155] (5.5,6.5);
    \draw [ultra thick,
    ] (4.2,4.9) to [out=60,in=-130] (5.5,6.5) to [out=50,in=-75] (6,10) to [out=105,in=-30] (5,12) to [out=150,in=60] (2.8,11.6) to [out=-120,in=150] (3.5,10) to [out=-30,in=-165] (5.7,9.9);
    \draw [ultra thick,
    ](9.9,13.3) to [out=105,in=-25] (8,15) to [out=155,in=55] (-1,14) to [out=-125,in=125] (-.5,7.5) to [out=-55,in=150] (4,4.5) to [out=-30,in=165] (6,3.5) to [out=-15,in=135] (7.9,2.4) to [out=-45,in=30] (7.5,1) to [out=-150,in=-90] (3,2) to [out=90,in=-120] (3.8,4.1);
    \draw [ultra thick,postaction={decorate}] (6.3,10.1) to [out=15,in=-135] (8,12) to [out=45,in=-165] (10,13) to [out=15,in=180] (13,13.5) to [out=0,in=75] (16,11) to [out=-105,in=30] (13.4,8.2);
    \draw [ultra thick,postaction={decorate}] (12.6,7.8) to [out=-150,in=15] (11,7.5) to [out=-165,in=65] (9,5.5) to [out=-115,in=135] (9.1,3.6) to [out=-45,in=-165] (11,3)to [out=15,in=-90] (13.5,6) to [out=90,in=-60] (13,8) to [out=120,in=-65] (11,10.5) to [out=115,in=-75] (10.1,12.7);
    \end{scope}

\end{tikzpicture}
    \caption{Simplification in the Homfly computation process that 
    applies at crossings that need to be smoothed.}
    \label{fig:shorthand}
\end{figure}

Let us take a moment to explain why the multiplier $\bar\Delta$ appears between $B$ and $D$.
Please refer to Figure \ref{fig:shorthand}. If we ``properly'' applied the skein relation \eqref{eq:negative_skein} to the edge (i.e., crossing) $e$ of $B$, we would obtain two new graphs (diagrams) $H$ and $D$. Here $H$ depicts a distant union of two links (in this case, two unknots) and $D$ shows the connected sum of the same two links. Thus, as is well known, their Homfly polynomials differ by a factor of $\frac{v^{-1}-v}z$ and we do see that, in this case,
\begin{multline*}
P_{\includegraphics[totalheight=8pt]{skeinneg.eps}}=v^{-2}\,P_{\includegraphics[totalheight=8pt]{skeinpos.eps}}\,-\,v^{-1}z\,P_{\includegraphics[totalheight=8pt]{skeinsmooth.eps}}
=v^{-2}(v^{-1}z^{-1}-vz^{-1})\,P_{\includegraphics[totalheight=8pt]{skeinsmooth.eps}}\,-\,v^{-1}z\,P_{\includegraphics[totalheight=8pt]{skeinsmooth.eps}} \\
=(-v^{-1}z+v^{-3}z^{-1}-v^{-1}z^{-1})\,P_{\includegraphics[totalheight=8pt]{skeinsmooth.eps}}.
\end{multline*}
In Lemma \ref{l:handling_edges_that_need_to_be_removed} we will see that a similar phenomenon plays out in general, for all externally active edges of both signs.

As to the rest of the graph/diagram
$D$, one may realize that each edge/crossing has either been addressed already, or it cannot be removed from the graph. That, in turn, is either because the edge is reached at its head, or because removing it would make the graph disconnected. Hence this branch of the computation tree ends in $D$.
Notice that the graph $D$ 
is identical 
to $K_1$ of Figure \ref{fig:61_Jaeger_kvazifai}---we did not even have to reverse any internally passive edges because there were not any.
Also, note that $D$ is a diagram of the unknot.

Stepping back to the diagram $C$, after reversing $c$, the process must also include $b$ as it is reached at its head. The edge $a$ does offer both options, so we perform the skein relation \eqref{eq:positive_skein} at the corresponding crossing. After smoothing (diagram $E$), the walk visits the edges $b$ and $c$ for the second time, then $d$ for the first time but $d$ must be kept, as does $f$, whereas $e$ must be cut for the same reason as in the case of $B$. This leads us to the end stage $F$, which corresponds to the Jaeger quasi-tree $K_2$ with the internally passive edge $c$ reversed. Notice that $F$ also has a trivially embedded boundary. 

Finally, let us examine the diagram $G$. If we continue the walk thus far, the next edge is $e$ and we reach it at its head. The parallel edge $f$ is reached at its tail but we cannot remove it. Indeed, by removing $f$, we would get to $e$ again, then to $d$, at its head, and then we would arrive back to the basepoint without visiting all of the boundary. 
Thus indeed $f$ must be kept and traversed, after which we revisit $a$, $b$, and $c$ and come to $d$ for the first time. Since cutting it would finish the walk too early, we 
keep $d$. The walk coninues with the edges $f$, $e$, and $d$, each of which has been addressed before, and then 
the process stops. The graph $G$ corresponds to the Jaeger quasi-tree $K_3$, via the reversal of the internally passive edges $a$ and $c$. 

Importantly, $G$ also represents
a trivial knot. To see this, note that if we did smooth the crossing corresponding to $f$ (which is the only crossing at which the diagram is not descending) then we would have obtained two unknots that lie over one another; $G$ itself is the connected sum of these two unknots.

This phenomenon, namely that meeting a certain \emph{combinatorial} stopping condition is sufficient to guarantee the \emph{topological} triviality of each end stage, is absolutely key to this paper. Our process will generate these unknot diagrams typically without them being descending. (They will be connected sums of several descending diagrams.) Much of Section \ref{s:formula_proof} is devoted to hashing out the details. 

Since all three leaves of the computation tree are indeed diagrams of the unknot, it yields the already stated result
\begin{multline*}\label{eq:6_1=homfly}
P_{6_1}(v,z)=vz\bar\Delta+v^3z\bar\Delta+v^4=(1+v^2)\left(\begin{aligned}
    &-z^2\\+v^{-2}&-1
\end{aligned}\right)+v^4 \\
=\begin{array}{cccc}
    &-z^2&-v^2z^2&\\
    +v^{-2}&&-v^2&+v^4
\end{array}.\end{multline*}

\section{Preliminaries}\label{s:prelim}

After the informal introduction, we give the precise definitions.

\subsection{Links and graphs}\label{ss:links_and_graphs}

For $n$ a natural number, an \emph{$n$-component link} is a smooth embedding, into the standard Euclidean $3$-space $\R^3$, of a disjoint union of $n$ copies of the standard unit circle $S^1$. Note that the embedding induces an \emph{orientation} of its image. A \emph{knot} is a one-component link. In this paper knots and links are considered up to orientation-preserving ambient isotopy, that is, we work with the usual notion of an oriented link. 

A \emph{link diagram} is a generic orthogonal projection of a link to a plane. The projection fails to be an embedding only at finitely many transverse double points called \emph{crossings}, at which we retain over- and undercrossing information. A crossing is \emph{positive} if the undercrossing strand passes from the right to the left of the other strand, and \emph{negative} otherwise; see Figure \ref{fig:posnegsmooth}. Two link diagrams represent the same oriented link if and only if they can be connected with a finite sequence of the well known Reidemeister moves.

The \emph{Seifert circles} of a link diagram are obtained when each crossing is smoothed in the unique way that respects the orientation. They form the vertex set of the \emph{Seifert graph}, in which two Seifert circles are connected with one edge for each crossing that they share.  
Notice that the two strands of any crossing, when we think of them as connections between the two incident Seifert circles, are directed oppositely. We orient each edge of the Seifert graph as induced by the direction of the strand that crosses beneath the other strand. (See Figure \ref{fig:posnegsmooth}.) 
For an edge $e$ of the digraph $G$, we denote by $h_G(e)$ the head, and by $t_G(e)$ the tail of $e$. (We may omit the subscript if $G$ is clear from the context.)
We also record the signs of the crossings as signs associated to the corresponding edges.
The notions of sign and direction, for an edge of the Seifert graph, are equivalent in the sense that they both flip if we exchange the roles of over- and undercrossing strands at a crossing. 
We note that Seifert graphs are always bipartite, though this property will not play a role in this paper.

It is common to distinguish between \emph{type~I} and \emph{type~II} Seifert circles, where being type~I means that the circle does not separate parts of the diagram from each other. As the type~II circles do separate, they become cut vertices in the Seifert graph.

\begin{ex}
    The standard diagram of the knot $6_1$ has $5$ Seifert circles. The upper right panel of Figure \ref{fig:6_1_link_diagram_and_Seifert} shows them, as well as the edges of the Seifert graph. 
\end{ex}

\subsection{Ribbon graphs}\label{ssec:ribbon_graphs}

We will need yet more structure on our Seifert graphs, namely, the natural ribbon structure induced by the link diagram. 
We will first introduce ribbon graphs as an abstract notion and then explain 
why a Seifert graph is always an example,
as well as how  
that relates to the more concrete (from a knot theory point of view) concept of a Seifert surface. 

\subsubsection{Abstract ribbon graphs}
A \emph{ribbon graph} $G$ is a graph where for each vertex, 
a cyclic permutation of the incident edges is given. For a vertex $u\in V(G)$ and an edge $e$ incident to $u$, we denote by $\next(u,e)$ the edge that follows $e$ in the cyclic order around $u$. For us, $G$ will typically be a directed graph, in which  
case we will talk about ribbon digraphs. The orientation and the ribbon structure are completely unrelated, though.

\begin{ex}
    The lower left panel of Figure \ref{fig:6_1_link_diagram_and_Seifert} shows a digraph. If we take the cyclic order induced by the positive orientation of the plane around each vertex, then this becomes a ribbon digraph. For example, we have $\next(y,d)=f$ and $\next(y,f)=a$.
\end{ex}

\begin{remark}
    It is possible, just like above, to present any ribbon graph with a planar drawing,
    possibly with edges crossing, so that around each vertex, the cyclic ordering is induced by the positive 
    orientation of the plane.
    But in many cases, especially in the context of knot theory, it is more convenient to draw ribbon graphs in such a way that around some vertices, the ribbon structure is induced by the positive orientation, while at other vertices, by the negative orientation.
\end{remark}

Ribbon graphs are equivalent to cellularly embedded graphs into closed oriented surfaces\footnote{Some authors consider ribbon graphs in a more general sense, allowing ``twisted ribbons.'' That is they drop the condition, stated below, that the discs corresponding to vertices and the rectangles corresponding to edges be fused in an orientation-preserving way, resulting in graphs cellularly embedded into non-orientable surfaces. We do not go to this generality in this paper.}: For a graph cellularly embedded into an oriented surface, the positive orientation of the surface determines a cyclic ordering of the incident edges around each vertex. Conversely, to each ribbon graph $G$, one can canonically associate a closed oriented surface $\Sigma_G$ such that the graph is cellularly embedded into the surface, and the positive orientation gives the ribbon structure. 
This is done in the following way: replace each vertex with a `small' oriented disc, whose boundary is also oriented by the usual `outside normal first' rule. Then, thicken each edge into an oriented rectangle, and glue these rectangles, along their `short' sides, to the oriented boundaries of the appropriate discs in the order given by the ribbon structure. We do this in such a way that the orientations of the discs and rectangles match. This way, we get a surface $F_G$ with some $S^1$ boundary components, altogether called the \emph{boundary of the ribbon graph} $G$.
To get the closed surface $\Sigma_G$ of $G$, we glue a disc to each boundary component. 
The \emph{genus} of $G$ is defined as the genus of $\Sigma_G$.
We denote number of boundary components of $G$ by $b(G)$.

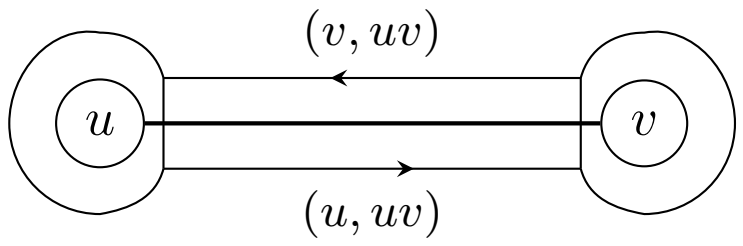
\begin{figure}
    \centering
    \begin{tikzpicture}[decoration={markings,
    mark= at position 0.6 with {\arrow{stealth}}},scale=.6,shift={(0,0)},>=stealth]
    \tikzstyle{v}=[circle,draw]
    \node [v] (u) at (0,0) {$u$};
    \node [v] (v) at (6,0) {$v$};
    \draw [-,thick] (u) to (v);
    \draw [-,postaction={decorate}] (5.3,0.5) to (0.7,0.5);
    \draw [-,postaction={decorate}] (0.7,-0.5) to (5.3,-0.5);
    \draw [-] (5.3,-0.5) to (5.3,0.5);
    \draw [-] (0.7,0.5) to (0.7,-0.5);
    \draw [-] (0.7,0.5) to [out=105,in=0] (0,1) to   [out=170,in=90] (-1,0) to [out=-90,in=180] (0,-1) to [out=10,in=-105] (0.7,-0.5);
    \draw [-] (5.3,0.5) to [out=75,in=180] (6,1) to   [out=10,in=90] (7,0) to [out=-90,in=0] (6,-1) to [out=170,in=-75] (5.3,-0.5);
    \node at (3,-1) {\small $(u,uv)$};
    \node at (3,1) {\small $(v,uv)$};
    \end{tikzpicture}
    \caption{The two half-edges corresponding to an edge (when the orientation is counterclockwise).}
\label{fig:half-edges_of_a_ribbon_graph}
\end{figure}

There is also a more abstract way to think about the boundary of $F_G$ for a ribbon graph $G$.
In the above construction of $F_G$, each edge has two ``sides'' that are the two ``long'' sides of the corresponding rectangle (the ones that are not glued to the discs associated to vertices). We identify these sides with vertex-edge pairs: For an edge $e$ incident to a vertex $u$, we let the vertex-edge pair $(u,e)$ correspond to the side of $e$ 
whose induced orientation (again by the outside normal first rule) points away from $u$. (See Figure \ref{fig:half-edges_of_a_ribbon_graph}.)
We also refer to a vertex-edge pair $(u,e)$, where $e$ is incident to $u$, as a \emph{half-edge} (of $e$).

With this convention, we can describe the boundary components of $G$ using half-edges (and isolated points if there are any). Namely, the boundary component of a half-edge $(u,uw)$ is the set of half-edges $(u,uw),\,\next(w,uw)=:wz,\,\next(z,wz), \dots$, listed until we arrive back to $(u,uw)$.

\begin{ex}\label{ex:6_1_boundary}
    Consider one more time the ribbon digraph 
    depicted in the lower left panel of Figure \ref{fig:6_1_link_diagram_and_Seifert}. 
    Being the Seifert graph of a \emph{knot} diagram, it
    has one boundary component,
    which can be described as the (cyclic) sequence of half-edges $(u,c)$, $(w,b)$, $(x,a)$, $(y,e)$, $(z,f)$, $(y,a)$, $(x,b)$, $(w,c)$, $(u,d)$, $(y,f)$, $(z,e)$, $(y,d)$.
\end{ex}

\subsubsection{Duality}

Any ribbon graph $G$ has a dual $G^*$. 
The vertex set of $G^*$ may be identified with the set of discs that make up $\Sigma_G-F_G$; in other words, with the boundary components of $G$. The edge sets of $G$ and $G^*$ are in bijection: For each edge $e\in G$, the dual $G^*$ has an edge $e^*$ that connects the boundary components lying on the two sides of $e$. In other words, $e^*$ connects the boundary components containing the two half-edges of $e$. In yet another language, we may think of the edges of $G^*$ as the co-cores (median segments parallel to the `short' sides) of the rectangles in the ribbon surface. 

The dual $G^*$ carries a natural ribbon structure. Indeed, for a vertex of $G^*$ (that corresponds to a boundary component of $G$, oriented by the outside normal first rule), take the order in which the corresponding primal edges follow each other along the boundary component. This choice results in $\Sigma_{G^*}=-\Sigma_G$, meaning that there is a natural orientation-reversing identification\footnote{Each edge and its dual edge correspond to the same rectangle, with the roles of the `short' and `long' sides interchanged.} between the surfaces, but that does imply $(G^*)^*=G$.

If $G$ is a directed graph then we get an orientation for $G^*$ as well, using the following rule: 
at the center point of each rectangle inside $F_G$, the pair
\begin{equation}\label{eq:iranyitas}
(\text{direction vector of edge},\text{direction vector of dual edge})
\end{equation}
is a positive basis of the tangent plane. Alternatively, this means that the tail of $e^*$ is the boundary component of $G$ containing the half-edge $(t(e),e)$ and the head of $e^*$ is the boundary component of $G$ containing $(h(e),e)$. Note that 
we have $(G^*)^*=G$ even 
in the directed sense, 
because when the same rectangle is viewed as part of $F_{G^*}$, it is considered with the reverse orientation. 

\begin{figure}[h]
	\centering
	\begin{tikzpicture}[decoration={markings,
			mark= at position 0.6 with {\arrow{stealth}}}]
    \tikzstyle{e}=[circle,scale=.9,draw,color=e]
    \tikzstyle{v}=[circle,scale=.9,draw,color=blue]
        \begin{scope}[scale=.27,shift={(0,0)},>=stealth]
			\node [v] (u) at (10,0) {$u$};
			\node [v] (v) at (0,0) {$v$};
			\node [v] (w) at (2,8) {$w$};
			\node [e] (x) at (5,4) {$x$};
			\node [e] (y) at (5,0) {$y$};
			\node [e] (z) at (5,-4) {$z$};
			\node [e] (t) at (-3,4) {$t$};
			\draw [<-,thick] (u) to (x);
			\draw [->,thick] (u) to (y);
			\draw [->,thick] (u) to (z);
			\draw [->,thick] (v) to (x);
			\draw [<-,thick] (v) to (y);
			\draw [->,thick] (v) to (z);
			\draw [->,thick] (w) to (t);
			\draw [<-,thick] (w) to (x);
			\draw [ultra thick,color=red,postaction={decorate}] (11.5,0) to [out=90,in=-10] (7.2,2.2) to [out=170,in=0] (5,2.5) to [out=180,in=10] (2.6,2.05) to [out=-170,in=50] (-1.3,1.6)  to [out=-130,in=-180] (2.8,-2.6);
			\draw [ultra thick,color=gray,postaction={decorate}] (-0.5,5.5) to [out=-90,in=40] (-2,2.5) to [out=-160,in=-90] (-4.2,4) to [out=90,in=180] (-0.5,6) to [out=0,in=180] (2,6) to [out=0,in=170] (4.2,5.65) to [out=-10,in=140] (6.3,4.8) to [out=-40,in=110] (7.2,2.6);
			\draw [ultra thick,color=gray] (7.3,1.8) to [out=-80,in=80] (7.3,0.3);
			\draw [ultra thick,color=gray] (7.2,-0.3) to [out=-110,in=-70] (2.85,-0.3);
			\draw [ultra thick,color=gray] (2.7,0.3) to [out=110,in=-100] (2.55,1.8);
			\draw [ultra thick,color=gray] (2.6,2.3) to [out=80,in=-100] (3.6,5.5);
			\draw [ultra thick,color=gray,postaction={decorate}] (3.7,6) to [out=80,in=0] (2,9.5) to [out=180,in=90] (-0.5,6.3);
			\draw [ultra thick,color=red] (3.5,-2.55) to [out=10,in=170] (6.8,-2.4) to [out=-10,in=-90] (11.5,0);
			\draw [ultra thick] (7.05,-2.1) to [out=90,in=-60] (7.25,0) to [out=120,in=0] (5,1.2) to [out=180,in=70] (2.7,0) to [out=-110,in=100] (3.1,-2.4) to [out=-80,in=180] (5,-5.5) to [out=0,in=-90] (7.05,-2.75);
			\draw[->,ultra thick] (5.2, -5.5) to (4.8, -5.5);
			\node[circle,draw,fill,scale=.3] at (11.5,0) {};
            \node at (-.5,-5) {$G$};
		\end{scope}
        \begin{scope}[scale=.35,shift={(15,0)},>=stealth]
			\node [circle,draw,color=red,fill=red] (p) at (10,1) {};
			\node [circle,draw,fill] (f) at (7,-2) {};
			\node [circle,draw,color=gray,fill=gray] (sz) at (4,6) {};
			\draw [<-,thick] (p) to[out=100,in=10] (sz);
			\draw [<-,thick] (p) to[out=-160,in=-120] (sz);
			\draw [<-,thick] (p) to[out=-80,in=-10] (f);
			\draw [->,thick] (p) to[out=-110,in=50] (f);
			\draw [->,thick] (sz) to[out=-135,in=150] (f);
			\draw [->,thick] (sz) to[out=-170,in=180] (f);
			\draw [<-,thick] (sz) to[out=-10,in=45] (7.5,2) to[out=-135,in=-90] (sz);
			\draw [<-,thick] (sz) to[out=-30,in=45] (6.5,3) to[out=-135,in=-70] (sz);
			\node at (2.5,-3.5) {$G^*$};
		\end{scope}
	\end{tikzpicture}
	\caption{A ribbon graph and its dual, see Example \ref{ex:duality_ex}.}
	\label{fig:duality_example}
\end{figure}
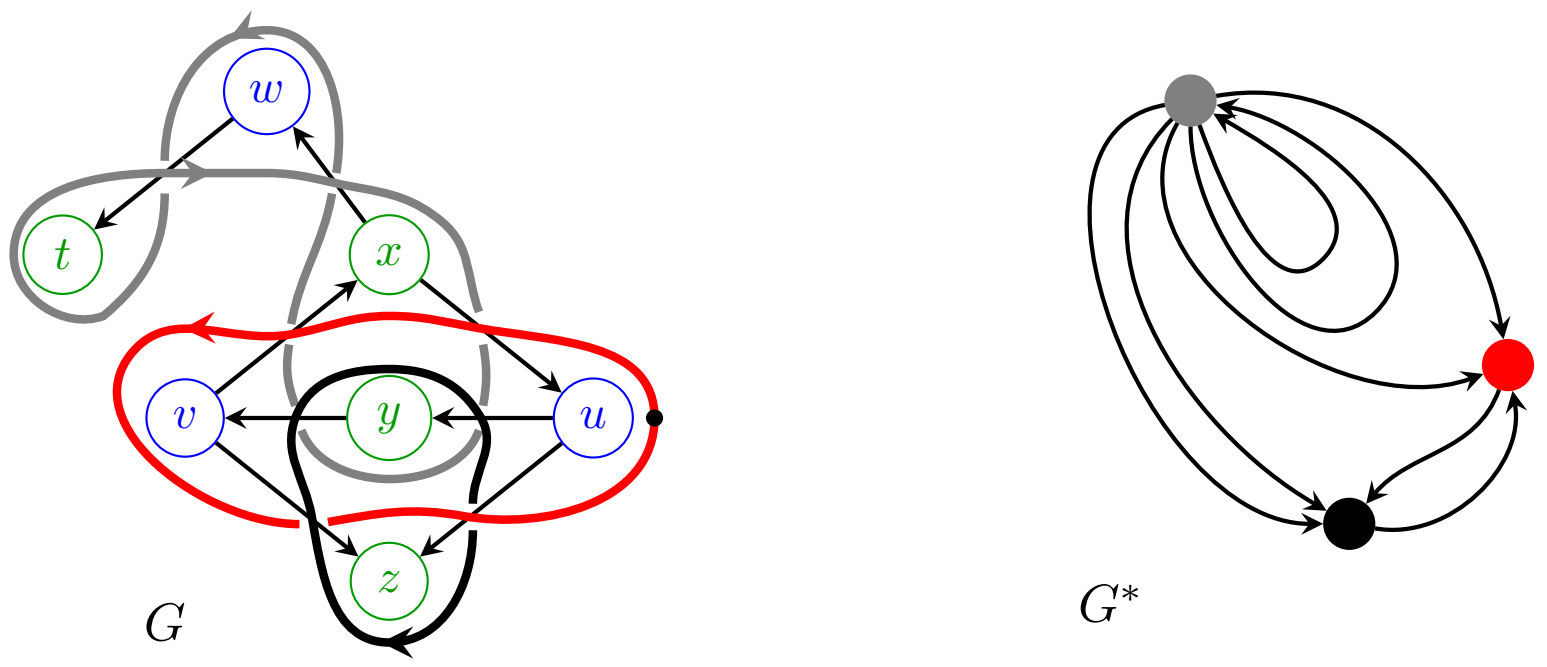

\begin{ex}\label{ex:duality_ex}
    Figure \ref{fig:duality_example} shows a ribbon graph $G$ on the left, where the ribbon structure is induced by the counterclockwise orientation around the blue vertices $\{u,v,w\}$ and by the clockwise orientation around the green vertices $\{x,y,z,t\}$. (In fact $G$ 
    is the Seifert graph of a link diagram, as explained in the next section.) $G$ has three boundary components, indicated by the red, gray, and black curves, and the ribbon structure induces orientations of these boundary components as depicted in the figure. Consequently $G^*$ has three vertices, drawn on the right in the same colors. The ribbon structure of $G^*$ is induced by the counterclockwise orientation around all three vertices.
    Indeed $G$ and $G^*$, as ribbon graphs, are planar in this case.
\end{ex}

\subsubsection{Seifert graphs as ribbon digraphs}
The Seifert graph of a link diagram has a natural ribbon structure:
For any vertex, that is
Seifert circle, traverse it in the direction given by the link. 
Take the edges in the order in which we see the corresponding crossings.
When we talk about Seifert graphs in this paper, we will always think of them as signed, directed graphs with a ribbon structure.

\begin{ex}
    Figure \ref{fig:6_1_link_diagram_and_Seifert} shows a diagram of the $6_1$ knot, with a patrticular orientation, and the upper right panel depicts its Seifert graph. For the Seifert circle denoted by $y$, the diagram induces the counterclockwise orientation. This gives the cyclic ordering $a,e,d,f$ to the edges incident to the vertex corresponding to $y$. The lower left panel of Figure \ref{fig:6_1_link_diagram_and_Seifert} presents the same ribbon digraph with a more conventional drawing, where vertices are represented by nodes, and the ribbon structure is induced by the positive orientation around each vertex. (Note that we had to draw the graph with a crossing in order to represent the ribbon structure this way.)
\end{ex}

When $G$ is the Seifert graph of an oriented link diagram and its ribbon structure is induced as above, then 
\begin{enumerate}
\item the boundary of $F_G$ and the link, as oriented curves, may be identified, and
\item the surface with boundary $F_G$ has a particular embedding in $\R^3$, known as the \emph{Seifert surface} of the diagram. 
\end{enumerate}
To build the latter, first isotope the diagram so that away from the crossings the Seifert circles appear as round circles; then realize the discs in $F_G$ as hemispheres on one side of the projection plane; finally add the rectangles as thin half-twisted strips near the plane. 
This embedding cannot be extended to an embedding of the closed surface $\Sigma_G$. In this paper we will focus on the abstract surface $\Sigma_G$ (and $F_G$) much more than this concrete embedding. 

\begin{remark}
It is well known that a Seifert graph is planar, that is, it can be drawn on the plane without crossings. However, the above defined ribbon structure of a Seifert graph $G$ is in general not planar, i.e., the corresponding canonical surface $\Sigma_G$ is typically not the sphere. 
\end{remark}

\begin{ex}
    As noted in Example \ref{ex:6_1_boundary}, the Seifert graph of the standard diagram of $6_1$ has one boundary component, simply because $6_1$ is a knot. The lower right panel of Figure \ref{fig:6_1_link_diagram_and_Seifert} shows the corresponding Seifert surface 
    using a slightly different presentation, where instead of hemispheres one embeds flat disks, all parallel to the projection plane but at various distances from it.
\end{ex}

Recall that if $G$ is the Seifert graph of a link diagram, then each edge $e$ corresponds to a crossing in the diagram, and the two half-edges correspond to the two strands of the crossing.
Our convention for orienting the Seifert graph means that the half-edge $(h(e),e)$ corresponds to the strand crossing above, and the half-edge $(t(e),e)$ corresponds to the strand crossing below. 

If $G$ is the Seifert graph of a link diagram, then 
the boundary components of $G$, that is 
the vertices of $G^*$, are the components of the link. Specifically, for a knot, $G^*$ has just one vertex.
The edges of $G^*$ appear in $\R^3$ as segments properly embedded in the Seifert surface $F_G$; namely, they are chords that connect the over- and understrands at each crossing, oriented toward the overstrand.

\subsection{Minors and special subgraphs of ribbon graphs} 

\subsubsection{Deletion}
We say that a ribbon graph is \emph{connected} if the underlying abstract graph is connected. We say that a subgraph $G'$ of a ribbon graph $G$ is \emph{spanning} if it is spanning in the abstract sense, that is, if $G'$ contains each vertex of $G$ and they have the same number of connected components.

If $G'$ is a subgraph of a ribbon graph $G$, we get an induced ribbon structure for $G'$ as well by appropriately restricting each cyclic order. 
To be concrete, when
we delete one edge $uv$, this works in the following way.
The cyclic order around any vertex $w\notin \{u,v\}$ is not modified. For $u$, if $\next_G(u,e)=uv$, then let $\next_{G-uv}(u,e)=\next_G(u,uv)$, otherwise let $\next_{G-uv}(u,e)=\next_G(u,e)$, and similarly for $v$. 

\subsubsection{Contraction}
Furthermore, when $G'$ is derived from $G$ by contracting an edge $e$, we once again get a ribbon structure for $G'$.
Contraction in a ribbon graph is slightly more complicated than deletion. 
The main point is that we want deletion and contraction to be dual operations of each other, just like in the abstract setting. So we define $G/e:=(G^*-e^*)^*$. 

One can check that this is equivalent to the following description:
If $e=uv\in G$ is not a loop, then the contraction works as one would intuitively imagine: The vertex set of $G/e$ is obtained by gluing $u$ and $v$ into a vertex $w$. The edges incident to $w$ are the edges incident to either $u$ or $v$, except for $e$. For some $f\neq e$ incident to $u$, let $\next_{G/e}(w,f)=\next_G(u,f)$ if $\next_G(u,f)\neq e$. If $\next_G(u,f)=e$, then let $\next_{G/e}(w,f)=\next_G(v,e)$. We use the analogous definition for edges incident to $v$ in $G$. In this case the surfaces $F_G$ and $F_{G/e}$ are naturally homeomorphic.

If $e$ is a loop attached to the vertex $w$, then in $G/e$, the vertex $w$ is substituted by two vertices $u$ and $v$. The edge $e$ occurs twice in the cyclic order of the edges incident to $w$, and these two occurrences determine two intervals. Let one of the intervals be the set of edges incident to $u$ and the other interval be the set of edges incident to $v$ (in both cases, the cyclic order is obtained from the linear order of the interval). 
In topological language, the ribbon surface $F_G$ undergoes surgery where the annulus, comprising the disk of $w$ and the rectangle of $e$, is replaced by two disks, representing $u$ and $v$.
This applies even if one of the intervals mentioned above is empty. Then, either $u$ or $v$ becomes an isolated vertex in $G/e$. Note that in particular, the contraction of a loop edges $e$ can result in $G/e$ having one more connected component than $G$.

For an example, see Figure \ref{fig:Jaeger_elimination_example}. This figure shows three Seifert graphs. The Seifert graph in the third panel is obtained by contracting the edges $\{ux,xv, vy, yu\}$ one-by-one in the graph of the second panel.

Note that even though a contraction modifies the vertex set, the half-edges of $G/e$ can still be identified with half-edges of $G$. Moreover, the boundary components of $G$ are in a one-to-one correspondence with the boundary components of $G/e$, by 
thinking of them as sets of half-edges of $G$ and then 
restricting to half-edges of $G/e$. (This follows easily from the definition $G/e=(G^*-e^*)^*$, since the boundary components of $G/e$ correspond to the vertices of $G^*-e^*$ and those are the same as the vertices of $G^*$.)

We refer to a ribbon graph, obtained from $G$ by a sequence of deletions and contractions, as a \emph{minor} of $G$.
It turns out that when a minor is formed, the order of the deletions and contractions does not matter \cite{joandiain}.

\subsubsection{Separating cycles}
We say that a subgraph $C$ of the ribbon graph $G$ is a \emph{separating cycle} of $G$, if $C$ is an inclusionwise minimal subset of edges such that $\Sigma_G-C$ has more components than $\Sigma_G$. 
Here by $\Sigma_G-C$ we mean the result of removing the interiors of all rectangles corresponding to edges of $C$, as well as the interiors of all discs corresponding to the incident vertices of $G$.
Equivalently, in more graph theoretic terms, $C$ is a separating cycle if the corresponding set of edges $C^*$ is an elementary cut in $G^*$, or put in yet another way, if $C$ is an inclusionwise minimal set of edges such that $G/C$ has more connected components than $G$.

\begin{ex}
    The ribbon graph $G$ in the lower left panel of Figure \ref{fig:6_1_link_diagram_and_Seifert} has no separating cycles, since its dual graph has one single vertex. However, if we remove the edge $c$, the set $\{e,f\}$ becomes a separating cycle in $G-c$.

    The ribbon graph $G$ on the left panel of Figure \ref{fig:duality_example} (see Example \ref{ex:duality_ex} for the description of the ribbon structure) has $\{ux,xv,vz,zu\}$ as a separating cycle. Note that the corresponding edge set is a cut in the dual $G^*$ (separating the gray vertex from the red and black vertices).
\end{ex}

\begin{claim}\label{cl:separating_cycle_in_subgraph}
	Let $G$ be a ribbon graph and $G'\subseteq G$ a subgraph.
	If $C\subseteq G'$ is a separating cycle in $G$, then some $C'\subseteq C$ is a separating cycle in $G'$.
\end{claim}

\begin{proof}
    We think of $G'$ as having the same vertex set as $G$ as the addition of isolated vertices does not modify whether some $C\subseteq G'$ is a separating cycle.
    
    Let us denote $F=E(G)-E(G')$. 
    Then we have $(G')^*=G^*/F^*$. 
    If $C$ is a separating cycle in $G$, then $C^*$ is a cut in $G^*$. Since $C\subset G'$, we have $F \cap C=\varnothing$. 
    $C^*$ stays a cut after contracting an edge set disjoint from it. (Indeed, the contraction of an edge either merges two vertices on the same side of the cut, or splits a vertex into two vertices. In both cases, $C^*$ is still the set all of edges connecting vertices on the two sides.) If $C^*$ is still an elementary cut in $G^*/F^*$, then $C$ itself is a separating cycle in $G'$. If $C^*$ is not elementary in $G^*/F^*$, then it has a subset that is an elementary cut. 
\end{proof}

\subsubsection{Rooted digraphs}\label{ss:rooted_digraphs}
In many of our arguments we need to pick a point on the boundary of a ribbon graph. 
Informally, we call a ribbon graph $G$ \emph{rooted} if a boundary point of $G$ is fixed. The fixed point on the boundary is called the \emph{root}, or the \emph{basepoint}. We will always imagine that the root is where the ``rectangle'' of an edge meets the ``disc'' of a vertex (or on the boundary of a disc of an isolated vertex), which leads to the following formal definition.

\begin{defn}
A \emph{rooted ribbon graph} is a ribbon graph $G$ together with a half-edge $(v_0,e_0)$, where $e_0\in E(G)$ and $v_0$ is an endpoint of $e_0$, or a ribbon graph $G$ together with an isolated vertex $v_0$. We call $(v_0,e_0)$ (of $v_0$) the \emph{root} or \emph{basepoint}. Moreover, we refer to the boundary component of $G$ containing the root $(v_0,e_0)$ (or $v_0$) as the \emph{root boundary component} of the rooted ribbon graph. When treated as 
a vertex of $G^*$, we call this the root of $G^*$. 
\end{defn}

These concepts are independent of edge directions but may apply just as well to ribbon digraphs.

If $G$ is a rooted ribbon digraph, we can also think of a minor of $G$ as a rooted ribbon digraph. 
If we think of deletion and contraction as topological operations then it is fairly easy to see where the induced basepoint should be. The corresponding graph-theoretical formulation is as follows:
Let the root of $G$ be $(v_0,e_0)$ with $e_0=v_0u$. If $e\neq e_0$, then we take the root of both $G-e$ and $G/e$ to be $(v_0,e_0)$. (In the case of $G/e$ the vertex $v_0$ might have been modified if it is incident to $e$, but we can then consider the 
derived vertex.) If $e=e_0$, then for $G-e_0$, we take $(v_0,\next(v_0,e_0))$ as the root of $G-e_0$ (or $v_0$ if it becomes an isolated vertex). For $G/e_0$, let $w$ be the vertex of $G/e_0$ corresponding to $u$ that is incident to $\next_G(u,e_0)$. Let the root of $G/e_0$ be $(w,\next_G(u,e_0))$, or $w$ if $w$ becomes an isolated vertex. (That is, if $e_0=v_0u$ is not a loop, then $w$ is the vertex obtained by merging $v_0$ and $u$. If $e_0$ is a loop incident to $v_0$, then $w$ is the vertex corresponding to $v_0$ such that $\next_G(u,e_0)$ is incident to $w$.) 

\begin{ex}
    Consider the ribbon graph $G$ in the left panel of Figure \ref{fig:duality_example}. If we take $(u,ux)$ as the root, that corresponds to choosing the marked point near $u$ on the red boundary component for basepoint. In this case, the red component is the root boundary component of $G$, and for $G^*$, the root is the red vertex. 
\end{ex}

\subsubsection{Quasi-trees} \label{ss:quasi-trees_and_tours}
Spanning trees, that is connected and cycle-free spanning subgraphs, are among the most important objects
of abstract graph theory. In the context of ribbon graphs, 
their class has the following natural extension.

\begin{defn}
    A \emph{quasi-tree} is a spanning subgraph of a ribbon graph that has a connected boundary.
\end{defn} 

A spanning tree is always a quasi-tree, but if $\Sigma_G$ has positive genus then $G$ contains quasi-trees that are not trees. For example, Figure \ref{fig:61_Jaeger_kvazifai} shows three of the (altogether nine) quasi-trees of the Seifert graph of the standard $6_1$ diagram, two of whom are trees, but the third one has higher genus. In many senses quasi-trees play the same role for ribbon graphs that spanning trees play for graphs. In particular, while graphs have associated matroids 
and spanning trees are the bases of this so called graphic matroid, ribbon graphs induce a delta-matroid structure, where quasi-trees are the bases of the ribbon-graphic delta-matroid.

Fixing a quasi-tree enables us to define a cyclic order on the half-edges of a ribbon graph. This becomes a linear order once we fix a starting half-edge (root). 

\begin{defn}
    Let $G$ be a connected rooted ribbon graph and let $K$ be a quasi-tree of $G$. We define $<_K$ 
    by listing all half-edges as follows.
    Let the root half-edge $(v_0,e_0)$ be the smallest element in $<_K$. If the ordering is defined until the half-edge $(u,uv)$ and $uv\notin K$, then let the next half-edge be $\next(u,uv)$. If the ordering is defined until the half-edge $(u,uv)$ and $uv\in K$, then let the next half-edge be $\next(v,uv)$. If at some point the next half-edge would be $(v_0,e_0)$, we stop.
\end{defn}

The idea is just a slight extension of Bernardi's notion of the tour of a spanning tree \cite{Bernardi_first}.
Geometrically, just imagine tracing the boundary of $F_K$; this will line up all half-edges exactly when $K\subset G$ is a quasi-tree. We remark that the two halves of an edge may appear in $<_K$ with a significant separation, regardless of whether the edge is in $K$ or not.

\begin{prop}
$<_K$ is a linear ordering of the half-edges of $G$.
\end{prop}

\begin{proof}
    $<_K$ is defined by traversing the boundary of $K$, and listing the half-edges in the order in which they appear during this process, including the edges of $G-K$. Since $K$ is a quasi-tree, by definition all half-edges belonging to edges in $K$ are listed. In the meantime, the traversal also makes a full turn around (the disc corresponding to) each vertex of $G$, whence also for each edge $e\notin K$, both half-edges of $e$ are listed by the process.
\end{proof}

\begin{defn}
The \emph{tour} of the quasi-tree $K$, in the rooted ribbon graph $G$, is the list of half-edges of $G$ according to the order $<_K$, together with the information whether the given edge is in $K$ or not.
\end{defn}

The ordering $<_K$ also induces an ordering of the edges of $G$: We put $uv <_K xy$ if $uv$ has a half-edge, say  $(u,uv)$, such that $(u,uv) <_K (x,xy)$ and $(u,uv) <_K (y,xy)$. With a slight abuse, we will use the same notation for the orderings of the edges and half-edges.

For two quasi-trees $K$ and $L$ of $G$, and $e\in E(G)$, we will often say that ``the tours of $K$ and $L$ agree until $e$.'' This means that the orders $<_K$ and $<_L$,
as lists of half-edges starting with the root, are identical up to (and including) the first appearance of a half-edge of $e$ (and possibly beyond), and for $f<_K e$ we have $f\notin K\Delta L$, where $\Delta$ denotes symmetric difference. Note that $e\in K\Delta L$ is allowed.
We also say that $e$ is the \emph{first difference} between the quasi-trees $K$ and $L$ if the tours of $K$ and $L$ agree until $e$ but $e\in K\Delta L$.

\begin{ex} 
Figure \ref{fig:61_Jaeger_kvazifai} shows three quasi-trees of a rooted ribbon graph where the ribbon structure is induced by the positive orientation of the plane around each vertex, and the root half-edge is $(u,c)$ (we indicate this by the little arrow beneath $u$). In the middle panel, the non-dotted edges form the quasi-tree $K_2$. This induces the half-edge order $(u,c)$, $(w,b)$, $(x,a)$, $(x,b)$, $(w,c)$, $(u,d)$, $(y,f)$, $(z,e)$, $(z,f)$, $(y,a)$, $(y,e)$, $(y,d)$. Consequently, we have $c <_{K_2} b <_{K_2} a <_{K_2} d <_{K_2} f <_{K_2} e$.
The tours of $K_2$ and $K_3$ agree until $g$ for each $g \leq_{K_2} a$. The first difference between $K_2$ and $K_3$ is $a$.
\end{ex}

\section{Jaeger quasi-trees}\label{s:Jaeger_trees_def_and_formula}

In this section we define our key notion of Jaeger quasi-tree. This first requires a discussion of \emph{naive Jaeger quasi-trees}, which form a superclass of Jaeger quasi-trees, and the simplicity of whose definition makes them useful for us.

\begin{defn}[naive Jaeger quasi-tree] 
A quasi-tree $K$ of the connected rooted ribbon digraph $G$ is called a \emph{naive Jaeger quasi-tree} (for the root $(v_0,e_0)$) if for each $e\notin K$, we have $(t_G(e),e)<_K (h_G(e),e)$.
\end{defn}

In Definition \ref{def:jaegquastree} below, Jaeger quasi-trees are defined by adding one more condition. In cases arising from knot theory this has to do with simplifying the diagram in certain situations, and requiring that the quasi-tree is also a quasi-tree after the simplification. The general definition is abstracted from this requirement. 
Since taking a closer look at naive Jaeger quasi-trees will help make this extra ingredient sound more natural, we start with a necessary condition for the existence of naive Jaeger quasi-trees. 

\begin{remark}
We worked with the notion of a ``Jaeger tree'' in several earlier papers, perhaps most notably in \cite{semibalanced} where we showed that they give rise to a dissection of the root polytope for a wide class of directed bipartite graphs.
In the terminology of this paper, those ``Jaeger trees'' 
are the spanning trees of the ribbon graph that are naive Jaeger. In Section \ref{ssec:naive_J_trees_are_typically_J_trees} we will show that (with a mild assumption on the underlying graph), naive Jaeger spanning trees are in fact Jaeger (in the terminology of the current paper).
\end{remark}

\begin{defn}[bad separating cycle]
	A separating cycle $C$ of a rooted ribbon digraph $G$ is a \emph{bad separating cycle} if the corresponding cut $C^*$ of $G^*$ is such that each edge of $C^*$ is oriented toward the component of $G^*-C^*$ containing the root.

Equivalently, cf.\ Figure \ref{fig:half-edges_of_a_ribbon_graph} and \eqref{eq:iranyitas}, this means that among the two components of $\Sigma_G-C$, the one containing $(v_0,e_0)$ also contains the half-edge $(h(e),e)$ for each $e\in C$. If $G$ is not connected, then we consider $\varnothing$ as a bad separating cycle of $G$.
\end{defn}

Recall that half-edges of $G$ are defined as segments along the boundary of $F_G$. Thus it makes sense to say that $\Sigma_G-C$ contains a half-edge, even if the corresponding edge is part of $C$.

\begin{ex}
    If one takes $(u,ux)$ as the root for the graph on Figure \ref{fig:duality_example}, then $C=\{ux,xv,vy,yu\}$ becomes a bad separating cycle. (See Example \ref{ex:duality_ex} for the description of the ribbon structure.) In this case, the red node becomes the root in the dual, and indeed, $C^*$ is a cut that is oriented towards the side containing the red node.
\end{ex}

\begin{prop}\label{prop:no_naive_Jaeger_quasi-tree_if_separating_directed_cut}
	Let $G$ be a rooted ribbon digraph. If $G$ has a bad separating cycle $C$, then $G$ has no naive Jaeger quasi-trees.
\end{prop}

\begin{proof}
    If $G$ is not connected, then it has no quasi-trees.
	Let $G$ be connected, and suppose for a contradiction that $G$ has a bad separating cycle $C$ and a naive Jaeger quasi-tree $K$.
	If $C$ was part of $K$ then by Claim \ref{cl:separating_cycle_in_subgraph} a subset of $C$ would be a separating cycle in $K$ as well, but quasi-trees cannot contain such.
    Thus there needs to be an edge in $C-K$, and we may let $f$ be the first such edge according to $<_K$. 
	As the base component of $\Sigma_G-C$ contains $(h(f),f)$ we have $(h(f),f) <_K (t(f),f)$,
    contradicting the assumption that $K$ is a naive Jaeger quasi-tree.
\end{proof}

Bad separating cycles will play an important role throughout the paper, especially the so-called ``nearest'' ones. 

\begin{defn}
    We call a bad separating cycle $C$ of the rooted ribbon digraph $G$ a \emph{nearest bad separating cycle} in $G$, if for any other bad separating cycle $C'$, the connected component of $G^*-C^*$ containing the root is a subset of the connected component of $G^*-(C')^*$ containing the root. 
\end{defn}

\begin{prop}
	Let $G$ be a rooted ribbon digraph.
	If $G$ has a bad separating cycle, then it has a unique nearest bad separating cycle. 
\end{prop}

\begin{proof}
	Let $C_1, \dots, C_k$ be the set of all bad separating cycles of $G$. 
    Let $(B_0^i,B_1^i)$ be the shores of the cuts $C_i^*$ in $G^*$, with $B_0^i$ containing the root. Since $C_i$ is a bad separating cycle, each edge in $C_i^*$ has its head in $B_0^i$.

    Take $B_0:=B_0^1 \cap \dots \cap B_0^k$. Then the cut $C^*$ of $G^*$ with shores $(B_0, V(G^*)-B_0)$ is such that each edge has its head in $B_0$. Indeed, if an edge $f^*$ has one endpoint in $B_0$ and the other endpoint in $V(G^*)-B_0$, then one endpoint is in $B_0^i$ for each $i=1, \dots, k$, and the other endpoint is in $B_1^j$ for some $j$. This implies $f^*\in C_j^*$, whence $f^*$ has its tail in $B_1^j$, and in particular, in $V(G^*)-B_0$.

    We obtain that the dual of $C^*$ is a bad separating cycle and such that $G^*-C^*$ has a minimal root component among bad separating cycles of $G$.
\end{proof}

Let us also remark that finding the nearest bad separating cycle is computationally tractable.

\begin{prop}\label{prop:bad_sep_cycle_computation}
    Let $G$ be a rooted ribbon digraph. It can be decided in polynomial time in $|E(G)|$ whether $G$ has a bad separating cycle, and if the answer is yes, the nearest bad separating cycle can be computed in polynomial time in $|E(G)|$.
\end{prop}

\begin{proof}
    One can compute $G^*$ in polynomial time. Let $s$ be the vertex of $G^*$ that corresponds to the root boundary component of $G$. To check if $G$ has a bad separating cycle, one has to check if $G^*$ contains a cut with each edge directed toward the side containing $s$. The existence of such a cut can be decided by a breadth first search from $s$: If every vertex of $G^*$ is reached from $s$, then $G$ has no bad separating cycle. If only a subset $S\subsetneq V(G^*)$ is reached, then the cut determined by $(S,V-S)$ is the nearest bad separating cycle in $G$.
\end{proof}

After this detour, we define Jaeger quasi-trees. This will be our most important definition. A Jaeger quasi-tree will always be a naive Jaeger quasi-tree, but we need a slightly stronger condition, summarized in 
a procedure called the \emph{Jaeger test}. 
Let $G$ be a connected rooted ribbon digraph, and $K$ one of its quasi-trees.
The test takes the edges of $G$ in the order $<_K$ determined by the tour of $K$, and checks if certain conditions are fulfilled. 
If some edge does not satisfy the conditions, then the test stops, and concludes that $K$ is not a Jaeger quasi-tree. If each edge $e\in G$ passes the test, then $K$ is said to be a Jaeger quasi-tree. 

To state the condition that the edges need to satisfy, we recursively define a series of auxiliary graphs $\AuxContr{G,K,e_0}, \AuxContr{G,K,e_1}, \dots, \AuxContr{G,K,e_k}$, where $e_0, e_1, \dots, e_k$ are the first few edges of $G$ according to $<_K$ 
and they, with the possible exception of $e_k$, pass the Jaeger test.
Each $\AuxContr{G,K,e}$ will be a rooted ribbon digraph, in fact a reorientation of a minor of $G$, 
and it will have the property that if an edge $f\in G$ is deleted when we form $\AuxContr{G,K,e}$, or if it is present but its orientations differ in $G$ and in $\AuxContr{G,K,e}$, then $f <_K e$ (this need not be true for the edges that are contracted).

Let us introduce some more notation. If $G'$ is a reoriented minor of a ribbon digraph $G$, then by $e\in G'$ we mean that $e$ is neither deleted nor contracted in $G'$, but its orientations might be different in $G$ and in $G'$. The orientation of $e$ is indicated through the usage of $t_{G'}(e)$ and $h_{G'}(e)$.

\medskip
\noindent
\textbf{Jaeger test:}
We define $\AuxContr{G,K,e_0}=G$.

Suppose that we have defined $\AuxContr{G,K,e}$, and let $f$ be the edge following $e$ in $<_K$. If $e$ is the last edge according to $<_K$, then we use the symbol $f=\infty$. The following case analysis is called the \emph{Jaeger test of $K$ for the edge $e$}, and we define $\AuxContr{G,K,f}$ based on its result. At this point, before reading on, all we know about $\AuxContr{G,K,e}$ is that it is the result of deleting some edges of $G$, reversing some other edges of $G$ (both these types strictly precede $e$ in $<_K$) and contracting some further edges of $G$ (which may include $e$ or even some edges after $e$ in $<_K$).

\begin{enumerate}[label=(\alph*)]
\item\label{e_notin_K_contracted} If $e\notin K$ and $e$ was contracted when we formed $\AuxContr{G,K,e}$, then we stop, and conclude that 
$e$ fails the Jaeger test.
In this case $\AuxContr{G,K,f}$ remains undefined for $f$ and for all subsequent edges in $<_K$.

\item\label{e_notin_K_head} If $e\notin K$, the edge $e$ is not contracted in $\AuxContr{G,K,e}$, and $(h_G(e),e)<_K (t_G(e),e)$, then we stop, and conclude that 
$e$ fails the Jaeger test.
In this case $\AuxContr{G,K,f}$ remains undefined for $f$ and for all subsequent edges in $<_K$.

\item\label{e_notin_K_tail} If $e\notin K$, the edge $e$ is not contracted in $\AuxContr{G,K,e}$ and $(t_G(e),e)<_K (h_G(e),e)$, then $e$ passes the Jaeger test and we let $\AuxContr{G,K,f}=\AuxContr{G,K,e}-e$. 

\item\label{e_in_K_contracted} If $e\in K$ and $e$ is contracted in $\AuxContr{G,K,e}$, then $e$ passes the Jaeger test and we let $\AuxContr{G,K,f}=\AuxContr{G,K,e}$. 

\item\label{e_in_K_head} If $e\in K$, the edge $e$ is not contracted in $\AuxContr{G,K,e}$ and $(h_G(e),e)<_K (t_G(e),e)$, then $e$ passes the Jaeger test and we let $\AuxContr{G,K,f}=\AuxContr{G,K,e}$. 

\item\label{e_in_K_tail_active} If $e\in K$, the edge $e$ is not contracted in $\AuxContr{G,K,e}$ and $(t_G(e),e)<_K (h_G(e),e)$, moreover, $\AuxContr{G,K,e}-e$ has a bad separating cycle, and the nearest bad separating cycle of $\AuxContr{G,K,e}-e$ is $C$ (note that $C$ might be the empty set if $e$ is a bridge in $\AuxContr{G,K,e}$), then $e$ passes the Jaeger test and we let $\AuxContr{G,K,f}=\AuxContr{G,K,e}/C$.

\item\label{e_in_K_tail_passive} If $e\in K$, the edge $e$ is not contracted in $\AuxContr{G,K,e}$ and $(t_G(e),e)<_K (h_G(e),e)$, moreover, $\AuxContr{G,K,e}-e$ has no bad separating cycles, then $e$ passes the Jaeger test and we let $\AuxContr{G,K,f}$ be obtained from $\AuxContr{G,K,e}$ by reversing $e$. 
\end{enumerate}

We think of each $\AuxContr{G,K,e}$ as a rooted ribbon digraph, where the root is induced from the root of $G$ as explained in Subsection \ref{ss:rooted_digraphs}.

\begin{defn}[Jaeger quasi-tree]\label{def:jaegquastree}
Let $G$ be a connected rooted ribbon digraph. A quasi-tree $K$ of $G$ is a Jaeger quasi-tree of $G$ if $K$ passes the Jaeger test for each edge $e\in G$.
\end{defn}

We denote the set of Jaeger quasi-trees of a rooted ribbon digraph $G$ by $\jaeg(G)$. If we want to emphasize that $G$ has root $(v_0,e_0)$, we use $\jaeg(G,v_0,e_0)$.

Let us focus on a special class of Jaeger quasi-trees that will be especially important. A Jaeger quasi-tree is called a \emph{Jaeger tree} if it is cycle-free (that is, if its underlying undirected graph is a tree). These are exactly the Jaeger quasi-trees with the minimum number of edges. We denote the set of Jaeger trees of $G$ by $\jaegtree(G)$ (and by $\jaegtree(G,v_0,e_0)$ if we want to emphasize the root).

We say that a quasi-tree $K$ \emph{passes the Jaeger test until $e$} if $K$ passes the Jaeger test for each $f <_K e$. In other words, $\AuxContr{G,K,f}$ is defined for each $f\leq_K e$.

\begin{remark}\label{rem:osszehuzott_kor_bekerül_a_kvazifaba}
Note that if a Jaeger quasi-tree $K$ passes the Jaeger test for some edge  $e\in K$ via case \ref{e_in_K_tail_active}, then for the nearest bad separating cycle $C$ of $\AuxContr{G,K,e}-e$, we need $C\subseteq K$, otherwise at some later point we would encounter an edge falling under case \ref{e_notin_K_contracted}. 
\end{remark}

\begin{figure}[h]
	\centering
	\begin{tikzpicture}[decoration={markings,
			mark= at position 0.6 with {\arrow{stealth}}}]
    \tikzstyle{e}=[circle,scale=.9,draw,color=e]
    \tikzstyle{v}=[circle,scale=.9,draw,color=blue]
    \begin{scope}[scale=.27,shift={(0,0)},>=stealth]
    \node [v] (u) at (10,0) {$u$};
    \node [v] (v) at (0,0) {$v$};
    \node [v] (w) at (2,8) {$w$};
    \node [e] (x) at (5,4) {$x$};
    \node [e] (y) at (5,0) {$y$};
    \node [e] (z) at (5,-4) {$z$};
    \node [e] (t) at (-3,4) {$t$};
    \draw [<-,thick] (u) to (x);
    \draw [->,thick] (u) to (y);
    \draw [->,thick] (u) to (z);
    \draw [->,thick] (v) to (x);
    \draw [<-,thick] (v) to (y);
    \draw [->,thick,dashed] (v) to (z);
    \draw [->,thick] (v) to (t);
    \draw [->,thick] (w) to (t);
    \draw [<-,thick] (w) to (x);
    \draw [thick,postaction={decorate},color=gray] (11.5,0) to [out=90,in=-10] (7.2,2.2) to [out=170,in=0] (5,2.5) to [out=180,in=10] (2.6,2.05) to [out=-170,in=-10] (-1.25,2.1);
    \draw [thick,postaction={decorate},color=gray] (-2,2.2) to [out=170,in=-90] (-4.2,4) to [out=90,in=180] (-0.5,6) to [out=0,in=180] (2,6) to [out=0,in=170] (4.2,5.65) to [out=-10,in=140] (6.3,4.8) to [out=-40,in=110] (7.2,2.6);
    \draw [thick,color=gray] (7.3,1.8) to [out=-80,in=80] (7.3,0.3);
    \draw [thick,postaction={decorate},color=gray] (7.2,-0.3) to [out=-110,in=-70] (2.85,-0.3);
    \draw [thick,color=gray] (2.7,0.3) to [out=110,in=-100] (2.55,1.8);
    \draw [thick,postaction={decorate},color=gray] (2.6,2.3) to [out=80,in=-100] (3.6,5.5);
    \draw [thick,postaction={decorate},color=gray] (3.7,6) to [out=80,in=0] (2,9.5) to [out=180,in=90] (-0.5,6.3);
    \draw [thick,postaction={decorate},color=gray] (-0.5,5.5) to [out=-90,in=70] (-1.6,2) to [out=-110,in=120] (-1.5,-0.3) to [out=-50,in=-180] (2.8,-2.6);
    \draw [thick,color=gray] (3.5,-2.55) to [out=10,in=170] (6.8,-2.4) to [out=-10,in=-90] (11.5,0);
    \draw [thick,color=gray] (7.05,-2.1) to [out=90,in=-60] (7.25,0) to [out=120,in=0] (5,1.2) to [out=180,in=70] (2.7,0) to [out=-110,in=100] (3.1,-2.4) to [out=-80,in=180] (5,-5.2) to [out=0,in=-90] (7.05,-2.75);
    \draw[->,thick,color=gray] (4.9, 1.2) to (5.1, 1.2);
    \node[circle,draw,fill,scale=.3] at (11.5,0) {};
	\draw[->,thick] (11.5, 0) to (11.5, 1);
    \node at (0.5,-5) {$G$ and $K$};
    \end{scope}    
	
        \begin{scope}[scale=.27,shift={(17,0)},>=stealth]
			\node [v] (u) at (10,0) {$u$};
			\node [v] (v) at (0,0) {$v$};
			\node [v] (w) at (2,8) {$w$};
			\node [e] (x) at (5,4) {$x$};
			\node [e] (y) at (5,0) {$y$};
			\node [e] (z) at (5,-4) {$z$};
			\node [e] (t) at (-3,4) {$t$};
			\draw [<-,thick] (u) to (x);
			\draw [->,thick] (u) to (y);
			\draw [->,thick] (u) to (z);
			\draw [->,thick] (v) to (x);
			\draw [<-,thick] (v) to (y);
			\draw [->,thick] (v) to (z);
			\draw [->,thick] (w) to (t);
			\draw [<-,thick] (w) to (x);
			\draw [thick,postaction={decorate},color=gray] (11.5,0) to [out=90,in=-10] (7.2,2.2) to [out=170,in=0] (5,2.5) to [out=180,in=10] (2.6,2.05) to [out=-170,in=50] (-1.3,1.6)  to [out=-130,in=-180] (2.8,-2.6);
			\draw [thick,postaction={decorate},color=gray] (-0.5,5.5) to [out=-90,in=40] (-2,2.5) to [out=-160,in=-90] (-4.2,4) to [out=90,in=180] (-0.5,6) to [out=0,in=180] (2,6) to [out=0,in=170] (4.2,5.65) to [out=-10,in=140] (6.3,4.8) to [out=-40,in=110] (7.2,2.6);
			\draw [thick,color=gray] (7.3,1.8) to [out=-80,in=80] (7.3,0.3);
			\draw [thick,postaction={decorate},color=gray] (7.2,-0.3) to [out=-110,in=-70] (2.85,-0.3);
			\draw [thick,color=gray] (2.7,0.3) to [out=110,in=-100] (2.55,1.8);
			\draw [thick,postaction={decorate},color=gray] (2.6,2.3) to [out=80,in=-100] (3.6,5.5);
			\draw [thick,postaction={decorate},color=gray] (3.7,6) to [out=80,in=0] (2,9.5) to [out=180,in=90] (-0.5,6.3);
			\draw [thick,color=gray] (3.5,-2.55) to [out=10,in=170] (6.8,-2.4) to [out=-10,in=-90] (11.5,0);
			\draw [thick,color=gray] (7.05,-2.1) to [out=90,in=-60] (7.25,0) to [out=120,in=0] (5,1.2) to [out=180,in=70] (2.7,0) to [out=-110,in=100] (3.1,-2.4) to [out=-80,in=180] (5,-5.2) to [out=0,in=-90] (7.05,-2.75);
			\draw[->,thick,color=gray] (4.9, 1.2) to (5.1, 1.2);
			\node[circle,draw,fill,scale=.3] at (11.5,0) {};
			\draw[->,thick] (11.5, 0) to (11.5, 1);
            \node at (-.5,-5) {$\AuxContr{G,K,vt}-vt$};
		\end{scope}
		\begin{scope}[scale=.27,shift={(34,0)},>=stealth]
			\node [v] (uv) at (2,0) {$u'$};
			\node [v] (w) at (2,8) {$w$};
			\node [e] (xy) at (5,4) {$x'$};
			\node [e] (z) at (5,-4) {$z$};
			\node [e] (t) at (-3,4) {$t$};
			\draw [->,thick] (uv) to (t);
			\draw [->,thick] (uv) to [out=-70,in=120] (3, -2.3) to [out=-60,in=140] (z);
			\draw [->,thick] (uv) to [out=0,in=60] (6.2, -2.3) to [out=-120,in=60] (z);			\draw [->,thick] (w) to (t);
			\draw [<-,thick] (w) to (xy);
			\draw [thick,postaction={decorate},color=gray] (8.5,0) to [out=90,in=0] (5,1.2) to [out=180,in=0] (0,2);
			
			\draw [thick,postaction={decorate},color=gray] (-1,2) to [out=170,in=-90] (-4.2,4) to [out=90,in=180] (-0.5,6) to [out=0,in=180] (2,6) to [out=0,in=170] (4.2,5.65) to [out=-10,in=120] (6.3,4.8) to  [out=-70,in=0] (5,2.3) to [out=170,in=-110] (3.6,5.5);
			\draw [thick,postaction={decorate},color=gray] (3.7,6) to [out=80,in=0] (2,9.5) to [out=180,in=90] (-0.5,6.3);
			\draw [thick,postaction={decorate},color=gray] (-0.5,5.5) to [out=-90,in=90] (-0.5,2) to [out=-90,in=180] (2.8,-2.6);
			\draw [thick,color=gray] (3.5,-2.55) to [out=10,in=170] (6.8,-2.4) to [out=-10,in=-90] (8.5,0);
			\draw [thick,color=gray] (6.1,-2) to [out=110,in=0] (5,-1.2) to [out=180,in=100]  (3.1,-2.4) to [out=-80,in=180] (5,-5.2) to [out=0,in=-90] (6.2,-2.75);
			\draw[->,thick,color=gray] (4.9, -1.2) to (5.1, -1.2);
			
            \node[circle,draw,fill,scale=.3] at (8.5,0) {};
			\draw[->,thick] (8.5, 0) to (8.5, 1);
            \node at (1,-5) {$\AuxContr{G,K,tw}$};
		\end{scope}
	\end{tikzpicture}
	\caption{Illustration for the Jaeger test.}
	\label{fig:Jaeger_elimination_example}
\end{figure}
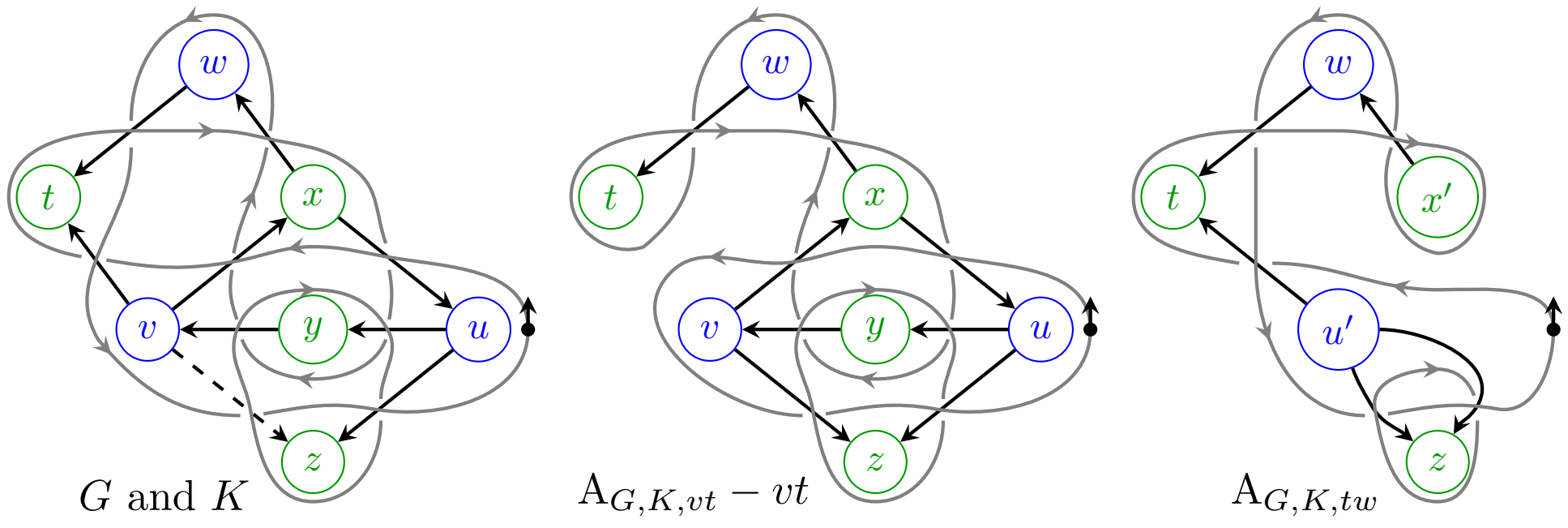

\begin{remark}\label{rem:reason_for_contraction}
    We will see later (in Section \ref{s:formula_proof}) that if $G$ is a Seifert graph of some link diagram, and $K$ is a Jaeger quasi-tree of $K$, then each $\AuxContr{G,K,e}$ is also a Seifert graph of a modified diagram. For most steps of the Jaeger process, this is clear. The only case that is not obvious is when we obtain $\AuxContr{G,K,f}$ as $\AuxContr{G,K,e}/C$ due to case \ref{e_in_K_tail_active}. This contraction operation corresponds to ``pulling apart'' two parts of the link that have a trivial relationship to each other, where this triviality is ensured by the fact that $C$ is a bad separating cycle in $\AuxContr{G,K,e}-e$. (See the third panel of Figure \ref{fig:Jaeger_elimination_example} for an example.) Hence,  $\Link(\AuxContr{G,K,f})=\Link(\AuxContr{G,K,e}/C)$ is isotopic to $\Link(\AuxContr{G,K,e})$. This means that the contraction simplifies the diagram without modifying the actual link.
\end{remark}

\begin{ex}\label{ex:Jaeger_test} 
The left panel of Figure \ref{fig:Jaeger_elimination_example} shows a ribbon digraph $G$, where the non-dashed edges $K=G-vz$ form a Jaeger quasi-tree for the root $(u,ux)$. Note that $G$ is the Seifert graph of a complicated diagram of the negative Hopf link. We have drawn the link in grey, in part to help visualize the ribbon structure.

The first step of the Jaeger test considers $ux$, where $(u,ux) <_K (x,ux)$ and $u$ is the head of $ux$, hence for the next edge $xv$, we get $\AuxContr{G,K,xv}=G$, and similarly, $\AuxContr{G,K,vt}=G$. The first interesting step is the subsequent edge $vt$ where $(v,vt) <_K (t,vt)$ and $v$ is the tail. The second panel of Figure \ref{fig:Jaeger_elimination_example} shows $\AuxContr{G,K,vt}-vt=G-vt$. We can see that the cycle $\{ux,xv, vy, yu\}$ forms a nearest bad separating cycle in $G-vt$, whence case \ref{e_in_K_tail_active} applies. We get $\AuxContr{G,K,tw}=\AuxContr{G,K,vt}/\{ux,xv, vy, yu\}$, depicted on the third panel.
Continuing with the process, we see that each edge passes the Jaeger test. 
\end{ex}

\begin{remark}
    Intuitively, we require two things for a quasi-tree $K$ to be a Jaeger quasi-tree: the first one is that each non-edge is first seen at its tail in the tour of $K$. This is a necessary ``side-effect'' of our strategy of applying skein relations on crossings that we reach on the understrand. 
    
    The second condition is that if for the currently examined edge $e$, that we reach at its tail, the graph $\AuxContr{G,K,e}-e$ has a bad separating cycle, and if the nearestsuch bad separating cycle is $C$ then $C\subseteq K$.  
    As we mentioned in Remark \ref{rem:reason_for_contraction}, in this case contracting $C$ results in a Seifert graph with an isotopic link, and less crossings. Hence it makes sense to pass to the Seifert graph $\AuxContr{G,K,e}/C$  
    to avoid unnecessary branches in the computation tree.

    We note that requiring the second condition for Jaeger quasi-trees is not only a computational simplification. In Section \ref{ssec:naive_Jaeger_trees_are_not_all_trivial} we show that our main theorem would not remain true if we used all naive Jaeger quasi-trees instead of just Jaeger quasi-trees.
\end{remark}

\begin{defn}[tight edges of a quasi-tree]\label{def:type0_type1_and_2_edges_in_K}
We call an edge $e\in K$ in a quasi-tree $K$ \emph{tight} if $K$ passes the Jaeger test for $e$ via case \ref{e_in_K_head}, \ref{e_in_K_tail_active} or \ref{e_in_K_contracted}. In other words, $e\in K$ is not reversed in the auxiliary graph during the Jaeger test. We call $e$ a \emph{type~$0$} edge if $K$ passes the Jaeger test for $e$ via case \ref{e_in_K_head}, we call $e$ a
\emph{type~$1$} edge if $K$ passes the Jaeger test for $e$ via case \ref{e_in_K_tail_active}, and we say that $e$ is of \emph{type~$2$} if $K$ passes the Jaeger test for $e$ via case \ref{e_in_K_contracted}.
\end{defn}

\begin{remark}
    In Section \ref{s:activities}, we will show that the tight edges of a Jaeger quasi-tree $K$ are in fact exactly the internally active edges of $K$. Hence, along the branch of the Homfly computation tree leading to $K$, we do not perform skein relations on these edges. This is why these edges do not get reversed in the auxiliary graph.
\end{remark}

The following observation will be important.

\begin{prop}\label{prop:type1_edge_becomes_bridge}
    Let $G$ be a rooted ribbon digraph and $K$ a quasi-tree of $G$ that passes the Jaeger test for each edge $f$ such that $f\leq_K e$. Here we suppose that $e\in K$ a type~$1$ edge of $K$. Then for each $f>_K e$ such that $\AuxContr{G,K,f}$ is defined, the edge $e$ is a bridge in $\AuxContr{G,K,f}$.
\end{prop}

\begin{proof}
Since $e$ is a type~$1$ edge of $K$, the graph $\AuxContr{G,K,e}-e$ has a bad separating cycle. Let $C$ be the nearestsuch cycle. Let $f$ be the edge following $e$ in $<_K$. Then $(\AuxContr{G,K,e}-e)/C=\AuxContr{G,K,f}-e$ is not connected. On the other hand $\AuxContr{G,K,f}$ is connected because its subgraph corresponding to $K$ is a quasi-tree. Indeed, all the contracted edges are part of $K$, also none of the deleted edges belong to $K$, whence the boundary of $K$ is not affected by our operations.
Therefore $e$ is a bridge in $\AuxContr{G,K,f}$. It is easy to check that $e$ being a bridge is preserved in all subsequent auxiliary graphs, too.
\end{proof}

Let us introduce some more notation. 

\begin{defn}
Let $K$ be a quasi-tree of $G$ that passes the Jaeger test until $e$.
We define $K^{\AuxContr{G,K,e}}$ 
as the subgraph corresponding to $K$ in $\AuxContr{G,K,e}$. That is, if some edge got reversed, we reverse it, and if some edge got contracted, we contract it. 
\end{defn}

Notice that if $K$ is a Jaeger quasi-tree, then $K^{\AuxContr{G,K,e}}$ is a quasi-tree in $\AuxContr{G,K,e}$. 
This is true since the contracted edges of $\AuxContr{G,K,e}$ are all part of
$K$, whereas its deleted edges are not; that is, the tour of $K^{\AuxContr{G,K,e}}$ is a restriction of the tour of $K$ to half-edges of $\AuxContr{G,K,e}$.

The following is immediate from the definitions.

\begin{prop}\label{prop:Jaeger_is_naive_Jaeger}
    Jaeger quasi-trees are also naive Jaeger quasi-trees. 
\end{prop}

\begin{coroll}\label{cor:no_Jaeger_if_bad_sep}
    If $G$ has a bad separating cycle, then $G$ has no Jaeger quasi-trees.
\end{coroll}

\begin{proof}
    Follows directly from Proposition \ref{prop:Jaeger_is_naive_Jaeger} and Proposition \ref{prop:no_naive_Jaeger_quasi-tree_if_separating_directed_cut}.
\end{proof}

Crucially, the converse is also true.

\begin{thm} \label{thm:if_no_oriented_sep_cycle_then_exists_Jaeger}
    If a connected rooted ribbon digraph $G$ has no bad separating cycles
    (i.e., there is not a separating directed cycle $C$ in $G$ such that the cut $C^*$ of $G^*$ is oriented toward the base component), 
    then $G$ contains a Jaeger quasi-tree. 
\end{thm}

Before giving the proof, let us make a few remarks. Theorem \ref{thm:if_no_oriented_sep_cycle_then_exists_Jaeger} will be an essential tool in our arguments. For instance, it will enable us to show that any link has a diagram whose Seifert graph contains a Jaeger quasi-tree (with any root; see Theorem \ref{thm:each_link_has_diagram_with_Jaeger_quasi-trees}). At this point, let us remark the following corollary.

\begin{coroll}\label{cor:acyclic_has_Jaeger_quasitrees}
    If $G$ is an acyclic ribbon digraph, then for any root $(v_0.e_0)$ we have $\jaeg(G,v_0,e_0)\neq \varnothing$, that is, $G$ has Jaeger quasi-trees for any root.
\end{coroll}
\begin{proof}
    Since any bad separating cycle is 
    directed, $G$ has no bad separating cycles with respect to $(v_0,e_0)$. Now Theorem \ref{thm:if_no_oriented_sep_cycle_then_exists_Jaeger} implies that $\jaeg(G,v_0,e_0)\neq\varnothing$.
\end{proof}

As to the proof that follows:
It is easy to show that the conditions of Theorem \ref{thm:if_no_oriented_sep_cycle_then_exists_Jaeger} imply the existence of a naive Jaeger quasi-tree. Indeed, if $G$ has no bad separating cycles, then in $G^*$ each vertex is reachable from the root. Hence $G^*$ has an arborescence, and it is easy to see that the complement of that is a naive Jaeger quasi-tree of $G$ (which \emph{keeps} the highest possible number of edges among quasi-trees).
However, these naive Jaeger quasi-trees are not always (proper) Jaeger quasi-trees, and it is not so easy to see when they are Jaeger. Hence we employ a different, in a way opposite, strategy and
construct the quasi-tree greedily, always \emph{deleting} edges unless we are forced not to.
The following lemma will be necessary.

\begin{lemma}\label{lem:contractability_of_bad_sep_cycle}     Let $G$ be a connected rooted ribbon digraph without a bad separating cycle. Let $e$ be an edge of $G$. 
    If $G-e$ has a bad separating cycle, and its nearest bad separating cycle is $C$, then $G/C$ does not have any bad separating cycles and $e$ is a bridge in $G/C$.
\end{lemma}

\begin{proof}
    We claim that $C$ is not a separating cycle in $G$. If $C$ were a separating cycle in $G$, then $C^*$ would be a cut in $G^*$. Then, since $e^*\notin C^*$, the set $C^*$ would still be a cut in $(G-e)^*=G^*/e^*$. Since in $G^*/e^*$, each edge of $C^*$ points toward the component containing the root, this is also true for $C^*$ in $G^*$. Hence $C$ would be a bad separating cycle in $G$, which is a contradiction.
    
    Since $C$ is a separating cycle in $G-e$, the graph $(G-e)/C=G/C-e$ is disconnected. On the other hand, since $C$ is not a separating cycle in $G$, the graph $G/C$ is connected. Hence $e$ is indeed a bridge in $G/C$.

    Let $D_0$ be the connected component of $(G-e)^*-C^*$ containing the root, and let $D_1$ be the other connected component. As $C$ is a bad separating cycle in $G-e$, each edge of $C^*$ has its head in $D_0$.

    Suppose for a contradiction that $G/C$ has a bad separating cycle $C'$. Then, $(C')^*$ is an elementary cut in $(G/C)^*=G^*-C^*$.
    Since $e$ is a bridge in $G/C$, we have $e\notin C'$. Thus, $(C')^*$ is still a cut in $(G^*-C^*)/e^*=(G-e)^*-C^*$. Recall that $(G-e)^*-C^*$ has two connected components $D_0$ and $D_1$. Hence $(C')^*$ either has a subset $(C'')^*\subseteq (C')^*$ that is a cut within $D_0$, or a subset that is a cut within $D_1$, or both.
    
    From this point on, we look at two cases separately.

    Case 1: If $(C')^*$ has a subset $(C'')^*$ that is a cut within $D_0$. Let $D'_0$ be the connected component of 
    $D_0-(C'')^*$ containing the root. By the assumption that $C'$ is a bad separating cycle in $G/C$, we know that each edge $f^*\in (C'')^*$ has its head in $D'_0$. Recall that $D'_0$ is a subgraph of $(G-e)^*-C^*$, wherefore we can also think of it as a subgraph of $(G-e)^*$.
    Take the edges of $(G-e)^*$ connecting $V(D'_0)$ to $V((G-e)^*)-V(D'_0)$, and let this edge set be $S^*$. We have $S^*\subseteq (C'')^*\cup C^*$, but in both cases, the head of each edge is in $D'_0$. That is, $S^*$ is a cut of $(G-e)^*$ directed toward $D'_0$. Hence $S$ is a bad separating cycle in $G-e$, and because $V(D'_0) \subsetneq V(D_0)$, it is strictly smaller than $C$. This contradicts the assumption that $C$ was the nearest bad separating cycle in $G-e$.

    Case 2: 
    If $(C')^*$ has a subset $(C'')^*$ that is a cut within $D_1$.
    Then, by undoing the contraction of $e^*$ in $(G^*-C^*)/e^*$, we glue a vertex of $D_0$ to a vertex of $D_1$, whence we get that $(C'')^*$ is also a cut within $G^*-C^*=(G/C)^*$, moreover, one of the shores of $(C'')^*$ entirely contains $D_0$. 
    Let $D'_0$ and $D'_1$ be the shores of the cut $(C'')^*$ within $(G/C)^*$, with $D'_0 \supseteq D_0$. By the assumption that $C''$ is a bad separating cycle in $G/C$, each edge of $(C'')^*$ points toward $D'_0$. 
    
    Recall that $D'_0$ and $D'_1$ are subgraphs of $(G/C)^*=G^*-C^*$, wherefore we can also think of them as subgraphs of $G^*$.
    Let $S^*$ be the set of edges of $G^*$ connecting $D'_0$ to $D'_1$. These are either edges of $(C'')^*$ or of $C^*$, but in both cases, the head of the edge is in $D'_0$. (If $f^*\in C^*\cap S^*$ then this is true because the endpoint of $f^*$ in $D'_1 \subseteq D_1$ needs to be the tail.) Hence, $S$ is a bad separating cycle in $G$, which is again a contradiction.
\end{proof}

\begin{proof}[Proof of Theorem \ref{thm:if_no_oriented_sep_cycle_then_exists_Jaeger}]
We seek to build up a Jaeger quasi-tree $K$ edge-by-edge. In the beginning, for each edge $e\in E(G)$, it is undecided whether it is in $K$ or not.
Notice, however, that no matter what the future $K$ will be, we have $\AuxContr{G,K,e_{0}}=G$, 
and the first half-edge of $<_K$ is $(v_0,e_0)$. Our first aim is to decide whether $e_0\in K$ or $e_0\notin K$ in such a way that $K$ passes the Jaeger test for $e_0$. Notice that whether $K$ passes the Jaeger test for $e_0$ depends only on $\AuxContr{G,K,e_0}=G$, which is well-defined at this point, and it is independent of the rest of $K$. Moreover, notice that once the status of whether $e_0\in K$ or $e_0\notin K$ is decided, the next half-edge $(v_1,e_1)$ of $<_K$ is well-defined regardless of the rest of $K$. (Of course, unless $e_1=e_0$, it is yet undefined whether $e_1\in K$ or not.) 

We inductively continue this reasoning:  
In the general step, take the first half-edge $(v,e)$ in the 
already existing portion of the
tour of $K$ (in other words, according to $<_K$) such that whether $e\in K$ or $e\notin K$ has not been decided. It follows that the other half-edge $(w,e)$ of $e$ has not yet occurred in the tour, since in that case, the status whether $e$ is in $K$ or not would have already been decided, and the half-edge following $(v,e)$ would also be well-defined. 
In any case, $<_K$ is well-defined up to $e$, which implies that $\AuxContr{G,K,e}$ is also well-defined.

Suppose that in previous steps, we made choices such that $K$ passed the Jaeger test for the edges preceding $e$ in $<_K$, and $\AuxContr{G,K,e}$ has no bad separating cycles. (This is true by assumption for the base case $e=e_0$.) We show that we can choose the status of $e$ such that $K$ also passes the Jaeger test for $e$, and $\AuxContr{G,K,e}$ still has no bad separating cycles.

If $e$ is contracted in $\AuxContr{G,K,e}$, then let $e\in K$. Then, regardless of how $K$ is built up for the rest of the edges, the edge $e$ will pass the Jaeger test with case \ref{e_in_K_contracted}, and for the next edge $f$ of the tour of $K$ we have $\AuxContr{G,K,f}=\AuxContr{G,K,e}$, which does not contain a bad separating cycle. 

If $e$ is not contracted in $\AuxContr{G,K,e}$ and $v=h_G(e)$, then $(h_G(e),e) <_K (t_G(e),e)$. Choose $e\in K$, and (no matter how we choose the rest of the edges of $K$) $\AuxContr{G,K,f}=\AuxContr{G,K,e}$. 
With this choice, $K$ passes the Jaeger test for $e$ (with case \ref{e_in_K_head}), and $\AuxContr{G,K,f}$ has no bad separating cycles.

If $e$ is not contracted in $\AuxContr{G,K,e}$ and $v=t_G(e)$, then we have $(t_G(e),e) <_K (h_G(e),e)$.
If $\AuxContr{G,K,e}-{e}$ has no bad separating cycle, then put $e\notin K$. In this case, regardless of how the rest of $K$ is chosen, we have $\AuxContr{G,K,f}=\AuxContr{G,K,e}-e$.
Hence, $K$ passes the Jaeger test for $e$ (with case \ref{e_notin_K_tail}), and $\AuxContr{G,K,e}$ has no bad separating cycle. 

If $\AuxContr{G,K,e}-{e}$ has a bad separating cycle, then 
put $e\in K$. Regardless of how the rest of $K$ is chosen, we are in case \ref{e_in_K_tail_active}. Hence $\AuxContr{G,K,f}=\AuxContr{G,K,e}/C$ 
for the nearest bad separating cycle $C$ of $\AuxContr{G,K,e}-e$ (which is also independent of how the rest of $K$ is chosen). $K$ passes the Jaeger test for $e$ with case \ref{e_in_K_tail_active}, and by Lemma \ref{lem:contractability_of_bad_sep_cycle}, $\AuxContr{G,K,f}=\AuxContr{G,K,e}/C$ has no bad separating cycles.

We have thus proved that we can choose the status of $e$ such that $K$ passes the Jaeger test for $e$, and for the edge $f$ following $e$ in $<_K$, the graph $\AuxContr{G,K,f}$ does not have a bad separating cycle.

It remains to show that using the above strategy, we eventually reach each half-edge, in other words, $K$ is a quasi-tree. 
Suppose for a contradiction that $K$ is not a quasi-tree, that is, we arrive back to the half-edge $(v_0,e_0)$ without seeing each half-edge of $G$. (It is clear that the first half-edge that we see for the second time needs to be $(e_0,v_0)$ since the tour of $K$ is also deterministic backward.)

Take 
the set of edges $K$ and 
the auxiliary graph $\AuxContr{G,K,\infty}$ that we have when we arrive back to $(v_0,e_0)$. 
In the tour of $K$, when we reach a contracted edge, it is always an edge of $K$.
Hence, the tour of $K^{\AuxContr{G,K,\infty}}$ in $\AuxContr{G,K,\infty}$ is a restriction of the tour of $K$ in $G$ to edges of $\AuxContr{G,K,\infty}$ (where the orientations of some edges might differ in the two graphs). 
If we contract some edges in a ribbon graph, the resulting graph has the same boundary components as the original (with the caveat that some components may end up surrounding isolated vertices).
Hence if $K$ is not a quasi-tree of $G$, then also $K^{\AuxContr{G,K,\infty}}$ is not a quasi-tree of $\AuxContr{G,K,\infty}$. 

Suppose for a contradiction that $K^{\AuxContr{G,K,\infty}}$ is not a quasi-tree in $\AuxContr{G,K,\infty}$.  
Take those edges of $\AuxContr{G,K,\infty}$ where one half-edge was traversed in the tour of $K^{\AuxContr{G,K,\infty}}$ and the other one was not. As $\AuxContr{G,K,\infty}$ is connected (indeed, a disconnected graph has $\varnothing$ as a bad separating cycle), this is a nonempty set of edges, and it forms a separating cycle $C$ of $\AuxContr{G,K,\infty}$. If some edge $e\in C$ is reached at its tail in the tour of $K^{\AuxContr{G,K,\infty}}$, then $e$ can only be a type~$1$ edge of $K$ (considered in case \ref{e_in_K_tail_active} of the Jaeger test). 
Now by Proposition \ref{prop:type1_edge_becomes_bridge}, $e$ is a bridge in $\AuxContr{G,K,\infty}$. Since a cycle cannot contain a bridge, $C$ does not contain any type~$1$ edges. This means that $C$ is a bad separating cycle in $\AuxContr{G,K,\infty}$, a contradiction.
\end{proof}

\section{Activities in quasi-trees}\label{s:activities}

In this section we give the formal definitions of the activities that occur in Theorem \ref{thm:Homfly_formula}. We give a more general definition, where $\mathcal{F}$ is a family of quasi-trees, and the activities of a quasi-tree $K\in \mathcal{F}$ are defined with respect to $\mathcal{F}$. For most of the paper, we will use this definition with $\mathcal{F}$ being the set of Jaeger quasi-trees of a ribbon digraph with respect to some fixed root.

\begin{defn}[internal passivity and activity]\label{def:internal_activity_of_quasi-trees}
Let $G$ be a rooted ribbon graph, and $\mathcal{F}$ a family of quasi-trees of $G$. 

For a quasi-tree $K\in \mathcal{F}$, an edge $e\in K$ is called \emph{internally passive} in $K$ (with respect to $\mathcal{F}$) if it is the first difference between $K$ and another quasi-tree $K'\in \mathcal{F}$. That is, if there exists another quasi-tree $K'\in\mathcal{F}$ such that $e\notin K'$ and the tours of $K$ and $K'$ agree until $e$.

We call an edge $e\in K$ \emph{internally active} in $K$ (with respect to $\mathcal{F}$) if it is not internally passive.
We denote the number of internally active edges of $K$ by $ia_{\mathcal{F}}(K)$ and the number of internally passive edges if $K$ by $ip_{\mathcal{F}}(K)$.

We denote the number of internally active positive edges of $K$ by $ia^+_{\mathcal{F}}(K)$, the number of internally active negative edges of $K$ by $ia^-_{\mathcal{F}}(K)$, the number of internally passive positive edges of $K$ by $ip^+_{\mathcal{F}}(K)$, and the number of internally passive negative edges of $K$ by $ip^-_{\mathcal{F}}(K)$. 
\end{defn}

\begin{defn}[external passivity and activity]\label{def:external_activity_of_quasi-trees}
Let $G$ be a rooted ribbon graph, and $\mathcal{F}$ a family of quasi-trees of $G$. 

For a quasi-tree $K\in \mathcal{F}$, an edge $e\notin K$ is called \emph{externally passive} in $K$ (with respect to $\mathcal{F}$) if it is the first difference between $K$ and another Jaeger quasi-tree $K'\in \mathcal{F}$. That is, if there exists another quasi-tree $K'\in\mathcal{F}$ such that $e\in K'$ and the tours of $K$ and $K'$ agree until $e$.

We call an edge $e\notin K$ \emph{externally active} in $K$ (with respect to $\mathcal{F}$) if it is not externally passive. We denote the number of externally active edges of $K$ by $ia_{\mathcal{F}}(K)$ and the number of externally passive edges if $K$ by $ip_{\mathcal{F}}(K)$.

We denote the number of externally active positive edges of $K$ by $ea^+_{\mathcal{F}}(K)$, the number of externally active negative edges of $K$ by $ea^-_{\mathcal{F}}(K)$, the number of externally passive positive edges of $K$ by $ep^+_{\mathcal{F}}(K)$, and the number of externally passive negative edges of $K$ by $ep^-_{\mathcal{F}}(K)$. 
\end{defn}

In this paper, we will almost always talk about activities with respect to Jaeger quasi-trees, that is, $\mathcal{F}=\jaeg(G)$ (where we use the same root for defining the Jaeger quasi-trees as for determining the first differences of Jaeger quasi-trees).
However, we will occasionally talk about activities with respect to Jaeger trees (i.e., Jaeger quasi-trees that are trees).

Our goal in this section is to give structural characterizations for activities of Jaeger quasi-trees, using bad separating cycles in the auxiliary graphs. 

The following claim states a ``local'' property of the Jaeger test, and will be very important.

\begin{claim}\label{cl:Jaeger_test_is_local}
    Let $G$ be a connected rooted ribbon digraph and let $K$ and $L$ be quasi-trees of $G$. If the tours of $K$ and $L$ agree until $e$, moreover $K$ passes the Jaeger test for all edges $f<_K e$, then $L$ also passes the Jaeger test for these same edges, described by the equivalent condition $f <_L e$, and for each $f\leq_K e$, we have 
    $\AuxContr{G,K,f}=\AuxContr{G,L,f}$.
\end{claim}

\begin{proof}
    The statement follows by induction over the set of edges, as ordered by $<_K$. $\AuxContr{G,K,e_{0}}=\AuxContr{G,L,e_{0}}=G$ 
    by definition. In a subsequent step regarding the edge $e$, 
    for $f<_K e$, the graphs $\AuxContr{G,K,f}$ and $\AuxContr{G,L,f}$ are isomorphic 
    by the inductive hypothesis. Since the tours of $K$ and $L$ agree until $e$, each of these edges $f$ is such that $f\notin K\Delta L$. Hence we perform identical steps on the isomorphic graphs $\AuxContr{G,K,f}$ and $\AuxContr{G,L,f}$ and 
    we choose the next edge with the same rule for $<_K$ and $<_L$.
\end{proof}

The next two lemmas compare Jaeger quasi-trees of $G$ and of the auxiliary digraphs. 
They essentially say that the Jaeger quasi-trees of 
$\AuxContr{G,K,e}$ are the same as those Jaeger quasi-trees of $G$ whose tours agree with that of $K$ up to $e$.

\begin{lemma}\label{lem:Jaeger_in_G_and_in_B}
	Let $G$ be a connected rooted ribbon digraph, and $e$ be an edge of $G$, or let $e=\infty$.     
    If $K$ is a Jaeger quasi-tree of $G$, then $K^{\AuxContr{G,K,e}}$ is a Jaeger quasi-tree of $\AuxContr{G,K,e}$. 
\end{lemma}

\begin{proof}
    Suppose that $K$ is a quasi-tree of $G$ that passes the Jaeger test (in $G$) for each edge $f<_K e$. We show inductively that $K^{\AuxContr{G,K,e}}$ passes the Jaeger test in $\AuxContr{G,K,e}$ for each edge $f\in \AuxContr{G,K,e}$ with $f <_{K^{\AuxContr{G,K,e}}} e$, and for each $f\leq_{K^{\AuxContr{G,K,e}}} e$, the auxiliary graph $\AuxContr{\AuxContr{G,K,e},K^{\AuxContr{G,K,e}},f}$ is $\AuxContr{G,K,e}$. 
    
    Let $(v'_0, e'_0)$ be the first half-edge of $G$ according to $<_K$ such that $e'_0\in K$ and $e'_0$ is not contracted in $\AuxContr{G,K,e}$. Then, $(v'_0, e'_0)$ is the root half-edge of $\AuxContr{G,K,e}$, and $e'_0$ is the first edge in $<_{K^{\AuxContr{G,K,e}}}$. By definition, $\AuxContr{\AuxContr{G,K,e},K^{\AuxContr{G,K,e}},e'_{0}}=\AuxContr{G,K,e}$ so the induction hypothesis holds for $f=e'_0$. 

    Now suppose that for some $f<_{K^{\AuxContr{G,K,e}}} e$, we have $\AuxContr{\AuxContr{G,K,e},K^{\AuxContr{G,K,e}},f}=\AuxContr{G,K,e}$. 
    We show that in this case $K^{\AuxContr{G,K,e}}$ passes the Jaeger test for $f$, and for the edge $g$ following $f$ in the tour of $K^{\AuxContr{G,K,e}}$, we have $\AuxContr{\AuxContr{G,K,e},K^{\AuxContr{G,K,e}},g}=\AuxContr{G,K,e}$. 
    
    We know that $K$ passes the Jaeger test for $f$. 
    That cannot happen via
    case \ref{e_notin_K_tail} because then $f$ would not be present in $\AuxContr{G,K,e}$ (it would have been deleted).
    
    Similarly, if $K$ passed the test for $f$ via case \ref{e_in_K_contracted}, then $f$ would not be present in $\AuxContr{G,K,e}$ (it would have been contracted), wherefore this is not possible, either.
    
    If $K$ passes the test for $f$ via case \ref{e_in_K_head}, then $f$ is either present in $\AuxContr{G,K,e}$ with the same orientation as in $G$ or, potentially, it has been contracted in some step after $f$. In the latter case we have $f\notin\AuxContr{G,K,e}$, which is a contradiction. In the former case $K^{\AuxContr{G,K,e}}$ passes the Jaeger test for $f$ in $\AuxContr{G,K,e}$ via case \ref{e_in_K_head}, and the auxiliary graph is not modified. Thus, the inductive assumption stays true for $g$.

    If $K$ passes the test for $f$ via case \ref{e_in_K_tail_active} then, by Proposition \ref{prop:type1_edge_becomes_bridge}, $f$ is a bridge in subsequent auxiliary graphs; in particular, $f$ is a bridge in $\AuxContr{G,K,e}$. Hence $\AuxContr{\AuxContr{G,K,e},K^{\AuxContr{G,K,e}},f}-f=\AuxContr{G,K,e}-f$ has a bad separating cycle, namely $\varnothing$. We claim that $\varnothing$ is also the nearest bad separating cycle in $\AuxContr{G,K,e}-f$. Indeed, if the nearestseparating cycle were some other subgraph $C$, then $C$ would need to be a separating cycle of the component of $\AuxContr{G,K,e}-f$ containing the root. In this case $C$ would still be a bad separating cycle in $\AuxContr{G,K,e}$, which is not possible. We conclude that case \ref{e_in_K_tail_active} applies in the Jaeger test of $K^{\AuxContr{G,K,e}}$ for $f$, and $\AuxContr{\AuxContr{G,K,e},K^{\AuxContr{G,K,e}},g}=\AuxContr{G,K,e}/\varnothing=\AuxContr{G,K,e}$. 
	
    Finally, if $K$ passes the test for $f$ via case \ref{e_in_K_tail_passive} then $f$ gets reversed and, possibly, even contracted in a later step. But if that happens then $f$ is not present in $\AuxContr{G,K,e}$, which contradicts our assumption. Else $K^{\AuxContr{G,K,e}}$ passes the Jaeger test for $f$ in $\AuxContr{G,K,e}$ via case \ref{e_in_K_head}. The auxiliary graph is not modified, so the induction hypothesis holds for the next edge $g$ in $<_{K^{\AuxContr{G,K,e}}}$.

    By the above argument, $K^{\AuxContr{G,K,e}}$ passes the Jaeger test in $\AuxContr{G,K,e}$ for all $f$ such that $f<_{K^{\AuxContr{G,K,e}}} e$, moreover $\AuxContr{\AuxContr{G,K,e},K^{\AuxContr{G,K,e}},e}=\AuxContr{G,K,e}$.
    
    Note that at this point we have the same auxiliary graph when performing the Jaeger test of $K$ for $e$ (in $G$) as when performing the Jaeger test of $K^{\AuxContr{G,K,e}}$ for $e$ (in $\AuxContr{G,K,e}$). 
    From this point on, the Jaeger tests for $K$ and for $K^{\AuxContr{G,K,e}}$ run identically.  
    Hence if $K$ is a Jaeger quasi-tree of $G$, then indeed $K^{\AuxContr{G,K,e}}$ is a Jaeger quasi-tree of $\AuxContr{G,K,e}$.
\end{proof}

\begin{defn}
    Let $G'$ be the reorientation of a minor of $G$. (For example, $G'=\AuxContr{G,K,e}$ for some edge $e$ and Jaeger quasi-tree $K$ of $G$.) Suppose that $G'$ is obtained from $G$ by contracting the edges $S_c\subseteq E(G)$, deleting the edges  $S_d\subseteq E(G)$, and reversing the edges $S_r \subseteq E(G)$.
    
    Let $J$ be a quasi-tree of $G'$. We define $J^{G}$ to be the subgraph $J^{G}=J \cup S_c$ of $G$. 
    That is, we take the edges of $J$ (that might have different orientations in $G$ and in $G'$) together with the contracted edges $S_c$. It is easy to see from the definition of contraction that $J^G$ is a quasi-tree of $G$.
\end{defn}

\begin{lemma}\label{lem:Jaeger_of_B_is_active_and_Jaeger} 
    Let $G$ be a connected rooted ribbon digraph with an edge $e$ and quasi-tree $K$ that passes the Jaeger test for each edge $f<_K e$, and let $L'$ be a Jaeger quasi-tree of $\AuxContr{G,K,e}$. 
\begin{enumerate}
    \item[$(i)$]  
    $L:=(L')^G$ is a Jaeger quasi-tree of $G$, and $L'=L^{\AuxContr{G,L,e}}$. 

     \item[$(ii)$] We have $f\in L'$ for each $f\in\AuxContr{G,K,e}$ such that $f<_{L'} e$, and the tours of $K$ and $L$ agree up to $e$. 

    \item[$(iii)$] Each edge $f\in\AuxContr{G,L,e}$ with $f<_{L'} e$ is internally active in $L'$. An edge $g\in\AuxContr{G,L,e}$ with $g \geq_{L'} e$ is internally/externally active/passive in $L'$ if and only if $g$ is internally/externally active/passive in $L$. If an edge $g \geq_L e$ is contracted in $\AuxContr{G,L,e}$, then $g$ is internally active in $L$, as well as in each Jaeger quasi-tree whose tour agrees with the tour of $L$ up to $e$.
\end{enumerate}
\end{lemma}

\begin{proof}
    $(i)$: We prove by induction that for each $f\leq_{L'} e$, we have $\AuxContr{\AuxContr{G,K,e},L',f}=\AuxContr{G,K,e}$, and for each $g<_{L'} f$ we have $g\in L'$ and $g\in K^{\AuxContr{G,K,e}}$. 
    
    Let $(v',e')$ be the first half-edge of $G$ according to $<_K$ such that $e'\in \AuxContr{G,K,e}$ (that is, it does not get deleted or contracted). Then, $(v',e')$ is the root of $\AuxContr{G,K,e}$, and $e'$ is the first edge in both $<_{K^{\AuxContr{G,K,e}}}$ and $<_{L'}$. By definition, $\AuxContr{\AuxContr{G,K,e},L',e'}=\AuxContr{G,K,e}$ so the assumption holds for $f=e'$. Suppose that the assumption holds for some $f <_{L'} e$.
    We show that in this case $f\in L' \cap K_{\AuxContr{G,K,e}}$ and, moreover, the edge $f'$ following $f$ in $<_{L'}$ satisfies $\AuxContr{\AuxContr{G,K,e},L',f'}=\AuxContr{G,K,e}$.

    Let $f=uw$ and suppose that $(u,f) <_{L'} (w,f)$. Since, by the induction hypothesis, the tours of $L'$ and $K^{\AuxContr{G,K,e}}$ agree until $f$, we also have $(u,f) <_{K^{\AuxContr{G,K,e}}} (w,f)$, and $f <_{K^{\AuxContr{G,K,e}}} e$, which imply $f <_K e$. This means that $f$ passes the Jaeger test for $K$. 
    As $f\in \AuxContr{G,K,e}$, it passes the Jaeger test with one of the cases \ref{e_in_K_head}, \ref{e_in_K_tail_active} or \ref{e_in_K_tail_passive}. In each case we conclude $f\in K^{\AuxContr{G,K,e}}$.
    
    If $f$ passes the test for $K$ via case \ref{e_in_K_head} or \ref{e_in_K_tail_passive}, then we have $u=h_{\AuxContr{G,K,e}}(f)$, which means that for $L'$ to pass the Jaeger test for $f$, we need $f\in L'$. Then, in $\AuxContr{G,K,e}$, the edge $f$ passes the Jaeger test of $L'$ with case \ref{e_in_K_head}, which tells us that $\AuxContr{\AuxContr{G,K,e},L',f'}=\AuxContr{\AuxContr{G,K,e},L',f}=\AuxContr{G,K,e}$.

    If $f$ passes the Jaeger test for $K$ with case \ref{e_in_K_tail_active}, then for the edge $f''$ following $f$ in $<_K$, we see that $f$ becomes a bridge in $\AuxContr{G,K,f''}$ and this property stays true in subsequent auxiliary graphs, too. In particular, $f$ is a bridge in $\AuxContr{G,K,e}$. Because $K^{\AuxContr{G,K,e}}$ and $L'$ are quasi-trees in $\AuxContr{G,K,e}$, we again have to have $f\in L'\cap K^{\AuxContr{G,K,e}}$ and $f$ passes the Jaeger test for $L'$ via case \ref{e_in_K_tail_active} with $C=\varnothing$. Again, $\AuxContr{\AuxContr{G,K,e},L',f''}=\AuxContr{\AuxContr{G,K,e},L',f}=\AuxContr{G,K,e}$. 
    
    We have proved that for each $f\leq_{L'} e$, we have $\AuxContr{\AuxContr{G,K,e},L',f}=\AuxContr{G,K,e}$, and for each $g<_{L'} f$ we have $g\in L'$ and $g\in K^{\AuxContr{G,K,e}}$. 

    $(ii)$:  
    We claim that $L:=(L')^G$ is a quasi-tree of $G$.
    We can derive the half-edge order $<_{L'}$ by removing from $<_L$ the half-edges corresponding to contracted edges. Hence, each half-edge of $L'$ is reached in the tour of $L$. If $L$ is not a quasi-tree, then some contracted edges of $\AuxContr{G,K,e}$ form a separating cycle $C_L$ in $L$. 
    
    Take the last edge $f$ according to $<_K$ such that $(L')^{\AuxContr{G,K,f}}$ is not a quasi-tree. (Such an edge exists since $(L')=(L')^{\AuxContr{G,K,e}}$ is a quasi-tree, while $L=(L')^{G}=(L')^{\AuxContr{G,K,e_0}}$ is not.) Let $f'$ be the edge immediately following $f$ in $<_K$. Then, to obtain $\AuxContr{G,K,f'}$ from $\AuxContr{G,K,f}$, we need to contract a bad separating cycle $C$ of $\AuxContr{G,K,f}-f$. As $\AuxContr{G,K,f}$ has no bad separating cycles, $C$ is not a separating cycle in $\AuxContr{G,K,f}$. On the other hand, in order for $(L')^{\AuxContr{G,K,f'}}$ to be a quasi-tree, we need to contract each edge of $C_L$, which implies $C_L\subseteq C$. This however, as a separating cycle cannot be a subset of a non-separating cycle, is a contradiction. With this, we have proved that $L=(L')^G$ is a quasi-tree of $G$.

    Note that to get $L$ from $L'$, we need to undo contractions on the same set of edges, and also to reverse the same set of edges, as when we reconstruct $K$ from $K^{\AuxContr{G,K,e}}$. Since the tours of $L'$ and $K^{\AuxContr{G,K,e}}$ agree until $e$, the tours of $L$ and $K$ also agree until $e$. As $K$ passes the Jaeger test until $e$, by Claim \ref{cl:Jaeger_test_is_local}, $L$ also passes the Jaeger test for each $f <_L e$.

    Since we have now established $\AuxContr{\AuxContr{G,K,e},L',e}=\AuxContr{G,K,e}=\AuxContr{G,L,e}$, we see that the Jaeger property of $L'$ implies that $L$ is also a Jaeger quasi-tree. Indeed, over the remaining edges, the Jaeger tests run identically for the two quasi-trees. 
    
    $(iii)$: By Remark \ref{rem:osszehuzott_kor_bekerül_a_kvazifaba}, if some edge $g \geq_L e$ is contracted in $\AuxContr{G,L,e}$, then $g\in L$. If $N$ is a Jaeger quasi-tree of $G$ such that the tours of $L$ and $N$ agree until $e$, then by Claim \ref{cl:Jaeger_test_is_local}, we have  $\AuxContr{G,L,e}=\AuxContr{G,N,e}$. This implies $g\in N$ for any Jaeger quasi-tree of $G$ whose tour agrees with the tour of $L$ until $e$. In particular, for any Jaeger quasi-tree whose tour agrees with $L$ until $e$, the edge $g$ cannot be a first difference. We conclude that $g$ is indeed internally active in any such Jaeger quasi-tree. 
    
    It follows immediately from part $(ii)$ that each edge $f<_{L'} e$ is internally active in $L'$. If some edge $g \geq_{L'} e$ is internally (resp., externally) passive in $L'$, then there exists another Jaeger quasi-tree $L''$ of $\AuxContr{G,K,e}$ such that the tours of $L'$ and $L''$ agree until $g$ and $g\in L'-L''$ (resp., $g\in L''-L'$). In this case, by part $(i)$, $(L')^G$ and $(L'')^G$ are Jaeger quasi-trees of $G$ whose tours agree until $g$, and $g\in (L')^G-(L'')^G$ (resp., $g\in (L'')^G-(L')^G$), which implies that $g$ is indeed internally (resp., externally) active in $G$.
    
    If $g$ is internally (externally) passive in $L$, then from our previous argument we see that $g\in\AuxContr{G,L,e}$.
    We can then use Lemma \ref{lem:Jaeger_in_G_and_in_B} to conclude that $g$ is internally (externally) passive in $L'$.
\end{proof}

After these preparations, we are now in a position to give graph-theoretic characterizations for internally and externally active edges of Jaeger quasi-trees. Tight edges (of various types) were introduced in Definition \ref{def:type0_type1_and_2_edges_in_K}.

\begin{thm}\label{thm:internal_activity_condition}
	Let $G$ be a connected rooted ribbon digraph, $K$ a Jaeger quasi-tree of $G$, and $e\in K$.  
	The following are equivalent.
	\begin{itemize}
		\item[$(i)$] $e$ is internally active in $K$,
		\item[$(ii)$] 
        $e$ is a tight edge of $K$ (with respect to $G$). 	\end{itemize}
\end{thm}

\begin{proof}
	$(ii)\Rightarrow (i)$:
    Suppose for a contradiction that $(ii)$ holds but $(i)$ does not. That is, $e\in K$ is a tight edge, but $G$ has a Jaeger quasi-tree $L$ with $e\notin L$ such that the tours of $K$ and $L$ agree until $e$. 
    Then, by Claim \ref{cl:Jaeger_test_is_local}, we have $\AuxContr{G,K,e}=\AuxContr{G,L,e}$. 
        
    As the tours of $K$ and $L$ agree until $e$, if $e$ is of type~$0$, that is, we have $(h(e),e)<_{K} (t(e),e)$, then we also have $(h(e),e) <_L (t(e),e)$, but that is not allowed for a non-edge of a Jaeger quasi-tree (case \ref{e_notin_K_head} or \ref{e_notin_K_contracted}).

	If $e$ is of type~$1$, then $\AuxContr{G,K,e}-e=\AuxContr{G,L,e}-e$ has a bad separating cycle. Let $f$ be the edge immediately following $e$ in $<_L$. Since $e\notin L$, we have $\AuxContr{G,L,f}=\AuxContr{G,L,e}-e$. Now by Corollary \ref{cor:no_Jaeger_if_bad_sep}, the graph $\AuxContr{G,L,f}$ has no Jaeger quasi-trees. This contradicts Lemma \ref{lem:Jaeger_in_G_and_in_B}, which claims that $L^{\AuxContr{G,L,f}}$ is a Jaeger quasi-tree of $\AuxContr{G,L,f}$.
    
    If $e$ is a type~$2$ edge of $K$, then $e$ is contracted in $\AuxContr{G,K,e}=\AuxContr{G,L,e}$, which implies that we are in case 
    \ref{e_notin_K_contracted} for $e$ in the Jaeger test of $L$. This again contradicts the assumption that $L$ is a Jaeger quasi-tree of $G$.
    
	$(i)\Rightarrow (ii)$:
    Suppose that $(ii)$ does not hold. Then $e$ is not contracted in $\AuxContr{G,K,e}$, moreover, $(t_G(e),e) <_K (h_G(e),e)$ and $\AuxContr{G,K,e}-e$ has no bad separating cycles. 
    
    We claim that there exists some quasi-tree $K'$ in $G$ such that $e\notin K'$ and the tours of $K$ and $K'$ agree until $e$. (At this point, we do not yet claim $K'$ to be a Jaeger quasi-tree.)
    Consider $K-e$. To get the desired quasi-tree $K'$, it is enough to find an edge $f$ such that $f >_K e$, and either $f\in K$ and $K-e-f$ is a quasi-tree, or $f\notin K$ and $K-e+f$ is a quasi-tree. If there is an edge $f>_K e$ such that the two half-edges of $f$ belong to different boundary components of $K-e$, then $f$ fits our purpose.
    We show that if there is no such edge $f>_K e$, then $\AuxContr{G,K,e}-e$ has a bad separating cycle. The edges $f\in K-e$, such that the two half-edges of $f$ belong to different components of $K-e$, form a (possibly empty) separating cycle $C$ of $K$. By our assumption,
    each edge $g\in C$ has $g<_K e$.  
    
    Let $C'$ be the subgraph of $\AuxContr{G,K,e}$ corresponding to $C$, that is, obtained by contracting those edges of $C$ that are contracted in $\AuxContr{G,K,e}$. Each edge $g<_K e$ such that $g\notin K$ is deleted from $\AuxContr{G,K,e}-e$. 
    Hence $C'$ is a separating cycle in $\AuxContr{G,K,e}-e$. Each edge $g\in C$ satisfies $g<_K e$, whence $K$ passes the Jaeger test for each $g\in C$. If some edge $g\in K$ passed the Jaeger test via case \ref{e_in_K_tail_active}, then $g$ would be a bridge in $\AuxContr{G,K,e}$, implying $g\notin C'$. If $g$ passed the Jaeger test with case \ref{e_in_K_contracted}, then $g$ is contracted in $\AuxContr{G,K,e}$, hence $g\notin C'$. Finally, if $g$ passed the Jaeger test via case \ref{e_in_K_head} or \ref{e_in_K_tail_passive}, then $(h_{\AuxContr{G,K,e}}(g),g) <_{K^{\AuxContr{G,K,e}}} (t_{\AuxContr{G,K,e}}(g),g)$ would follow. Therefore $C'$ is indeed a bad separating cycle in $\AuxContr{G,K,e}-e$, a contradiction. 

    We conclude that there exists a quasi-tree $K'$ such that the tours of $K$ and $K'$ agree until $e$, but $e\notin K'$. Let $f'$ be the edge following $e$ in $<_{K'}$. 
    By Claim \ref{cl:Jaeger_test_is_local} we have $\AuxContr{G,K',e}=\AuxContr{G,K,e}$, 
    and $K'$ passes the Jaeger test for each edge $g<_{K'} e$. By assumption we have $(t_G(e),e) <_K (h_G(e),e)$, whence $(t_G(e),e) <_{K'} (h_G(e),e)$ and $e\notin K'$. Hence $K'$ also passes the Jaeger test for $e$, and $\AuxContr{G,K',f'}=\AuxContr{G,K',e}-e$. 
    
    Since $\AuxContr{G,K',f'}=\AuxContr{G,K,e}-e$ has no bad separating cycles, by Theorem \ref{thm:if_no_oriented_sep_cycle_then_exists_Jaeger}, $\AuxContr{G,K',f'}$ has a Jaeger quasi-tree $J$. By Lemma \ref{lem:Jaeger_of_B_is_active_and_Jaeger}, in this case $L:=J^G$ is a Jaeger quasi-tree of $G$, moreover, the tours of $L$ and $K'$ agree for all $g<_L f'$. 
    As the tours of $K$ and $K'$ agree for $g <_K e$, we see that the tours of $L$ and $K$ also agree for $g<_K e$. 
    Now from $e\in K-L$ we conclude that $(i)$ does not hold.
\end{proof}

\begin{thm}\label{thm:ext_activity_sufficient}
	Let $G$ be a connected rooted ribbon digraph, let $K$ be a Jaeger quasi-tree of $G$, and let $e\notin K$.
	The following are equivalent: 
    \begin{itemize}
    	\item[$(i)$] $e$ is externally active in $K$,
    	\item[$(ii)$] 
    	$(t_G(e),e)<_K (h_G(e),e)$ and the graph $G'$, obtained from $\AuxContr{G,K,e}$ by reversing $e$, has a bad separating cycle.
    \end{itemize}
\end{thm}

\begin{proof}[Proof of Theorem \ref{thm:ext_activity_sufficient}] 	$(ii)\Rightarrow (i)$:
    Suppose for a contradiction that $(ii)$ holds but $(i)$ does not. That is, $G$ has a Jaeger quasi-tree $L$ with $e\in L$ such that the tours of $K$ and $L$ agree until $e$.
    Then, by Claim \ref{cl:Jaeger_test_is_local}, we have $\AuxContr{G,K,e}=\AuxContr{G,L,e}$.
    Since $e\notin K$ and $K$ is Jaeger, we have $(t_G(e),e) <_K (h_G(e),e)$. Since the tours of $K$ and $L$ agree until $e$, we also have $(t_G(e),e) <_L (h_G(e),e)$.

    Since $e$ is a first difference between $K$ and $L$ and $e\in L$, by definition $e$ is not internally active in $L$. Hence by Theorem \ref{thm:internal_activity_condition}, $e$ is not a tight edge of $L$ (with respect to $G$). In other words, $e$ needs to pass the Jaeger test of $L$ with case \ref{e_in_K_tail_passive}. In particular, for the edge $f$ following $e$ in $<_L$ we have $\AuxContr{G,L,f}=G'$. 
    By Corollary \ref{cor:no_Jaeger_if_bad_sep} the graph $G'$ has no Jaeger quasi-trees, on the other hand by Lemma \ref{lem:Jaeger_in_G_and_in_B}, the subgraph $L^{\AuxContr{G,K,f}}$ is a Jaeger quasi-tree in $G'=\AuxContr{G,L,f}$, which is a contradiction.

    $(i)\Rightarrow (ii)$:
    Suppose that $(ii)$ does not hold.  
    We claim that there exists some quasi-tree $K'$ in $G$ such that $e\in K'$ and the tours of $K$ and $K'$ agree until $e$. (At this point, we do not yet claim $K'$ to be a Jaeger quasi-tree.)
    Consider the extended subgraph $K+e$, which has two boundary components. To get the desired quasi-tree $K'$, it suffices to find an edge $f$ such that $f >_K e$, and either $f\in K$ and $K+e-f$ is a quasi-tree, or $f\notin K$ and $K+e+f$ is a quasi-tree. Any edge $f>_K e$ such that the two half-edges of $f$ belong to different boundary components of $K+e$ satisfies this requirement.
    
    If there is no edge $f>_K e$ such that the two half-edges of $f$ belong to different boundary components of $K+e$, then $K+e$ has a (possibly empty) separating cycle $C$ consisting of edges whose two half-edges are in different boundary components of $K+e$, and each edge $g\in C$ has $g\leq_K e$. But this would mean that there is a bad separating cycle in $G'$. Indeed, let $C'$ be the set of edges of $C$ that are not contracted in $G'$. Since each edge $g<_K e$ that has $g\notin K$ is deleted from  $\AuxContr{G,K,e}$, the set $C'$ is a separating cycle in $\AuxContr{G,K,e}$ and thus also in $G'$. As each edge $g\in C$ has $g\leq_K e$, each edge $g\in C-e$ passes the Jaeger test for $K$. If some edge $g\in K$
    passed the Jaeger test via case \ref{e_in_K_tail_active}, then $g$ would be a bridge in $\AuxContr{G,K,e}$, implying $g\notin C'$. If $g$ passed the Jaeger test via case \ref{e_in_K_contracted}, then it would be contracted in $\AuxContr{G,K,e}$, whence again $g\notin C'$. If $g$ passed the Jaeger test via case \ref{e_in_K_head} or \ref{e_in_K_tail_passive}, then $(h_{\AuxContr{G,K,e}}(g),g) <_{K^{\AuxContr{G,K,e}}} (t_{\AuxContr{G,K,e}}(g),g)$ would hold. Finally, by the assumption that $(t_G(e),e) <_K (h_G(e),e)$, we have $(h_{G'}(e),e)<_K (t_{G'}(e),e)$.
    Hence $C'$ is indeed a bad separating cycle in $G'$, a contradiction. 

    We conclude that there exists a quasi-tree $K'$ such that the tours of $K$ and $K'$ agree until $e$, but $e\in K'$. Let $f$ be the edge following $e$ in $<_{K}$ and let $f'$ be the edge following $e$ in $K'$. 
    We claim that $G'=(\AuxContr{G,K,e})_{\overleftarrow{e}}=\AuxContr{G,K',f'}$ and $K'$ passes the Jaeger test for each edge $g<_{K'} f'$. To see this, note that by Claim \ref{cl:Jaeger_test_is_local} we have $\AuxContr{G,K',e}=\AuxContr{G,K,e}$ and that $K'$ passes the Jaeger test for all $g<_{K'} e$. 

    Lemma tells us that \ref{lem:Jaeger_in_G_and_in_B} $K^{\AuxContr{G,K,f}}$ is a Jaeger quasi-tree in $\AuxContr{G,K,f}=\AuxContr{G,K,e}-e$. Hence $\AuxContr{G,K,e}-e$ contains no bad separating cycles. 
    By assumption we have $(t_G(e),e) <_K (h_G(e),e)$, 
    whence $(t_G(e),e) <_{K'} (h_G(e),e)$. Therefore $K'$  passes the Jaeger test for $e$ as well, and $\AuxContr{G,K',f'}=(\AuxContr{G,K,e})_{\overleftarrow{e}}=G'$.  

    Since $G'$ has no bad separating cycles, Theorem \ref{thm:if_no_oriented_sep_cycle_then_exists_Jaeger} implies that $G'$ has a Jaeger quasi-tree $L'$. Let $L=(L')^{G}$.
    By Lemma \ref{lem:Jaeger_of_B_is_active_and_Jaeger} we have that $L:=(L')^G$ is a Jaeger quasi-tree of $G$, moreover, the tours of $L$ and $K'$ agree for all $g<_L f'$. 
    As the tours of $K$ and $K'$ agree for $g <_{K} e$, the tours of $L$ and $K$ also agree for $g<_K e$. 
    As $e\in K-L$, we conclude that $(i)$ does not hold.
\end{proof}

We note that the proof of part $(ii)\Rightarrow (i)$ implies the following alternative theorem, which is rather useful:
\begin{lemma}\label{lem:orthogonality}
    Let $G$ be a connected rooted ribbon digraph without bad separating cycles, let $K\in \jaeg(G)$, and $e$ an edge of $G$ such that $(t(e),e) <_K (h(e),e)$. Then the following two things cannot both happen:
    \begin{itemize}
        \item $\AuxContr{G,K,e}-e$ has a bad separating cycle,
        \item $(\AuxContr{G,K,e})_{\overleftarrow{e}}$ has a bad separating cycle.
    \end{itemize}
\end{lemma}
We also note that this result is a special case of the Farkas lemma of orthogonal matroids by Baker, Ding and Kim \cite[Lemma 2.12.]{BDK}.

On the computational side, we can say the following.
\begin{prop}\label{prop:Jaeger_and_activity_computation_in_polynomial_time}
    Given a connected rooted ribbon digraph and a quasi-tree $K$ of $G$, it can be decided in polynomial time (in $|E(G)|$) whether $K$ is a Jaeger quasi-tree, and if yes, the sets of internally/externally active/passive edges of $K$ can be determined.
\end{prop}

\begin{proof}
    During the Jaeger test for $K$, the only nontrivial decisions occur when an edge $e\in K$ is such that $(t_G(e),e) <_K (h_G(e),e)$, in which case one has to decide whether $\AuxContr{G,K,e}-e$ contains a bad separating cycle, and if so, determine the nearest bad separating cycle. This can be achieved in polynomial time by Propsition \ref{prop:bad_sep_cycle_computation}.
    We conlude that the Jaeger test can be run in polynomial time in $|E(G)|$.

    To determine the activities, by Theorems \ref{thm:internal_activity_condition} and \ref{thm:ext_activity_sufficient}, we need to check if an edge $e\in K$ is tight (this can be done in no extra time during the Jaeger test), and if for an edge $e\notin K$ the graph $(\AuxContr{G,K,e})_{\overleftarrow{e}}$ has a bad separating cycle. Again, this can be done by Proposition \ref{prop:bad_sep_cycle_computation}.
\end{proof}

\section{Jaeger quasi-trees in knot theory}\label{s:formula_proof}

In this section we prove our main theorem.
For a link diagram $D$, let us denote by $\Seifert(D)$ the Seifert graph of $D$. Seifert graphs are treated as signed, directed ribbon graphs which also happen to be bipartite.

\begin{customthm}{\ref{thm:Homfly_formula}} 
Let $D$ be a diagram of the oriented link $L$, and let $G=\Seifert(D)$. Let $(v_0,e_0)$ be an arbitrary half-edge of $G$ such that $\jaeg(G,v_0,e_0)\neq\varnothing$. Then the Homfly polynomial of $L$ is
\begin{align*}
P(v,z)=&\sum_{K\in\jaeg(G,v_0,e_0)} (v^2)^{ip^+(K)}\cdot (v^{-2})^{ip^-(K)}\cdot (vz)^{ep^+(K)}\cdot (-v^{-1}z)^{ep^-(K)} \\ &\cdot (vz+vz^{-1}-v^3z^{-1})^{ea^+(K)}\cdot (-v^{-1}z-v^{-1}z^{-1}+v^{-3}z^{-1})^{ea^-(K)}.    
\end{align*}
\end{customthm}

The exponents in the formula were introduced in Definitions \ref{def:internal_activity_of_quasi-trees} and \ref{def:external_activity_of_quasi-trees}. 
In what follows, we will alternate between formally assigning Homfly polynomials to links and to their particular diagrams.

We prove the theorem at the end of the section. For that, we need a couple of lemmas. Before stating them, let us say a few words on the condition $\jaeg(G,v_0,e_0)\neq\varnothing$ for Seifert graphs.

\subsection{Existence of Jaeger quasi-trees for Seifert graphs of links}
\label{ss:treeexistence}

Luckily, we can easily modify any link diagram so that it represents the same link and for any basepoint, there is not a bad separating cycle in its Seifert graph.

\begin{lemma}\label{lem:link_has_diagram_with_connected_acyclic_seifert}
    Each link has a diagram whose Seifert graph is connected and acyclic.
\end{lemma}

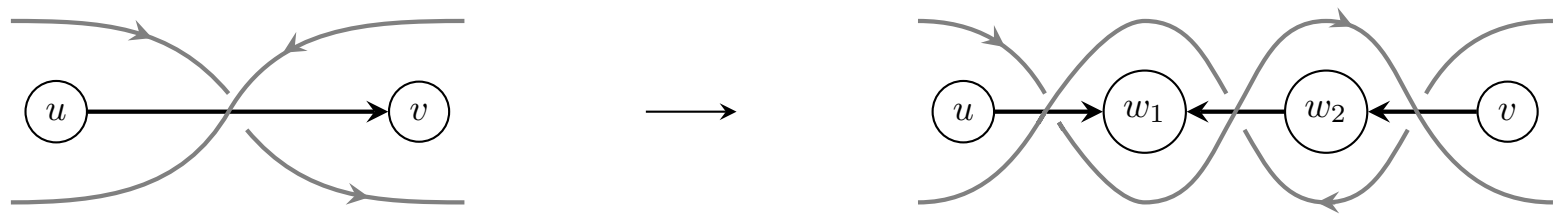
\begin{figure}
    \centering
    \begin{tikzpicture}[decoration={markings,
    mark= at position 0.6 with {\arrow{stealth}}}]
    \begin{scope}[scale=.6,shift={(0,0)},>=stealth]
    \tikzstyle{v}=[circle,scale=.7,draw]
    \node [v] (u) at (0,0) {$u$};
    \node [v] (v) at (4,0) {$v$};
    \draw [->,thick] (u) to (v);
    \draw [thick,postaction={decorate},color=gray] (-0.5,1) to [out=0,in=130] (1.9,0.2);    
    \draw [thick,postaction={decorate},color=gray] (2.1,-0.2) to [out=-50,in=180] (4.5,-1);
    \draw [thick,color=gray] (-0.5,-1) to [out=0,in=-120] (1.9,0) to [out=60,in=180] (4.5,1);
	\draw[->,thick,color=gray] (2.55, 0.68) to (2.5, 0.65);
    \draw[->] (6.5,0) to (7.5,0);
    \end{scope}
    \begin{scope}[scale=.6,shift={(10,0)},>=stealth]
    \tikzstyle{v}=[circle,scale=.7,draw]
    \node [v] (u) at (0,0) {$u$};
    \node [v] (v) at (6,0) {$v$};
    \node [v] (w1) at (2,0) {$w_1$};
    \node [v] (w2) at (4,0) {$w_2$};
    \draw [->,thick] (u) to (w1);
    \draw [->,thick] (w2) to (w1);
    \draw [->,thick] (v) to (w2);
    \draw [thick,postaction={decorate},color=gray] (-0.5,1) to [out=0,in=120] (0.9,0.2);    
    \draw [thick,postaction={decorate},color=gray] (1.1,-0.2) to [out=120,in=180] (2,-1) to [out=0,in=-120] (3,0) to [out=60,in=180] (4,1) to [out=0,in=120] (5,0) to [out=-60,in=180] (6.5,-1);
    \draw [thick,color=gray] (-0.5,-1) to [out=0,in=-120] (0.9,0) to [out=-120,in=180] (2,1) to [out=0,in=120] (2.9,0.2);
    \draw [thick,color=gray] 
    (3.1, -0.2) to [out=-60,in=180] (4,-1) to [out=0,in=-120] (4.9,-0.2);
    \draw [thick,color=gray] 
    (5.1, 0.2) to [out=60,in=180] (6.5,1);
	\draw[->,thick,color=gray] (4.1, -1) to (3.9, -1);
    \end{scope}
    \end{tikzpicture}
    \caption{Subdividing an edge of a Seifert graph using a Reidemeister II move.}
\label{fig:edge_subdivision_in_Seifert_graphs}
\end{figure}

\begin{proof}
    By applying several Reidemeister II moves it is easy to ensure that the diagram, and thus its Seifert graph, is connected.

    If $G$ has a directed cycle $C$, then take any edge $uv\in C$ (where, say, $u=t(uv)$), and substitute it with a path $uw_1$, $w_1w_2$, $w_2v$ such that $u=t(uw_1)$, $w_2=t(w_1w_2)$, and $v=t(w_2v)$. Then, $G'$ is the Seifert graph of a link diagram that is related to the original one via a Reidemeister II move. (See Figure \ref{fig:edge_subdivision_in_Seifert_graphs}.) With this step, we eliminated a directed cycle, and the newly constructed edges $uw_1$, $w_1w_2$ and $w_2v$ do not occur in any directed cycle because $w_1$ and $w_2$ have undirected degree 2, and $w_1$ has outdegree 0.
\end{proof}

Now Corollary \ref{cor:acyclic_has_Jaeger_quasitrees} directly implies the following.

\begin{thm}\label{thm:each_link_has_diagram_with_Jaeger_quasi-trees}
Each link has a diagram such that its Seifert graph $G$ is connected and acyclic, and in that case $G$ contains a Jaeger quasi-tree with respect to any basepoint. 
\end{thm}

For some important classes of diagrams, the acyclicity of the Seifert graph is automatic.

\begin{prop}\label{prop:pos_diag_Seifert_acyclic}
The Seifert graph of a positive diagram is acyclic.
\end{prop}

\begin{proof}
It is easy to check that any type~I Seifert circle of a positive diagram is either a source or a sink of the Seifert graph. Type~II Seifert circles are cut vertices with only incoming edges on one side and only outgoing edges on the other side. Such a graph is clearly acyclic.
\end{proof}

\subsection{Proof of Theorem \ref{thm:Homfly_formula}}

Our first lemma is a key observation about the role of bad separating cycles. 

\begin{lemma}\label{lemma:bad_sep_cycle_implies_dstant_union}
Let $D$ be an oriented link diagram with a basepoint, and let  $G=\Seifert(D)$ (considered as a rooted ribbon digraph).
Suppose that $G$ has a bad separating cycle $C$, and let $\mathcal{C}_0$ and $\mathcal{C}_1$ be the collections of boundary components of $G$ that are separated by $C$.
Then $D$ represents the distant union of the link formed by the components of $\mathcal{C}_0$ and the link formed by the components of $\mathcal{C}_1$.
\end{lemma}

\begin{ex}
    Consider the link diagram and its Seifert graph $G$ in the left panel of Figure \ref{fig:duality_example}. The basepoint is on the red component, next to the vertex $u$. Then, $\{ux,xv,vz,zu\}$ is a bad separating cycle in $G$ with $\mathcal{C}_0$ being formed by the red and the black components, while $\mathcal{C}_1$ contains just the gray component. Notice that indeed the link is the distant union of the gray knot and the link consisting of the red and black components.
\end{ex}

\begin{proof}
	The fact that $C$ is a separating cycle in $G$ means that its dual $C^*$ is an elementary cut in $G^*$. The boundary components contained in $\mathcal{C}_0$ and $\mathcal{C}_1$ are the vertices of $G^*$, as separated to the two shores of the cut $C^*$. 
    Assume that the basepoint lies along an element of $\mathcal C_0$.
    Then the fact that $C$ is a ``bad'' separating cycle means that each edge of $C^*$ points from $\mathcal{C}_1$ to $\mathcal{C}_0$. Equivalently, each crossing between a component of $\mathcal{C}_0$ and a component of $\mathcal{C}_1$ corresponds to some element of $C$, and the strand of the component of $\mathcal{C}_0$ crosses above. 
	 
	This indeed means that $D$ depicts the distant union of the link formed by the components of $\mathcal{C}_0$ and the link formed by the components of $\mathcal{C}_1$. 
	Indeed, after a suitable parallel translation in $\R^3$ orthogonal to the projection plane $\pi$, the two links can be separated by a plane parallel to $\pi$.
\end{proof}

The next lemma formalizes the phenomenon that was depicted in Figure \ref{fig:shorthand}.

\begin{lemma}\label{l:handling_edges_that_need_to_be_removed}
Let $D$ be an oriented link diagram with a basepoint, and let  $G=\Seifert(D)$ (considered as a rooted ribbon digraph).
Suppose that $K$ is a Jaeger quasi-tree of $G$ and $e\notin K$ is an externally active edge.
Let $f$ be the edge immediately following $e$ in $<_K$. Suppose that $\AuxContr{G,K,e}$ is the Seifert graph of a link diagram $D_{G,K,e}$. Then, $\AuxContr{G,K,f}$ is also the Seifert graph of a link diagram $D_{G,K,f}$. Moreover, if $e$ is a positive edge, then we have
$$P_{D_{G,K,e}}=(vz+vz^{-1}-v^3z^{-1})\cdot P_{D_{G,K,f}},$$
whereas if $e$ is a negative edge, then 
$$P_{D_{G,K,e}}=(-v^{-1}z-v^{-1}z^{-1}+v^{-3}z^{-1})\cdot P_{D_{G,K,f}}.$$
\end{lemma}

\begin{proof}
    If $e$ is externally active with respect to $K$, then $\AuxContr{G,K,f}=\AuxContr{G,K,e}-e$. Thus $\AuxContr{G,K,f}$ is the Seifert graph of the diagram that we get from $D_{G,K,e}$ by smoothing the crossing corresponding to $e$.

	Let $G'$ be the graph obtained from $\AuxContr{G,K,e}$ by reversing $e$. By Theorem \ref{thm:ext_activity_sufficient}, 
    $G'$ has a bad separating cycle $C$. 
    On the other hand, $\AuxContr{G,K,e}$ does not have any bad separating cycles, since by Lemma \ref{lem:Jaeger_in_G_and_in_B}, $K^{\AuxContr{G,K,e}}$ is a Jaeger quasi-tree in $\AuxContr{G,K,e}$, and by Corollary \ref{cor:no_Jaeger_if_bad_sep}, a rooted ribbon digraph with a bad separating cycle can not contain a Jaeger quasi-tree.
    As the orientation of $e$ is the only difference between $\AuxContr{G,K,e}$ and $G'$, it has to be that $e\in C$.

	Apply the skein relation to $D_{G,K,e}$ at $e$. 
    Changing the crossing results in a diagram $D'$ such that $\Seifert(D')=(\AuxContr{G,K,e})_{\overleftarrow{e}}=G'$, whereas smoothing it leads to the diagram $D_{G,K,f}$ that has $\Seifert(D_{G,K,f})=\AuxContr{G,K,e}-e=\AuxContr{G,K,f}$. Hence if $e$ is a positive edge,
    we get \begin{equation}\label{eq:external_activity_skein_positive}
    P_{D_{G,K,e}}=v^2\cdot P_{D'} + vz\cdot P_{D_{G,K,f}},  
    \end{equation}
    while if $e$ is a negative edge, we obtain  \begin{equation}\label{eq:external_activity_skein_negative}
    P_{D_{G,K,e}}=v^{-2}\cdot P_{D'} - v^{-1}z\cdot P_{D_{G,K,f}}.    
    \end{equation}

    The components of $D'$ correspond to the boundary components of $G'$, or in other words, to the vertices of the dual $(G')^*$. 
    Since $C$ is a bad separating cycle of $G'=\Seifert(D')$, by Lemma \ref{lemma:bad_sep_cycle_implies_dstant_union}, $D'$
    is the distant union of the link with components $L^0_0, \dots, L^0_k$ and the link with components $L^1_0, \dots, L^1_\ell$, where  $L^0_0, \dots, L^0_k$ are the boundary components of $G'$ on the on the side of $C$ containing the root, and $L^1_0, \dots, L^1_\ell$ are the boundary components on the other side of $C$.
    
    On the other hand, to get $D_{G,K,f}$ from $D'$, we pinch together two components $L^0_i$ and $L^1_j$. This is because, as we noted above, $e\in C$. 
    Now by a standard and straightforward application of the skein relation, this implies $P_{D'}=(v^{-1}z^{-1}-vz^{-1})P_{D_{G,K,f}}$. 
    Substituting this in \eqref{eq:external_activity_skein_positive} or \eqref{eq:external_activity_skein_negative}, depending on whether $e$ is a positive or a negative edge, we get the statement of the Lemma.
\end{proof}

The next lemma describes the cases when we decide not to apply the skein relation at a certain crossing. Its conclusion is fairly obvious from a knot theory standpoint, but care must be taken to ensure that the graph theory and the knot theory match up as expected.

\begin{lemma}\label{l:induction_internally_active_edges} 
Let $D$ be an oriented link diagram with a basepoint, and let  $G=\Seifert(D)$ (considered as a rooted ribbon digraph).
Suppose that $K$ is a Jaeger quasi-tree of $G$ and $e\in K$ is an internally active edge.
Let $f$ be the edge immediately following $e$ in $<_K$. Suppose that $\AuxContr{G,K,e}$ is the Seifert graph of a link diagram $D_{G,K,e}$. 
Then $\AuxContr{G,K,f}$ is also the Seifert graph of a link diagram $D_{G,K,f}$, moreover, $D_{G,K,e}$ and $D_{G,K,f}$ are diagrams of the same link.
\end{lemma}

\begin{proof}
    By Theorem \ref{thm:internal_activity_condition} we know that either $(h_G(e),e) <_K (t_G(e),e)$, or $e$ is contracted in $\AuxContr{G,K,e}$, or $\AuxContr{G,K,e}-e$ has a bad separating cycle.
    In the first two cases we have
    $\AuxContr{G,K,f}=\AuxContr{G,K,e}$ by definition, thus there is nothing to prove.

    If $\AuxContr{G,K,e}-e$ has a bad separating cycle, then for the nearest bad separating cycle $C$ of $\AuxContr{G,K,e}-e$, the next auxiliary graph $\AuxContr{G,K,f}$ is defined as $\AuxContr{G,K,e}/C$. It is enough to show that $\AuxContr{G,K,e}/C$ is the Seifert graph of a link diagram isotopic 
    to $D_{G,K,e}$.

    Take the diagram $\tilde{D}$ obtained from $D_{G,K,e}$ by smoothing the crossing corresponding to $e$. (Note that we have $\Seifert(\tilde{D})=\AuxContr{G,K,e}-e$.)

    As $C$ is a separating cycle in $\AuxContr{G,K,e}-e$, its dual $C^*$ is an elementary cut in $(\AuxContr{G,K,e}-e)^*$. Thus, the vertices of $(\AuxContr{G,K,e}-e)^*$ (equivalently, the components of $\tilde{D}$) can be divided into two sets $\mathcal{C}_0$ and $\mathcal{C}_1$ (with $\mathcal{C}_0$ containing the root component). 
    By Lemma \ref{lemma:bad_sep_cycle_implies_dstant_union}, in this case, $\tilde{D}$ is the distant union of the link formed by the components of $\mathcal{C}_0$ and the link formed by the components of $\mathcal{C}_1$.
    
    Let us return to the diagram $D_{G,K,e}$.
    We can view it as a collection of knots corresponding to elements of $\mathcal{C}_0$ and $\mathcal{C}_1$, where one element of $\mathcal{C}_0$ and one of $\mathcal{C}_1$ are joined by the crossing $e$. (In $\R^3$, the two strands forming the crossing intersect the separating plane of the previous paragraph.)
    By applying stereographic projection (which, mind, does not affect the Seifert graphs!), we may assume that $\mathcal{C}_0$ and $\mathcal{C}_1$ are drawn such that the crossing $e$ is incident to the unbounded face of each.
    
    Now, apply 
    what is essentially a $180^\circ$ degree rotation, about the crossing $e$ (more precisely, the line through $e$ that is perpendicular to $\pi$),
    to the components of $\mathcal{C}_1$. 
    On the one hand this is an isotopy of links. On the other hand, at the level of the diagram, we kept each crossing where both strands were from $\mathcal{C}_0$ or both strands were from $\mathcal{C}_1$, but we eliminated those crossings where one strand was from $\mathcal{C}_0$ and the other strand was from $\mathcal{C}_1$. That is, we obtained a link whose dual Seifert graph is exactly $(\AuxContr{G,K,e})^*-C^*$. In other words, the Seifert graph of this link is $\AuxContr{G,K,e}/C$. Hence we proved that this is indeed the Seifert graph of a link diagram isotopic to $D_{G,K,e}$, finishing the proof.
\end{proof}

The point of the rotation of the last paragraph of the proof is to remove $\mathcal C_1$ from underneath $\mathcal C_0$, which we can do because the existence of $C$ guarantees that there is no obstruction to it. We have already performed this kind of operation, for instance between the first and third panels of Figure \ref{fig:Jaeger_elimination_example}.
The bottom panel $G$ of Figure \ref{fig:computationtree} can be another example, with the role of $e$ being played by the crossing in the middle. By smoothing it we obtain two unknots $\mathcal C_0$ and $\mathcal C_1$, and in this case we may just shrink, in $G$, the part corresponding to $\mathcal C_1$ to eliminate the four crossings that form $C^*$. This works because $e$ is incident to the unbounded region of $\mathcal C_1$. In general though, it is possible for both $\mathcal C_0$ and $\mathcal C_1$ to surround $e$, which causes the slight complication in describing the necessary transformation.

Our last lemma clarifies why each end stage of our computation is indeed an unknot diagram.

\begin{lemma}\label{l:Jaeger_quasi-trees_link_trivial} 
    Let $D$ be an oriented link diagram with a basepoint, and let  $G=\Seifert(D)$ (considered as a rooted ribbon digraph).
    If $K\in \jaeg(G)$ is a Jaeger quasi-tree of $G$, then for each edge $e\in E(G)$, and also for $e=\infty$, the graph $\AuxContr{G,K,e}$ is the Seifert graph of a link diagram $D_{G,K,e}$. Furthermore, $D_{G,K,\infty}$ is a diagram of a trivial knot. 
\end{lemma}

\begin{proof}
    Let us use the notation $e_0 <_K e_1 <_K \dots <_K e_m$ for the edges of $G$. Run the Jaeger test of $K$, resulting in the sequence of auxiliary graphs $G=\AuxContr{G,K,e_0}$, $\AuxContr{G,K,e_1}, \dots,  \AuxContr{G,K,\infty}$. 
    By assumption, $G=\AuxContr{G,K,e_0}$ is a Seifert graph. Suppose that $i$ is such that $\AuxContr{G,K,e_j}=\Seifert(D_{G,K,e_j})$ is a Seifert graph for all $j\leq i$. If $e_i$ is an internally or externally active edge of $K$, then by Lemmas \ref{l:induction_internally_active_edges} and \ref{l:handling_edges_that_need_to_be_removed}, we see that $\AuxContr{G,K,e_{i+1}}=\Seifert(D_{G,K,e_{i+1}})$ is also a Seifert graph. If $e_i\notin K$ is an externally passive edge with respect to $K$, then we have $\AuxContr{G,K,e_{i+1}}=\AuxContr{G,K,e_i}-e_i$, which is the Seifert graph of the diagram obtained from $D_{G,K,e_i}$ by smoothing the crossing corresponding to $e_i$. If $e_i\in K$ is an internally passive edge in $K$, then Theorem \ref{thm:internal_activity_condition} implies that $e_i$ passes the Jaeger test via case \ref{e_in_K_tail_passive}, whence $\AuxContr{G,K,e_{i+1}}=(\AuxContr{G,K,e_i})_{\overleftarrow{e_i}}$. In this case $D_{G,K,e_{i+1}}$ can be obtained from $D_{G,K,e_i}$ by switching the under- and overstrands at the crossing corresponding to $e_i$. 
    With that the first half of the statement follows
    by induction, in particular we see that $\AuxContr{G,K,\infty}$ is indeed the Seifert graph of a link diagram $D_{G,K,\infty}$.
    
    Traverse $D_{G,K,\infty}$ starting from the basepoint. This equivalently means following the tour of $K^{\AuxContr{G,K,\infty}}=\AuxContr{G,K,\infty}$ 
    since all the non-edges of $K$ are already removed from $\AuxContr{G,K,\infty}$.
    Because $K$ is a quasi-tree, $D_{G,K,\infty}$ is a diagram of a knot.
    If $(h_{\AuxContr{G,K,\infty}}(e),e) <_{K^{\AuxContr{G,K,\infty}}} (t_{\AuxContr{G,K,\infty}}(e),e)$ for an edge $e\in K$, then for the crossing corresponding to $e$, the first-traversed strand of $D_{G,K,\infty}$ will cross above the second strand. If by chance,  
    the above is true for each edge $e\in \AuxContr{G,K,\infty}$,
    then $D_{G,K,\infty}$ is a descending diagram
    and as such it does represent the trivial knot.
    If there is an edge $e\in K$ such that $(t_{\AuxContr{G,K,\infty}}(e),e) <_{K^{\AuxContr{G,K,\infty}}} (h_{\AuxContr{G,K,\infty}}(e),e)$ then, by the construction of $\AuxContr{G,K,\infty}$, the edge $e$ is a type~$1$ edge of $K$, which implies that $e$ is a bridge in $\AuxContr{G,K,\infty}$. 

    Let $F\subset K$ be the set of type~$1$ edges of $K$. Take some $f\in F$. Since $f$ is a bridge in $\AuxContr{G,K,\infty}$, for each edge $g\in \AuxContr{G,K,\infty}$, either both half-edges of $g$ are between $(t_{\AuxContr{G,K,\infty}}(f),f)$ and $(h_{\AuxContr{G,K,\infty}}(f),f)$ in the order $<_{K^{\AuxContr{G,K,\infty}}}$, or neither half-edge is. Hence, if we smooth the crossing corresponding to $f$, we get two distant links, that is, two smaller links whose diagrams do not cross each other. 

    The occurrences of the half-edges $(t_{\AuxContr{G,K,\infty}}(f),f)$ and $(h_{\AuxContr{G,K,\infty}}(f),f)$ for $f\in F$ determine intervals in $<_{K^{\AuxContr{G,K,\infty}}}$. By the above remark, these intervals are laminar. That is, for $f_1, f_2\in F$, either the intervals corresponding to $f_1$ and $f_2$ are disjoint, or one of them contains the other. Take an $f\in F$ such that the corresponding interval does not contain the interval of any other element of $F$. Then, upon smoothing the crossing corresponding to $f$, we get two disjoint diagrams and among these,
    the one corresponding to the part of the tour between $(t_{\AuxContr{G,K,\infty}}(f),f)$ and $(h_{\AuxContr{G,K,\infty}}(f),f)$ is descending. 
    Now by applying a $360^\circ$ rotation to this part of $D_{G,K,\infty}$ (about an axis \emph{parallel} to the projection plane),
    we may switch the under- and overstrands of the crossing corresponding to $f$ without changing the rest of the diagram or the underlying knot.
    Moreover, this new diagram $D'$ has $\Seifert(D')=:G'=(\AuxContr{G,K,\infty})_{\overleftarrow{f}}$.  In $G'$, the set of edges with $(t_{G'}(e),e) <_{K^{G'}} (h_{G'}(e),e)$ is $F'=F-f$, and it is still true that the edges of $F'$ are all bridges in $G'$. We can continue with this argument until we arrive at a descending diagram 
    that is equivalent to
    $D_{G,K,\infty}$. Thus indeed, $D_{G,K,\infty}$ represents the trivial knot.
\end{proof}

Finally, let us put the ingredients together and establish our main theorem.

\begin{proof}[Proof of Theorem \ref{thm:Homfly_formula}]
    Throughout this proof we will think of $G$ as a rooted ribbon digraph, with root $(v_0,e_0)$. 
 
    We will use the notion of 
    the ``first edge that is excluded from at least one Jaeger quasi-tree'' in $G$. 
    This edge $e$ is defined as follows.
    For each Jaeger quasi-tree $K$, the first edge according to $<_K$ is $e_0$. If there exists any Jaeger quasi-tree $K$ such that $e_0\notin K$, then we take $e=e_0$. If $e_0$ is included in every Jaeger quasi-tree, then for each $K\in\jaeg(G)$, the next edge in $<_K$ is $e_1=\next(w,e_0)$ where $e_0=v_0w$. We can continue similarly. Eventually, we either conclude that each Jaeger quasi-tree contains all the edges of $G$ (which implies that there is in fact a single Jaeger quasi-tree that agrees with $G$), in which case we define $e=\infty$, or we find an edge $e$ that is the smallest edge according to each $<_K$ ($K\in\jaeg(G)$) among $\{f\in E(G)\mid \exists K\in\jaeg(G) \text{ s.t.\ } f\notin K\}$. In particular, the tours of any two Jaeger quasi-trees agree until $e$. Let us denote by $t(G)$ the number of edges preceding $e$ in any order $<_K$. (For example, $t(G)=0$ if $e=e_0$, and $t(G)=|E(G)|$ if $e=\infty$.) 

    Suppose for a contradiction that the statement of the theorem is not true, and take a minimal counterexample in the following sense: let $D$ be a diagram with a basepoint such that for $G=\Seifert(D)$ and the induced basepoint $(v_0,e_0)$ there are Jaeger quasi-trees, and the formula of the Theorem does not hold. Moreover, among these counterexamples, take one such that  $|E(G)|-t(G)$ is minimal, and among those, $|E(G)|$ is minimal.
    Let $e\in E(G)\cup \{\infty\}$ be the first edge that is excluded from at least one Jaeger quasi-tree, as explained in the previous paragraph.

    If $e=\infty$ then each edge of $G$ is contained in each Jaeger quasi-tree, so there is a unique Jaeger quasi-tree $K$ and $G=K$. It also follows that all edges $f\in K=G$ are internally active in $K$. Indeed, by definition, an internally passive edge is a first difference between $K$ and another Jaeger quasi-tree of $G$, which is impossible if there are no more Jaeger quasi-trees in $G$. Hence $K$ contributes $1$ to the right hand side of the formula.
    
    As to the left hand side $P_{D}(v,z)$, let us apply Lemma \ref{l:induction_internally_active_edges} repeatedly for $f$ and $\AuxContr{G,K,f}$ in the order $<_K$. We see that for each $f\in E(G)$, the auxiliary graph $\AuxContr{G,K,f}$ agrees with $\Seifert(D_{G,K,f})$ for an appropriate diagram $D_{G,K,f}$, and $P_{D}(v,z)=P_{D_{G,K,e_0}}(v,z)=\dots =P_{D_{G,K,\infty}}(v,z)$ because all these diagrams express the same link.
    By Lemma \ref{l:Jaeger_quasi-trees_link_trivial}, the link in question is the trivial knot. In particular, we have $P_{D_{G,K,\infty}}(v,z)=1$, meaning that the formula is correct in this case, contradicting our assumption.

    Now suppose that $e\neq\infty$, that is, $e\in E(G)$. We distinguish two possibilities: $e\notin K$ for each $K\in\jaeg(G)$, or $e\in K$ for some Jaeger quasi-trees $K$ of $G$ and $e\notin K$ for others.

    Case 1: $e\notin K$ for each $K\in\jaeg(G)$.
    In this case, for an arbitrary $K\in\jaeg(G)$, the edge $e$ is externally active in $K$, and (by the minimality of $e$) each edge $f<_K e$ has $f\in K$, and $f$ is internally active in $K$. 
    By applying Lemma \ref{l:induction_internally_active_edges} successively to the edges $f \leq_K e$ in the order $<_K$, we get that $\AuxContr{G,K,e}$ is the Seifert graph of a link $D_{G,K,e}$, and $P_{D_{G,K,e}}=P_{D}$.

    Let $g$ be the edge immediately following $e$ in $<_K$.
    Since $e$ is externally active in $K$ with respect to $G$, we can apply Lemma \ref{l:handling_edges_that_need_to_be_removed} to $\AuxContr{G,K,e}$ and obtain that $\AuxContr{G,K,g}=\AuxContr{G,K,e}-e=\Seifert(D_{G,K,g})$ for an appropriate diagram $D_{G,K,g}$. Furthermore, if $e$ is a positive edge, then
\begin{equation}\label{eq:removing_positive_externally_active}
    P_{D}(v,z)=P_{D_{G,K,e}}(v,z)=(vz+vz^{-1}-v^3z^{-1})\cdot P_{D_{G,K,g}}(v,z).
    \end{equation}
    If $e$ is a negative edge, then again by Lemma \ref{l:handling_edges_that_need_to_be_removed}, we have
\begin{equation}\label{eq:removing_negative_externally_active}
    P_{D}(v,z)=(-v^{-1}z-v^{-1}z^{-1}+v^{-3}z^{-1})\cdot P_{D_{G,K,g}}(v,z).
    \end{equation}

    By Lemmas \ref{lem:Jaeger_of_B_is_active_and_Jaeger} and \ref{lem:Jaeger_in_G_and_in_B} we obtain a natural bijection between sets of Jaeger quasi-trees, expressed as 
    $$\jaeg(G)=\{J^G \mid J\in\jaeg(\AuxContr{G,K,g})\}.$$
    Take some $J\in\jaeg(\AuxContr{G,K,g})$. By Lemma \ref{lem:Jaeger_in_G_and_in_B}, the activity status of each edge $f\in \AuxContr{G,K,g}$ is the same with respect to $J$ (and $\jaeg(\AuxContr{G,K,e})$) as with respect to $J^G$ (and $\jaeg(G)$). The edges of $J^G$ that are contracted in $J$ are internally active in $J^G$, so they do not contribute to the formula. The edge $e$ is externally active in $J^G$, and $e$ is not an edge of $\AuxContr{G,K,g}$, wherefore it gives no activity to $J$. Putting these together, we see that the formula holds for $D$ if and only if it holds for $D_{G,K,g}$, and thus the latter must also be a counterexample.
    
    However, again by Lemma \ref{lem:Jaeger_in_G_and_in_B}, for each $J\in\jaeg(\AuxContr{G,K,g})$, every edge $f\in J$ with $f<_J g$ is internally active in $J$.
    Hence $|E(\AuxContr{G,K,g})|-t(\AuxContr{G,K,g})\leq |E(G)|-t(G)$. (Indeed, this is the number of edges that come after the first difference of Jaeger quasi-trees. Their number can only decrease, due to the contractions in $\AuxContr{G,K,g}$.) Also, we have $|E(\AuxContr{G,K,g})|<|E(G)|$ because $e$ was deleted (if there were contractions as well, that further decreases $|E(\AuxContr{G,K,g})|$). This means that $D_{G,K,g}$ is a smaller counterexample than $D$, a contradiction.

    Case 2: We have $e\in L$ for some $L\in\jaeg(G)$ and $e\notin N$ for some $N\in\jaeg(G)$. 
    Fix an arbitrary Jaeger quasi-tree $K\in\jaeg(G)$. By the minimality of $e$, each edge $g <_K e$ is in $K$ and indeed, $g$ is internally active in $K$. 
    Applying Lemma \ref{l:induction_internally_active_edges} repeatedly to each $g<_K e$, we get that $\AuxContr{G,K,e}$ is the Seifert graph of a link diagram $D_{G,K,e}$, and $P_{D_{G,K,e}}(v,z)=P_{D}(v,z)$.

    Next, apply the skein relation to $e$ in the diagram $D_{G,K,e}$. Let $D_1$ be the diagram obtained by switching the crossing corresponding to $e$, and let $D_2$ be the diagram obtained by smoothing the crossing corresponding to $e$. We have $\Seifert(D_1)=(\AuxContr{G,K,e})_{\overleftarrow{e}}$ and $\Seifert(D_2)=\AuxContr{G,K,e}-e$.
    
    If $e$ is a positive edge, then 
    \begin{equation}\label{eq:skein_positive}
    P_{D_{G,K,e}}(v,z)=v^2 P_{D_1}(v,z) + vz P_{D_2}(v,z),
    \end{equation} while if $e$ is a negative edge, then \begin{equation}\label{eq:skein_negative}
    P_{D_{G,K,e}}(v,z)=v^{-2} P_{D_1}(v,z) - v^{-1}z P_{D_2}(v,z).
    \end{equation}
    
    Fix an arbitrary Jaeger quasi-tree $L\in\jaeg(G)$ such that $e\in L$, and an arbitrary Jaeger quasi-tree $N\in\jaeg(G)$ such that $e\notin N$.
    Let $f$ be the edge following $e$ in $<_L$,  
    and let $f'$ be the edge following $e$ in $<_N$. 
    By the choice of $e$, the tours of $L$ and $N$ agree up to $e$, hence $e$ is a first difference between $L$ and $N$. Thus, $e$ is internally passive in $L$ and externally passive with respect to $N$. Hence, by Theorem \ref{thm:internal_activity_condition}, $\AuxContr{G,L,e}-e$ has no bad separating cycles, implying $\AuxContr{G,L,f}=(\AuxContr{G,L,e})_{\overleftarrow{e}}=(\AuxContr{G,K,e})_{\overleftarrow{e}}=\Seifert(D_1)$ and $\AuxContr{G,N,f'}=\AuxContr{G,N,e}-e=\AuxContr{G,K,e}-e=\Seifert(D_2)$.
    
    By Lemmas \ref{lem:Jaeger_in_G_and_in_B} and \ref{lem:Jaeger_of_B_is_active_and_Jaeger}, $$\{L'\in\jaeg(G)\mid e\in L'\}=\{J^{G}\mid J\in\jaeg(\AuxContr{G,L,f})\}$$ and 
    $$\{N'\in\jaeg(G)\mid e\notin N'\}=\{J^{G}\mid J\in\jaeg(\AuxContr{G,N,f'})\}.$$ Hence we also have
    \begin{equation}\label{eq:quasi-tree_bontas}
    \jaeg(G)=\{J^{G}\mid J\in \jaeg(\AuxContr{G,L,f})\}\cup \{J^G\mid J\in \jaeg(\AuxContr{G,N,f'})\}.     
    \end{equation}
    
    Now let us compare the activities of some $J\in\jaeg(\AuxContr{G,L,f})$ and the corresponding $J^G\in\jaeg(G)$. If $J\in\jaeg(\AuxContr{G,L,f})$, then $e$ is internally passive in $J^G$ (with respect to $G$), because $e$ is a first difference between $J^G$ and $N$. On the other hand, by part $(iii)$ of Lemma \ref{lem:Jaeger_of_B_is_active_and_Jaeger}, the edge $e$ is 
    internally active in $J$ (with respect to $\AuxContr{G,L,f}$). 
    
    As to the rest of the edges, we may say the following: 
    Since all edges $g<_L e$ are contained in each Jaeger quasi-tree of $G$, all edges $g <_L e$ are internally active in $J^G$. By Lemma \ref{lem:Jaeger_of_B_is_active_and_Jaeger}, each edge $g\in\AuxContr{G,L,f}$ with $g <_J e$ is internally active in $J$. An edge $g\in\AuxContr{G,L,f}$ with $g \geq_{J} f$ is internally/externally active/passive in $L'$ if and only if $g$ is internally/externally active/passive in $J^G$. Moreover, if an edge $g \geq_L e$ is contracted in $\AuxContr{G,L,f}$, then $g$ is internally active in $J^G$, and it is not present in $J$.
    
    In other words, the Jaeger quasi-trees $J$ of $\AuxContr{G,L,f}$ correspond to those Jaeger quasi-trees $J^G$ of $G$ that contain $e$. 
    The internal activities of $J$ and $J^G$ might differ, but our formula does not depend on internal activities. On the other hand, we have a good control on internal passivities and external activities and passivities:
    Each $J\in\jaeg(\AuxContr{G,L,f})$ has one less internal passivity than $J^G$ (and the sign in which the decrease occurs is just the sign of $e$), while the external activities and passivities of $J$ and $J^G$ remain the same.  
    
    In terms of formulas, if $e$ is a positive edge, then
    $$ep^+_G(J^{G})=ep^+_{\AuxContr{G,K,f}}(J),\quad ep^-_G(J^{G})=ep^-_{\AuxContr{G,K,f}}(J),$$
    $$ea^+_G(J^{G})=ea^+_{\AuxContr{G,K,f}}(J),\quad ea^-_G(J^{G})=ea^-_{\AuxContr{G,K,f}}(J),$$
    $$ip^+_G(J^{G})=ip^+_{\AuxContr{G,K,f}}(J)+1,\quad ip^-_G(J^{G})=ip^-_{\AuxContr{G,K,f}}(J).$$
    For a negative edge $e$ we have a similar situation, with the only difference being that $ip^-_G(J^{G})=ip^-_{\AuxContr{G,K,f}}(J)+1$ and the positive internal passivity is not changed.

    If $J\in\jaeg(\AuxContr{G,N,f'})$ then $e$ is externally passive in $J^{G}$ (with respect to $G$), since $e$ is a first difference between $J^{G}$ and $L$. On the other hand, $e$ is not an edge of  $\AuxContr{G,N,f'}=\AuxContr{G,L,e}-e$. 
    Regarding the rest of the edges, we claim the following: 
    Again, each edge $g <_J^G e$ is internally active in $J^G$ (by the choice of $e$), and each edge $g\in \AuxContr{G,N,f'}$ with $g <_J e$ is internally active in $J$ (by Lemma \ref{lem:Jaeger_of_B_is_active_and_Jaeger}).
    Still by its choice, $e$ is the only edge that is present in $G$, but deleted from $\AuxContr{G,N,f'}$. By part $(iii)$ of Lemma \ref{lem:Jaeger_of_B_is_active_and_Jaeger}, if some edge $g\geq_N f'$ of $G$ is contracted in $\AuxContr{G,N,f'}$, then $g$ is internally active in $J^G$, and $g$ is not an edge of $J$. Again by $(iii)$ of Lemma \ref{lem:Jaeger_of_B_is_active_and_Jaeger}, the activity status of each edge $g\in\AuxContr{G,N,f'}$ with $g>_N e$ is the same for $J^{G}$ and for $J$. In other words, the Jaeger quasi-trees of $\AuxContr{G,N,f'}$ correspond to those Jaeger quasi-trees of $G$ that do not contain $e$, and compared to $G$, each of these Jaeger quasi-trees loses $1$ from one of its external passivities (the sign depends on the sign of $e$), but their external activities and internal passivities remain the same.  
    (As before, the internal activities might change because of contracted edges, but internal activities do not affect the formula, whereby this does not cause any problems.)

    In $\Seifert(D_1)=\AuxContr{G,L,f}=(\AuxContr{G,K,e})_{\overleftarrow{e}}$, each Jaeger quasi-tree contains $e$, whence $|E(\AuxContr{G,L,f})|-t(\AuxContr{G,L,f}) < |E(G)|-t(G)$. This means that $D_1$ is not a counterexample to the theorem. For $\Seifert(D_2)=\AuxContr{G,N,f'}=\AuxContr{G,K,e}-e$, we have $|E(\AuxContr{G,N,f'})|-t(\AuxContr{G,N,f'}) \leq |E(G)|-t(G)$, but $|E(\AuxContr{G,N,f'})|<|E(G)|$. Hence $D_2$ is also not a counterexample to the theorem. Using the skein relation (i.e., \eqref{eq:skein_positive} or \eqref{eq:skein_negative}, depending on whether $e$ is a positive or a negative edge), together with \eqref{eq:quasi-tree_bontas} and the above remarks on the activities, we obtain that the statement of the Theorem holds for $D$, which contradicts our assumption that $G$ is a counterexample.
\end{proof}

\section{Properties of the Homfly polynomial}\label{s:applications}

In this section we will use our main theorem to construct transparent proofs of some old and new properties of the Homfly polynomial.

\subsection{Maximum $z$-degree part of the Homfly polynomial}
Let us start with showing how Theorem \ref{thm:Homfly_formula} implies Morton's bound on the maximum $z$-degree of the Homfly polynomial. Then we turn to a formula for the top $z$-degree part of the Homfly polynomial that generalizes a result of Mura\-kami and the first author, and to the product formula of Murasugi--Przytyczki for the top $z$-degree part of the Homfly polynomial with respect to Mura\-sugi sum.

\begin{thm}\cite{morton1986seifert}\label{t:morton_z-bound}
    If $D$ is an oriented link diagram with Seifert graph $G$, then
    the maximum $z$-degree of $P_{D}$ is at most $|E(G)|-|V(G)|+1$.
\end{thm}

\begin{proof}
    To apply Theorem \ref{thm:Homfly_formula}, we need a diagram $D$ such that $\Seifert(D)$ has a Jaeger quasi-tree. However, we can achieve this using Reidemeister II moves as explained in Lemma \ref{lem:link_has_diagram_with_connected_acyclic_seifert}, and that operation does not change the value of $|E(\Seifert(D))|-|V(\Seifert(D))|+1$. 
    Hence we may suppose that $G=\Seifert(D)$ has a Jaeger quasi-tree.
    
    Consider an arbitrary Jaeger quasi-tree $K$ of $G$ and its contribution to $P_D$, which takes the form of a product. To this product, from each edge $e\in K$ comes a term with $z$-degree $0$ (the term is $1$ if the edge is internally active, it is $v^2$ if the edge is internally passive and positive, and $v^{-2}$ if the edge is internally passive and negative). On the other hand, out of each edge $e\notin K$, we get a term with $z$-degree $1$ (one of $vz, -v^{-1}z, (vz+vz^{-1}-v^3z^{-1})$, or $(-v^{-1}z-v^{-1}z^{-1}+v^{-1}z^{-1})$). 
    Thus we simply need 
    an upper bound on the number of non-edges of a Jaeger quasi-tree. Since any quasi-tree is a connected spanning subgraph of $G$, it needs to have at least $|V(G)|-1$ edges. Therefore $|E(G)|-|V(G)|+1$ is indeed an upper bound on the number of non-edges. 
\end{proof}

We also see that terms of $z$-degree $|E(G)|-|V(G)|+1$ are contributed exactly by Jaeger trees. Notice though that Jaeger trees may not exist, even after making sure that Jaeger quasi-trees do. It might also be that they exist but their highest $z$-degree contributions cancel each other out.
Thus the estimate of Theorem \ref{t:morton_z-bound} may fail to be sharp, depending on the diagram.

Continuing the same argument, we can also give a formula for the maximum $z$-degree part of the Homfly polynomial.

\begin{notation}
    For an oriented link diagram $D$ with Seifert graph $G$, we denote by $h_D(v)$ the polynomial in $v$ such that $z^{|E(G)|-|V(G)|+1}h_D(v)$ is the sum of those monomials in the Homfly polynomial of $D$ in which the degree of $z$ is $|E(G)|-|V(G)|+1$. Let us stress that, in light of Morton's bound, $h_D$ depends on $D$ very mildly. It is the same polynomial for all diagrams of an oriented link for which the bound is sharp, and it is $0$ for all other diagrams.
\end{notation}

Let us also use the notation $e_G^+(T)$ for the number of positive external edges of $T$. In other words, $e_G^+(T)$ is the number of positive edges in $G-T$. Similarly, let $e_G^-(T)$ be the number of negative external edges of $T$. Here $T$ can be any set of edges of the signed graph $G$, but mostly it will mean a spanning tree.

\begin{thm}\label{t:I_G_formula_using_trees}
    Let $D$ be an oriented link diagram with Seifert graph $G$, and let $(v_0,e_0)$ be an arbitrary half-edge of $G$ such that $\jaeg(G,v_0,e_0)\neq\varnothing$. Then 
    $$h_D(v)=\sum_{T \in\jaegtree(G,v_0,e_0)} (-1)^{e^-(T)}\,v^{2ip^+(T)-2ip^-(T)+e^+(T)-e^-(T)},$$
    where each activity is meant with respect to $\jaeg(G,v_0,e_0)$.
\end{thm}

\begin{proof}
    By the argument in the proof of Theorem \ref{t:morton_z-bound}, we need to aggregate the contributions of Jaeger trees to the part of the formula of Theorem \ref{thm:Homfly_formula} where the exponent of $z$ is the highest possible.
After noting that the parity of $ep^-(T)-ea^-(T)$, in the exponent of $-1$, is the same as that of $ep^-(T)+ea^-(T)=e^-(T)$, that gives exactly the formula in the statement of the Theorem.
\end{proof}

If all crossings in the diagram, that is to say all edges in the Seifert graph, are positive, then we get a particularly nice form and we recover an earlier result.

\begin{coroll}\cite{KM}\label{coroll:Homfly_top_for_positive_links}
    If $D$ is a positive diagram, and $G$ is the Seifert graph of $D$, then for any half-edge $(v_0,e_0)$ of $G$ 
    we have
    \begin{align*}
    h_D(v)&= v^{|E(G)-|V(G)|+1} \sum_{T \in\jaegtree(G,v_0,e_0)} v^{2ip^+_{\jaeg(G,v_0,e_0)}(T)} \\&= v^{|E(G)-|V(G)|+1} \sum_{T \in\jaegtree(G,v_0,e_0)} v^{2ip^+_{\jaegtree(G,v_0,e_0)}(T)}.
    \end{align*}
\end{coroll}

\begin{proof}
    By Proposition \ref{prop:pos_diag_Seifert_acyclic}, $G$ is acyclic, hence $\jaeg(G,v_0,e_0)\neq\varnothing$.
    The first equality follows directly from Theorem \ref{t:I_G_formula_using_trees}. The second one (stating that we can take internal passivities with respect to Jaeger trees only, instead of with respect to all Jaeger quasi-trees) is a consequence of the graph-theoretic statement below.
\end{proof}

\begin{prop}\label{prop:activity_of_a_tree_w_respect_to_trees_and_quasi-trees}
    If $G$ is an acyclic rooted ribbon digraph and $T\in\jaegtree(G)$, then an edge $e\in T$ is internally passive in $T$ with respect to $\jaegtree(G)$ if and only if it is internally passive with respect to $\jaeg(G)$.
\end{prop}

\begin{proof}
    Since $\jaegtree(G)\subseteq \jaeg(G)$, if $e$ is internally passive with respect to $\jaegtree(G)$, then it is also internally passive with respect to $\jaeg(G)$.

    For the converse, suppose that the Jaeger tree $T$ and its edge $e$ are such that there exists a Jaeger quasi-tree $K\in\jaeg(G)$ such that the tours of $T$ and $K$ agree until $e$, but $e\notin K$. We need to show that there is also a Jaeger tree with these properties.
    
    Since $T$ is a tree, $T-e$ has two connected components. Let us call their vertex sets $V_0$ and $V_1$ (with the root vertex $v_0\in V_0$). 
    Note that if an edge $e'\neq e$ connects $V_0$ and $V_1$, then it interlaces with $e$ in the tour of $T$. 
    (More precisely, the two occurrences of $e'$ in the tour interlace with the two occurrences of $e$.) 
    By the Jaeger property, $(t(e'),e') <_T (h(e'),e')$. Hence we have $e' <_T e$ if any only if $t(e')\in V_0$.
    
    As $K$ is connected, there needs to be an edge $f\in K-T$ connecting $V_0$ and $V_1$. Since the tours of $T$ and $K$ agree until $e$, we have $f>_T e$ and thus $t(f)\in V_1$. As explained in \cite[Lemma 6.4]{semibalanced}, in this case $e$ is internally passive with respect to the set $\naivejaegtree(G)$ of naive Jaeger trees. 
    (Essentially, the argument is that we can keep the part of $T$ within $G[V_0]$, choose the first edge $f$ (after $e$) according to $<_T$ from the cut $(V_0,V_1)$, and take a naive Jaeger tree in $G[V_1]$ by choosing the root to be the half-edge following $(h(f),f)$. 
    The argument in \cite{semibalanced} is for a special class of acyclic digraphs. However, the only instance where it uses the special property of the graph is through the fact that $G[V_1]$ has a naive Jaeger tree, which is also true for any acyclic graph. Indeed, if we greedily remove each edge that is reached at its tail, and whose removal does not leave the graph disconnected, then we get a naive Jaeger tree.) Since $G$ is acyclic, by Proposition \ref{prop:naive_Jaeger_tree_is_Jaeger_tree} we have $\naivejaegtree(G)=\jaegtree(G)$ and thus $e$ is indeed internally passive with respect to $\jaegtree(G)$.
\end{proof}

\begin{remark}\label{rem:homfly_top_relations}
    Let us explain the subtle relationship of Corollary \ref{coroll:Homfly_top_for_positive_links} to a result of Mura\-kami and the first author \cite{KM}. They showed that positive special alternating links satisfy
    $$h_D(v)=v^{|E(G)-|V(G)|+1} I_G(v^2),$$ 
    where $G$ is the Seifert graph of a (positive) special alternating diagram $D$ and $I_G$ is the so-called interior polynomial of $G$. The interior polynomial has several equivalent definitions. In particular, in \cite{hyperBernardi}, the authors of the current paper show that (for the Seifert graph $G$ of a special alternating diagram) $I_G$ can be obtained as the generating function of internal passivities of naive Jaeger trees of $G$. (In \cite{hyperBernardi}, this is defined as internal passivity with respect to naive Jaeger trees, that is, an edge of a naive Jaeger tree is internally passive if it is the first difference between two naive Jaeger trees.) As we show in Section \ref{s:Jaeger_and_naive_Jaeger} (see Proposition \ref{prop:naive_Jaeger_tree_is_Jaeger_tree}, and note that the Seifert graph of a special alternating diagram is always acyclic) naive Jaeger trees of $G$ are the same as Jaeger trees of $G$.
\end{remark}

Theorem \ref{thm:Homfly_formula} also enables us to give a simple proof of a known product formula for $h_L$, under Murasugi sums.

The notion of \emph{diagrammatic Murasugi sum} is defined as follows.
Suppose that $D_1$ and $D_2$ are oriented link diagrams. We will think of them as drawings on the sphere, in fact, for ease of exposition, on two disjoint spheres. Choose type~I Seifert circles $s_i$ in $D_i$ for $i=1,2$. Let $s_i$ bound the open disk $R_i$ of the sphere, where $R_i$ is disjoint from $D_i$. Remove the $R_i$ from the respective spheres and identify the boundaries of the remaining disks, that is $s_1$ and $s_2$, to form a new sphere, with the \emph{Murasugi sum} $D_1\star D_2$ drawn on it. Here the identification is to be performed so that the orientations of $s_1$ and $s_2$ match up, furthermore no crossing of $D_1$ along $s_1$ is identified with a crossing of $D_2$ along $s_2$. That is, we choose a way of merging the two cyclic orders of crossings, around $s_1$ and $s_2$, into a single cyclic order. This choice, as well as the choice of the $s_i$, obviously influences not only the diagram but also the link type of the outcome\footnote{A notable exception is the connected sum of two (oriented) knots, which can be represented with the Murasugi sum of arbitrarily chosen diagrams, along arbitrary type I Seifert circles, using any merging of the cyclic orders that avoids 
interlacement.}. Yet they are suppressed from the notation $D_1\star D_2$ and they do not appear in the following theorem, either.

Before the theorem let us note that alternatively, we may think of the $R_i$ as subsets of the respective Seifert surfaces, which we then identify to form a new surface and a new link as its boundary.
One may also observe that, apart from some trivial cases, the common image of $s_1$ and $s_2$ is a type II Seifert circle in the new diagram. All its other Seifert circles are inherited from $D_1$ or $D_2$, and the Seifert graph of the new diagram, as an abstract graph, is 
obtained from $\Seifert(D_1)$ and $\Seifert(D_2)$ by gluing together the vertices $s_1$ and $s_2$ into a new vertex $s$. In particular, $s$ is a cut vertex of $\Seifert(D_1 \star D_2)$.

\begin{thm}[Murasugi--Przytyczki \cite{MP_szorzodas}]
\label{t:product_for_Murasugi_sum}
    Let $D_1$ and $D_2$ be oriented link diagrams and consider any Murasugi sum $D_1\star D_2$ formed from them. Then we have \[h_{D_1\star D_2}=h_{D_1}\cdot h_{D_2}.\]
\end{thm}

By Theorem \ref{t:I_G_formula_using_trees}, the polynomial $h_D$ can be computed using Jaeger trees. We had Theorem \ref{thm:internal_activity_condition} characterizing internal passivities in Jaeger quasi-trees. For Jaeger trees, we can give a simpler characterization that will be useful for proving Theorem \ref{t:product_for_Murasugi_sum}. We start with this characterization.

\begin{prop}\label{prop:simplified_activity_description_for_trees}
    If $T$ is a Jaeger tree of the connected rooted ribbon digraph $G$, then an edge $e\in T$ is internally active in $T$ (with respect to $\jaeg(G)$) if and only if either $(h(e),e)<_T (t(e),e)$, or else $(t(e),e) <_T (h(e),e)$ and $e$ is a bridge in $\AuxContr{G,K,e}$.
\end{prop}

\begin{proof}
    By Theorem \ref{thm:internal_activity_condition}, the edge $e$ is internally active for $T$ if and only if $e$ is tight. 
    Recall that an $e$ edge is tight for $T$ if either $(h(e),e) <_T (t(e),e)$, or $(t(e),e) <_T (h(e),e)$ and $\AuxContr{G,T,e}-e$ has a bad separating cycle, or $e$ is contracted when we form $\AuxContr{G,T,e}$ (see Definition \ref{def:type0_type1_and_2_edges_in_K}). 
    If we have an edge $e\in T$ such that $\AuxContr{G,T,e}-e$ has a bad separating cycle $C$, then we need $C=\varnothing$ (in other words, $e$ must be a bridge in $\AuxContr{G,T,e}$), for otherwise (as pointed out in Remark \ref{rem:osszehuzott_kor_bekerül_a_kvazifaba}) we would have $C\subseteq T$, which is not possible for a tree. Hence, in particular, we do not contract any edges during the Jaeger test of $T$, wherefore we do not have any type~$2$ edges for $T$.
\end{proof}

\begin{proof}[Proof of Theorem \ref{t:product_for_Murasugi_sum}]
We may suppose that the graphs $G_i=\Seifert(D_i)$ are acyclic ($i=1,2$), for else we may perform some Reidemeister II moves to obtain diagrams $D'_i$, cf.\ Lemma \ref{lem:link_has_diagram_with_connected_acyclic_seifert}. These steps introduce some new Seifert circles and have the effect of substituting an edge of the Seifert graph with a path of length $3$. Importantly, they leave the existing Seifert circles intact. Hence we can construct the same Murasugi sum for $D'_1$ and $D'_2$ along $s_1$ and $s_2$, and the resulting diagram $D'_1 \star D'_2$ is also related to $D_1 \star D_2$ by a series of Reidemeister II moves. These do not alter either side of the equation in the Theorem. 

From this point we suppose that $G_1=\Seifert(D_1)$ and $G_2=\Seifert(D_2)$ are acyclic. This implies that $G=\Seifert(D_1\star D_2)$ is also acyclic.
Let $V_i=V(G_i)-s_i$ ($i=1,2$). Apart from these sets, $V(G)$ contains the cut vertex $s$ which is
the type II Seifert circle obtained by gluing $s_1$ and $s_2$.

Let $e_1$ be an arbitrary edge incident to $s$ such that the other endpoint of $e_1$ is in $V_1$. 
Let $e_2$ be the earliest edge following $e_1$ in the cyclic order around $s$ such that the other endpoint of $e_2$ is in $V_2$.
Take $(s,e_1)$ as the root of $G$. 
By Theorem \ref{t:I_G_formula_using_trees}, the polynomial  $h_{D_1\star D_2}$ can be written as a generating function of activities for the Jaeger trees of $G$.
Now, we claim that 
$$\jaegtree(G,s,e_1)=\{T_1\cup T_2\mid T_1\in\jaegtree(G_1,s,e_1),\,T_2\in \jaegtree(G_2,s,e_2)\}.$$
Proposition \ref{prop:naive_Jaeger_tree_is_Jaeger_tree} implies that Jaeger trees of $G$, $G_1$, and $G_2$, respectively, are the same as naive Jaeger trees in the respective graphs.

Let $T$ be an arbitrary 
Jaeger tree of $G$ (with root $(s,e_1)$). For any vertex $u$, the tour of any tree takes the half-edges $\{(u,f)\mid f \text{ is incident to } u\}$ in an order compatible with the ribbon structure around $u$. Hence for any $T\in\jaegtree(G,s,e_1)$, the first half-edge of $G_2$ considered in the tour of $T$ is $(s,e_2)$. Let $T_i$ ($i=1,2$) be the subset of $T$ in (the edge set of) $G_i$. The restriction of the tour of $T$ to half-edges of $G_i$ is the tour of $T_i$, with root $(s,e_i)$. As $T$ is a naive Jaeger tree, we see that $T_i\in\jaegtree(G_i,s,e_i)$.

Conversely, if $T_i\in\jaegtree(G_i,s,e_i)$, then $T=T_1\cup T_2$ is a tree, and the restriction of the tour of $T$ to $G_i$ is the tour of $T_i$ with root $(s,e_i)$, implying that $T$ is a naive Jaeger tree of $G$.

We also claim that if $T=T_1\cup T_2$ for $T_i\in\jaegtree(G,s,e_i)$, and $e\in T_i$ for some $i=1,2$, then $e$ is internally active in $T$ (with respect to $G$) if and only if $e$ is internally active in $T_i$ (with respect to $G_i$). We apply Proposition \ref{prop:simplified_activity_description_for_trees}. Note that $(h(e),e)<_T (t(e),e)$ if and only if $(h(e),e) <_{T_i} (t(e),e)$. Moreover,
$s$ remains a cut vertex in $\AuxContr{G,T,e}$, and 
the part of the auxiliary graph in $G_i$ coincides with
$\AuxContr{G_i,T_i,e}$. Hence $e$ is a bridge in $\AuxContr{G,T,e}$ if and only if $e$ is a bridge in $\AuxContr{G_i,T_i,e}$. This proves that indeed $e$ is internally active in $T_i$ if and only if $e$ is internally active in $T$.

Note also that external edges of $T$ in $G$ are exactly those edges that are external edges of $T_i$ in $G_i$ for $i=1$ or $i=2$.
The same holds for the sub-collections of positive external edges and negative external edges, too.

We conclude that any pair of Jaeger trees, say $T_1\in\jaegtree(G_1,s_1,e_1)$ and $T_2\in\jaegtree(G_2,s_2,e_2)$, gives rise to a Jaeger tree $T=T_1\cup T_2$ in $\jaegtree(G,s,e_1)$ in a bijective way. 
As to the ingredients in Theorem \ref{t:I_G_formula_using_trees},
we do have $ip^+_{\jaeg(G)}(T)=ip^+_{\jaeg(G_1)}(T_1)+ip^+_{\jaeg(G_2)}(T_2)$ and $e_G^+(T)=e_{G_1}^+(T_1)+e_{G_2}^+(T_2)$, and similarly for the negative versions.
\end{proof}

\begin{remark}
    Note that it was important in our argument to choose the root incident to $s$, because otherwise we would have no control on the first half-edge of the tour of $T$ within $G_2$. Similarly, it was crucial that we only had to talk about trees, because for quasi-trees, we again have no control on which half-edge starts the tour of $T$ within $G_2$.
\end{remark}

\subsection{Minimum $z$-degree of the Homfly polynomial}

Lickorish and Millett proved a formula on the minimum $z$-degree of the Homfly polynomial \cite{lickorish1987polynomial}.
Their lower bound also follows easily from Theorem \ref{thm:Homfly_formula}, and the characterization of externally active edges.
\begin{prop}\label{prop:min_z_degre_lower_bound}\cite{lickorish1987polynomial}
   Let $D$ be a link diagram with $b(D)$ components. Then, the minimal $z$-degree of the Homfly polynomial $P_D(v,z)$ is $1-b(D)$. 
\end{prop}
Let us use the notation $P_D(v,z)=\sum_{i}P^i_D(v)z^i$. 

\begin{proof}
Let $G=\Seifert(D)$ and recall that $b(G)$, the number of boundary components of $G$, is equal to $b(D)$.

By Theorem \ref{thm:Homfly_formula} the minimum $z$-degree monomial  contributed by a Jaeger quasi-tree $K$ has $z$-degree $ep(K)-ea(K)$.
We show that $ep(K)-ea(K)\geq 1-b(G)$ for each Jaeger quasi-tree $K$.

Let us examine $b(\AuxContr{G,K,e})$ throughout the Jaeger test of $K$. 
If $e\in K$, then for the edge $f$ following $e$ in $<_K$, we have $b(\AuxContr{G,K,f})=b(\AuxContr{G,K,e})$. (The two graphs are either identical, or related by an edge reversal, or by a contraction, but none of these changes the number of boundary components). 

If we have $e\notin K$ then $\AuxContr{G,K,f}=\AuxContr{G,K,e}-e$, and the number of boundary components either increases or decreases by one.
If $e$ is an externally passive edge with respect to $K$, then by Theorem \ref{thm:ext_activity_sufficient}, the graph $(\AuxContr{G,K,e})_{\overleftarrow{e}}$ has a bad separating cycle. Since $\AuxContr{G,K,e}$ has no bad separating cycles, $\overleftarrow{e}$ must be in the bad separating cycle of $(\AuxContr{G,K,e})_{\overleftarrow{e}}$. In particular, $e$ is in a separating cycle of $\AuxContr{G,K,e}$, and thus its two half-edges are in different boundary components. As a consequence, $b(\AuxContr{G,K,f})=b(\AuxContr{G,K,e}-e)=b(\AuxContr{G,K,e})-1$.

Since $K$ is a quasi-tree, we have $b(\AuxContr{G,K,\infty})=1$. Thus, starting from $G$, each external edge increases or decreases the number of boundary components by at most 1, and each externally active edge decreases it by exactly one; at the end, the number of components becomes 1.
We conclude that $b(G)-ea(K)+ep(K)\geq 1$, so $ep(K)-ea(K)\geq 1-b(G)$.

This also implies that each quasi-tree that contributes to $P^{1-b(D)}_D$ has at least $b(D)-1$ externally active edges, whence it is possible to write $P^{1-b(D)}_D=(v^{-1}-v)^{b(D)-1} p(v)$, where $p(v)$ is still a Laurent polynomial. We will argue that $p(1)=\pm 1$. Thus indeed, there are nonzero monomials with $z$-degree $1-b(D)$.

To show $p(1)=1$, it is enough to show that there exists a unique Jaeger quasi-tree $K$ such that $ea(K)=b(G)-1$, and $ep(K)=0$. Indeed, for such a Jaeger quasi-tree $K$ one has contribution $\pm 1$ to $p(1)$ (depending on the parity of the number of negative external edges). For any Jaeger quasi-tree that has larger external activity, there is a contribution of $\pm(v^{-1}-v)$ to $p(v)$, whence these quasi-trees contribute $0$ to $p(1)$.

Since each Jaeger quasi-tree $K$ has $ep(K)-ea(K)\geq 1-b(G)$, it is enough to show that there exists a unique Jaeger quasi-tree $K$ such that $ep(K)=0$.

We proceed similarly to the proof of Theorem \ref{thm:if_no_oriented_sep_cycle_then_exists_Jaeger}, but this time we include each edge into $K$, unless the currently observed edge $e$ is not contracted in $\AuxContr{G,K,e}$, we have $(t(e),e)<_K (h(e),e)$, and there is a bad separating cycle in $(\AuxContr{G,K,e})_{\overleftarrow{e}}$. In the latter case, we let $e\notin K$. By Lemma \ref{lem:contractability_of_bad_sep_cycle}, this way $\AuxContr{G,K,f}$ has no bad separating cycles. Hence using this rule, the Jaeger test is never violated, and no bad separating cycle is introduced in $\AuxContr{G,K,f}$. We need to check that $K$ is a quasi-tree, but this can be done with the exact same argument as in the proof of Theorem \ref{thm:if_no_oriented_sep_cycle_then_exists_Jaeger}. Now the construction of $K$ and Theorem \ref{thm:ext_activity_sufficient} imply that each $e\notin K$ is externally active. Because of Theorem \ref{thm:ext_activity_sufficient} if we want a quasi-tree where each external edge is active, this procedure is our only option.
\end{proof}

Morton and, in the special case of closed braid diagrams, Franks and Williams, also found diagrammatic lower and upper bounds on the $v$-degree of the Homfly polynomial. Those bounds do not seem to follow easily from Theorem \ref{thm:Homfly_formula}. 

\subsection{The total degree of the Homfly polynomial}

We turn to a question which, to the best of our knowledge, was not previously investigated, and where more interesting graph-theoretic notions start to appear.
In this section we prove an upper bound on the degree of the Homfly polynomial, which is tight if each crossing of the diagram is positive.

\begin{defn}[dijoin]
    Let $G$ be a directed graph. A \emph{directed join}, or \emph{dijoin} for short, of $G$ is a subset of the edges $F\subseteq E(G)$ such that $F$ contains an edge from each directed cut of $G$. By a directed cut we mean a non-empty set of edges $C$, arising as those connecting $V_0$ to $V_1$ for some partition $V_0\sqcup V_1=V(G)$, so that each edge of $C$ has its head in $V_1$.
\end{defn}

\begin{customthm}{\ref{thm:max_degree_dijoin}}
    Let $D$ be a link diagram such that $G=\Seifert(D)$ has a Jaeger quasi-tree with respect to some basepoint. (For instance, $G$ may be acyclic.)
    Then the degree of $P_{D}$ (where by the degree of a monomial we mean the sum of its $v$-degree and $z$-degree) is at most $2(|E(G)|-\nu(G))$, where $\nu(G)$ is the minimum size of a dijoin of $G$. If each edge of $G$ is positive, then this estimate is tight. 
\end{customthm}

    If $D$ is a diagram with only positive crossings then $\Seifert(D)$ is acyclic and the theorem can be applied.

\begin{proof}
Fix a basepoint such that $\jaeg(G)\ne\varnothing$ and 
    take an arbitrary Jaeger quasi-tree $K\in\jaeg(G)$.
    Regarding the formula of Theorem \ref{thm:Homfly_formula}, each external edge 
    contributes at most $2$ to the degree of the term of $K$, and each internally passive edge 
    contributes at most $2$ as well. On the other hand, each internally active edge 
    contributes $0$. (All of the above are exact amounts when the edge is positive.) Hence if we show that each Jaeger quasi-tree of $G$ has at least $\nu(G)$ internally active edges, then it follows that the degree of each monomial in $P_D$ is at most $2(|E(G)|-\nu(G))$.

    This property, in turn, follows if we
    prove that the internally active edges of any Jaeger quasi-tree of $G$ form a dijoin of $G$.

    Let $K$ be a Jaeger quasi-tree of $G$ and let $F\subseteq K$ be its set of internally active edges. To show that $F$ is a dijoin, we need to ascertain that $F\cap C \neq \varnothing$ for any directed cut $C$ of $G$.
    Let $C$ be such a  cut.
	
	Case 1: $C$ is directed toward its side containing the base vertex $b_0$.
    As $K$ is connected, $K\cap C$ is nonempty. 
    The first time that the tour reaches an edge $e\in C\cap K$, it has to do so at its head by our assumption. Hence $e$, being tight of type~$0$, is internally active in $K$. I.e., we have $e\in F\cap C$.
	
	Case 2: $C$ is directed away from its side containing $b_0$.
	If there is any edge in $C\cap K$ that the tour of $K$ first reaches at its head, then that edge is once again internally active in $K$. It remains to examine the case where each edge in $e\in C\cap K$ has $(t(e),e) <_K (h(e),e)$. Let $e_1, \dots, e_s$ be the edges in $C\cap K$ and suppose that 
    $e_1 <_K e_2 <_K \dots <_K e_s$. 

    Denote the two sides of the cut $C$ by $V_0$ and $V_1$, where $b_0\in V_0$. In other words, $C$ is the set of edges connecting $V_0$ to $V_1$ and each such edge has its tail in $V_0$.
    When the tour of $K$ traverses the half-edge  $(t(e_1),e_1)$, it reaches $V_1$. Since we assumed that $(t(e_2),e_2) <_K (h(e_2),e_2)$, moreover $t(e_2)\in V_0$, the tour of $K$ needs to return to $V_0$ through $(h(e_1),e_1)$ before traversing $(t(e_2),e_2)$.
    With analogous reasoning, we get that $(t(e_1),e_1) <_K (h(e_1),e_1) <_K (t(e_2),e_2)<_K (h(e_2),e_2) <_K \dots <_K (t(e_s),e_s) <_K (h(e_s),e_s)$.
    After traversing $(h(e_s),e_s)$, the tour of $K$ does not return to $V_1$ anymore. Hence for each edge spanned by $V_1$, both of its half-edges need to be traversed before $(h(e_s),e_s)$. 

    We claim that $e_s$ is internally active in $K$ (with respect to $G$). More precisely we show that it is tight of type~$1$, i.e., $\AuxContr{G,K,e_s}-e_s$ has a bad separating cycle. To establish this, consider the tour of $K^{\AuxContr{G,K,e_s}}$ in $\AuxContr{G,K,e_s}$, in particular its part between $(t(e_s),e_s)$ and $(h(e_s),e_s)$.  
    To see that $\AuxContr{G,K,e_s}-e_s$ has a bad separating cycle, it suffices to check that 
    \begin{enumerate}
        \item[(1)] there is no edge $uv\in \AuxContr{G,K,e_s}-K^{\AuxContr{G,K,e_s}}$ such that $(u,uv)$ and $(v,uv)$ interlace with $(t(e_s),e_s)$ and $(h(e_s),e_s)$, 
        \item[(2)] if there is an edge $uv\in K^{\AuxContr{G,K,e_s}}$, such that $(u,uv)$ and $(v,uv)$ interlace with $(t(e_s),e_s)$ and $(h(e_s),e_s)$, and $(t(e_s),e_s)<_K (u,uv) <_K (h(e_s),e_s)$ then $v=h_{\AuxContr{G,K,e_s}}(uv)$.
    \end{enumerate}
    Indeed if (1) and (2) hold then the edges $uv\in K^{\AuxContr{G,K,e_s}}$, such that $(u,uv)$ and $(v,uv)$ interlace with $(t(e_s),e_s)$ and $(h(e_s),e_s)$, form a separating cycle in $\AuxContr{G,K,e_s}-e_s$. Because for such an edge $uv$, the half-edge $(v,uv)$ is in the base component of $\AuxContr{G,K,e_s}-e_s$ and we assumed $v=h(uv)$, this separating cycle is bad .
    
    We proceed to show that (1) and (2) hold. 
    Let $uv$ be such that $(u,uv)$ and $(v,uv)$ interlace with $(t(e_s),e_s)$ and $(h(e_s),e_s)$.
    Suppose that we have $(t(e_s),e_s)<_K (u,uv) <_K (h(e_s),e_s)$. We then either have $(v,uv) <_K (t(e_s),e_s)$ or $(h(e_s),e_s) <_K (v,uv)$. Since the tour of $K$ is in $V_1$ between $(t(e_s),e_s)$ and $(h(e_s),e_s)$, we have $u\in V_1$.

    If $v\in V_0$ then $uv\in C$. In this case we have $uv\notin K$ because $K\cap C=\{e_1, \dots, e_s\}$ and none of these edges interlaces with $e_s$ in $<_K$. Since $C$ is a directed cut with each edge having its head in $V_1$, we have $u=h(uv)$. Now as $K$ is Jaeger, $(v,uv)=(t(uv),uv)$ precedes $(u,uv)$ in the tour of $K$. Since between the two traversals of $e_s$ the tour is in $V_1$, we deduce that $(v,uv)$ also precedes the first traversal of $e_s$. Thus $uv$ gets deleted from the auxiliary graph before reaching $e_s$, i.e., $uv\notin \AuxContr{G,K,e_s}$.

    If $v\in V_1$ then $uv$ is spanned by $V_1$. Hence by the time of the traversal of $(h(e_s),e_s)$, both half-edges of $uv$ have been dealt with. In other words we are in the situation $(v,uv) <_K (t(e_s),e_s)<_K (u,uv) <_K (h(e_s),e_s)$, in particular $uv <_K e_s$. Once again if $uv\notin K$ then this means $uv\notin \AuxContr{G,K,e_s}$, because $uv$ is already deleted.
    
    Now suppose $uv\in K$. If $uv$ were a type~$1$ edge of $K$ then $uv$ would be a bridge in $\AuxContr{G,K,e_s}$. However a bridge does not interlace with any other edge in $<_{K^{\AuxContr{G,K,e_s}}}$, in particular $e_s$ and $uv$ also could not interlace in $<_K$. 
    If $uv$ is a type~$2$ edge of $K$ then $uv\notin \AuxContr{G,K,e_s}$, whence this is also impossible. Hence either $uv$ is of type~$0$, or it is not a tight edge of $K$. In both cases, we have $(h_{\AuxContr{G,K,e_s}}(uv),uv) <_{K^{\AuxContr{G,K,e_s}}} (t_{\AuxContr{G,K,e_s}}(uv),uv)$, in other words, $h_{\AuxContr{G,K,e_s}}(uv)=v$. This finishes the proof that the active edges form a dijoin for each Jaeger quasi-tree, which implies that the inequality of the theorem holds.

    Finally we show that the inequality becomes an equality if each edge of $G$ is positive. 
    By Theorem \ref{t:morton_z-bound} the maximum $z$-degree is $|E(G)|-|V(G)|+1$, and thus $h_D(v)\ne0$. 
    Moreover, by Corollary \ref{coroll:Homfly_top_for_positive_links} and Remark \ref{rem:homfly_top_relations}, if $G$ has all positive edges, then 
    \[h_D(v)=v^{|E(G)|-|V(G)|+1} \cdot\sum_{T \in\jaegtree(G,v_0,e_0)} v^{2ip^+(T)} =v^{|E(G)|-|V(G)|+1}\, I_G(v^2),\] 
    where $I_G$ is the so-called interior polynomial of $G$. By \cite[Theorem 1.1]{fokszam}, for a connected directed graph $G$, the polynomial $I_G$ has degree $|V(G)|-1-\nu(G)$. Hence, by taking the leading term of $I_G$, we obtain a term in $P_L$ with $z$-degree $|E(G)|-|V(G)|+1$, and $v$-degree $(|E(G)|-|V(G)|+1)+2(|V(G)|-1-\nu(G))$. Altogether this gives a term with total degree $2(|E(G)|-\nu(G))$, meaning that the inequality is indeed tight if $G$ has all positive edges. 
\end{proof}

\begin{remark}
    We note that in \cite{fokszam} the degree of $I_G$ is controlled using Ehrhart theory, but it is also possible to directly construct a Jaeger tree in $G$ that has internal activity $\nu(G)$.
\end{remark}

Let us state another version of the theorem. It provides a weaker bound, but drops the assumption on the Seifert graph of the diagram.

\begin{defn}[feedback arc set]
    Let $G$ be a directed graph. A \emph{feedback arc set} of $G$ is a subset of the edges $F\subseteq E(G)$ such that $F$ contains at least one edge from each directed cycle of $G$.
\end{defn}

\begin{thm}
    Let $D$ be an arbitrary link diagram and let $G=\Seifert(D)$.
    Then the degree of $P_{D}$ is at most $2(|E(G)|-\nu(G)+f(G))$, where $\nu(G)$ is the minimum size of a dijoin of $G$, and $f(G)$ is the minimum size of a feedback arc set of $G$.
\end{thm}

\begin{proof}
    We will use Theorem \ref{thm:max_degree_dijoin}. To do so, we need to have a basepoint $(v_0,e_0)$ such that $\jaeg(G,v_0,e_0)\neq\varnothing$, and for that it is enough to make the Seifert graph acyclic. 
    As explained in Theorem \ref{thm:each_link_has_diagram_with_Jaeger_quasi-trees} and Lemma \ref{lem:link_has_diagram_with_connected_acyclic_seifert}, we can achieve this by performing Reidemeister II moves on the diagram. At the graph level, let us pick a minimal cardinality feedback arc set $F\subseteq E(G)$, and substitute each edge $\overrightarrow{uv}\in F$ with a path of length $3$, namely ($\overrightarrow{uw_1}, \overleftarrow{w_1w_2}, \overleftarrow{w_2v}$), as depicted in Figure \ref{fig:edge_subdivision_in_Seifert_graphs}. The resulting Seifert graph $G'$ is clearly acyclic, and $|E(G')|=|E(G)|+2f(G)$.
    
    We claim that $\nu(G')=\nu(G)+f(G)$.
    Since a directed cycle is always disjoint from a directed cut, any minimal dijoin of $G$ is disjoint from $F$. For each edge $\overrightarrow{uv}\in F$, the graph $G'$ contains a new sink vertex $w_1$, whence $\nu(G')\geq \nu(G)+f(G)$. On the other hand, the new edge $\overrightarrow{uw_1}$ covers the newly created directed cuts $\{\overrightarrow{uw_1}, \overrightarrow{w_2w_1}\}$ and $\{\overrightarrow{uw_1}, \overrightarrow{vw_2}\}$, implying $\nu(G')\le \nu(G)+f(G)$.
    Finally, by applying Theorem \ref{thm:max_degree_dijoin}, we get that the degree of $P_D$ is at most $2(|E(G')|-\nu(G'))=2(|E(G)-\nu(G)+f(G))$.
\end{proof}

\section{On extensions of the Homfly polynomial}\label{s:abstract_homfly}

Theorem \ref{thm:Homfly_formula} is about links and Seifert graphs of links, but the corresponding formula makes sense for arbitrary signed, rooted ribbon digraphs. 
\begin{defn}
    For a signed, rooted ribbon digraph $G$ with root $(v_0,e_0)$, we call the polynomial 
    \begin{align*}
P_{G,v_0,e_0}(v,z)=&\sum_{K\in\jaeg(G,v_0,e_0)} (v^2)^{ip^+(K)}\cdot (v^{-2})^{ip^-(K)}\cdot (vz)^{ep^+(K)}\cdot (-v^{-1}z)^{ep^-(K)} \\ &\cdot (vz+vz^{-1}-v^3z^{-1})^{ea^+(K)}\cdot (-v^{-1}z-v^{-1}z^{-1}+v^{-3}z^{-1})^{ea^-(K)}.    
\end{align*}
the abstract Homfly polynomial of the rooted ribbon digraph $G$, with respect to the root $(v_0,e_0)$.
\end{defn}

Theorem \ref{thm:Homfly_formula} together with the well-definedness of the (knot-theoretic) Homfly polynomial implies that for Seifert graphs of link diagrams, as far as there is a Jaeger quasi-tree for each choice of root, the abstract Homfly polynomial is independent of the choice of root. 
However, the well-definedness of the Homfly polynomial of oriented links is a highly non-trivial result. It would be interesting to see a direct, combinatorial proof of root-independence for Seifert graphs of link diagrams. 
On the other hand, as we show in Section \ref{ss:root_dependence_of_abstract_homfly}, the abstract Homfly polynomial $P_{G,v_0,e_0}$ \emph{can} in fact depend on the basepoint, even if the stronger sense that we can have two basepoints where there are Jaeger quasi-trees with respect to both of them, but the abstract Homfly polynomials are different. Hence any combinatorial proof for the root-independence of $P_{\Seifert(D)}$ needs to use some special properties, ones that Seifert graphs of links exhibit, in a crucial way.

As a positive result, we present another non-trivial class of ribbon graphs where the abstract Homfly polynomial is independent of the root. There will be a slightly different, but still topological reason for this, namely that the polynomial agrees with a certain variant of the Homfly polynomial, developed for oriented links in a thickened surface.

\begin{prob}
    When is the abstract Homfly polynomial independent of the root? 
\end{prob}

\begin{prob}
Can one modify the definition of the abstract Homfly polynomial in such a way that it gives a Homfly polynomial of a link in an (appropriately defined) thickened surface for any ribbon digraph?
\end{prob}

\subsection{Homfly polynomial in thickened surfaces and root-independence for edge-subdivided graphs}\label{ss:Homfly_well_def_on_subdivided_graphs}

We consider a class of (signed) ribbon digraphs derived from
arbitrary undirected ribbon graphs $G=(V,E)$. Namely, let $D_G$ be the signed ribbon digraph obtained from $G$ by subdividing each edge with a new vertex, orienting the so obtained edges toward the subdividing points, and taking each edge to be positive. At vertices of $G$, the graph $D_G$ inherits the ribbon structure from $G$. As the subdividing vertices have degree $2$, the ribbon structure around them is unique. 
We claim that for this class of signed ribbon digraphs, the abstract Homfly polynomial is independent of the choice of root. 
This class of ribbon graphs 
has prior appearances
in the literature of the Homfly polynomial. Indeed, in the case when the ribbon structure of $G$ is planar, Jaeger \cite{Jaeger_Tutte_and_link} computed the Homfly polynomial of the `chainmail' link whose Seifert graph is 
$D_G$. His formula amounts to a change of variables in 
the Tutte polynomial of $G$. 

This was generalized by Moffatt \cite{moffatt_knot_inv} to the case of $D_G$ where $G$ is a general ribbon graph. 
In that context it is natural to consider the analogous chainmail link $L_{D_G}$ in the thickened surface $\Sigma_G\times I$. (Indeed, it is easy to check that the components of $L_{D_G}$ are bijective with the boundary components of $G$, i.e., with the vertices of $G^*$. The components are linked together in a pattern given by $G^*$ and if $G^*$ lacks loop edges then all components are trivial.)
Moffatt expressed the Homfly polynomial of $L_{D_G}$ 
using  
the Bollob\'as--Riordan polynomial (more precisely, he used a two-variable evaluation which is sometimes called the Tutte polynomial of a ribbon graph). We do not introduce the Homfly polynomial of oriented links in a thickened surface here, as we will not need the definition (we refer the interested reader to \cite{moffatt_knot_inv}). We mention that the 
classical Homfly polynomial corresponds to the case where the surface is the sphere.

We claim that similarly, the abstract Homfly polynomial of $D_G$ can be given as an evaluation of the Bollob\'as--Riordan polynomial of the ribbon graph $G$. In particular, it is independent of the choice of root, and it agrees with a Homfly polynomial in a thickened surface. To spell out the details, we do need some definitions.

The Bollob\'as--Riordan polynomial of an undirected (and unsigned) ribbon graph $G$ is defined as
$$R_G(x,y,z)=\sum_{A\subseteq E(G)}(x-1)^{k(A)-k(G)}\,y^{|A|-|V(G)|+k(A)}\,z^{2k(A)+|A|-|V(G)|-f(A)},$$
where $k(A)$ denotes the number of connected components of the graph $(V(G),A)$, and $f(A)$ is the number of boundary components of the ribbon graph $(V(G),A)$.

Moffatt introduced the dichromate polynomial $Z_G$ of ribbon graphs \cite{moffatt2025constructing}. It is obtained from the Bollob\'as--Riordan polynomial as 
\begin{equation}\label{eq:Z_G_in_terms_of_R_G}
    Z_G(x,y)=x^{|V(G)|-k(G)}\,y^{k(G)}\,R_G\left(\frac{x+y}{x},xy,\frac{1}{y}\right).
\end{equation}

Our result is the following:

\begin{thm}\label{thm:abstract_homfly_well_def_elosztott_graf}
    Let $G$ be an arbitrary (undirected) ribbon graph, and let $D_G$ be the signed ribbon digraph obtained from $G$ by subdividing each edge with a new point, orienting the edges toward the subdividing points, and defining each new edge to be positive. Then, for an arbitrary root $(v_0,e_0)$ of $D_G$, we have
    \begin{equation*}
    P_{D_G,v_0,e_0}(v,z)=
    \frac{z}{v^{-1}-v}\cdot\left(\frac{1}{vz}\right)^{|E(G)|}\cdot Z_G\left(\frac{v}{z}, \frac{v^{-1}-v}{z} \right).
\end{equation*}
    In particular, the abstract Homfly polynomial of $D_G$ is independent of the choice of root.
\end{thm}

Before giving the proof, let us further discuss related results that suggested the above formula.
As mentioned above, one can think of $L_{D_G}$ as a link embedded in $\Sigma_G\times I$, where $\Sigma_G$ is the canonical surface of the ribbon graph $G$, and $I$ is an interval, see \cite{moffatt_knot_inv}. Moffatt \cite[Theorem 3.3]{moffatt_knot_inv} proved that 
\footnote{More precisely, Moffatt phrased his formula in terms of the Bollob\'as--Riordan polynomial $R_G$, but it can easily be written in this form using \eqref{eq:Z_G_in_terms_of_R_G}. Note that in \cite{moffatt_knot_inv}, the definition of $R_G$ is shifted in the first variable compared to the definition here, and also that Moffatt uses a different convention for defining the Homfly polynomial, which necessitates the substitution $v=x^{-1}$ in the first indeterminate.}
\begin{equation*}
    P^{\Sigma_G}_{L_{D_G}}(x,y)=
    \frac{y}{x^{-1}-x}\cdot\left(\frac{1}{xy}\right)^{|E(G)|}\cdot Z_G\left(\frac{x}{y}, \frac{x^{-1}-x}{y} \right).
\end{equation*}
Comparing this to Theorem \ref{thm:abstract_homfly_well_def_elosztott_graf}, we get the following corollary.

\begin{coroll}
    Let $G$ be an arbitrary ribbon graph, and let $D_G$ be the signed ribbon digraph obtained from $G$ by subdividing each edge with a new point, orienting the edges toward the subdividing points, and defining each edge to be positive. Then, for an arbitrary root $(v_0,e_0)$ of $D_G$, the abstract Homfly polynomial of $D_G$ agrees with the Homfly polynomial $P^{\Sigma_G}_{L_{D_G}}$ (in the thickened surface $\Sigma_G \times I$).
\end{coroll}

To prove Theorem \ref{thm:abstract_homfly_well_def_elosztott_graf}, we will rely on a result of P\'ah\'an and the second author.
In \cite{pahan}, they give an embedding activity definition for $Z_G$. We will show that the terms of their formula correspond exactly to the terms in $P_{D(G),v_0,e_0}$. (In fact, the results of the current paper were the motivation for \cite{pahan}. We also note that the formula of \cite{pahan} is an analogue of the activity generating function definitions for $Z_G$ by Morse \cite{Morse_interlace} and Moffatt \cite{moffatt2025constructing}, where a fixed ordering of the edges is used to define the activities, whereas by contrast, \cite{pahan} uses only a fixed half-edge and the ribbon structure for the same purpose.)

Let us recall the formula of \cite{pahan}.
Fix a basepoint for $G$ (that is, an arbitrary vertex and an incident edge).
Let $K\subset G$ be an arbitrary quasi-tree. Similarly to Definitions \ref{def:internal_activity_of_quasi-trees} and \ref{def:external_activity_of_quasi-trees}, we say that an edge $e\in K$ is \emph{internally embedding-active} (with respect to the set of quasi-trees of $G$), if there exists another quasi-tree $L$ of $G$ such that the tours of $K$ and $L$ agree until $e$, but $e\in K-L$. Similarly, $e\notin K$ is \emph{externally embedding-active} if there exists another quasi-tree $L$ of $G$ such that the tours of $K$ and $L$ agree until $e$ and $e\in L-K$.
We denote the number of internally embedding-active edges of the quasi-tree $K$ by $ia^{em}(K)$ and the 
number of externally embedding-active edges by $ea^{em}(K)$. We denote the set of quasi-trees of $G$ by $\quasi-trees(G)$.
With that, the authors of \cite{pahan} prove that 
\begin{equation}\label{Z_G_via_embedding_activities}
Z_G(a,b)=b\cdot \sum_{K\in\quasi-trees(G)} a^{|K|}\,(1+a^{-1}b)^{ia^{em}(K)}\,(1+ab)^{ea^{em}(K)}.   
\end{equation}

To prove Theorem \ref{thm:abstract_homfly_well_def_elosztott_graf}, we compare quasi-trees of $G$ to Jaeger quasi-trees of $D_G$.
We will use the following notations: $V(D_G)=\{v_x\mid x\in V(G)\cup E(G)\}$ and $E(D_G)=\{\overrightarrow{v_uv_{uw}}\mid uw\in E(G)\}$. We also denote the edge $\overrightarrow{v_uv_{uw}}$ by $e_{u,uw}$.
Notice that if $F\subseteq E(D_G)$ is a spanning subgraph of $D_G$, then for each edge $uw\in E(G)$ the set $F\cap \{e_{u,uw}, e_{w,uw}\}$ consists of either 1 or 2 edges. 

\begin{defn}
    Let $F\subset E(D_G)$ be a spanning subgraph of $D_G$. We define 
    $$\pi(F)=\{uw\in E(G) \mid |F\cap \{e_{u,uw}, e_{w,uw}\}|=2\}.$$
\end{defn}

We have the following simple relationship between the quasi-trees of $G$ and the Jaeger quasi-trees of $D_G$:

\begin{prop}\label{prop:quasi-trees_of_G_and_DG}
    If $J$ is a Jaeger quasi-tree of $D_G$ (with respect to some fixed basepoint), then $\pi(J)$ is a quasi-tree of $G$.
    Conversely, for each quasi-tree $K$ of $G$ there is a unique Jaeger quasi-tree $J\in\jaeg(D_G)$ such that $\pi(J)=K$. 
\end{prop}

\begin{proof}
    Let $L$ be an arbitrary quasi-tree of $D_G$. Remove from $L$ the edges $e_{u,uw}\in L$ such that $e_{w,uw}\notin L$, and call the resulting (non-spanning) subgraph $L'$. As $d_L(v_{uw})=1$, the subgraph $L'$ still has one boundary component. The boundary of $\pi(L)$ within $G$ essentially agrees with the boundary of $L'$, moreover $\pi(L)$ is a spanning subgraph in $G$, whence $\pi(L)$ is a quasi-tree of $G$.

    Conversely, if we take $K':=\{e_{u,uw}, e_{w,uw}\mid uw\in K\}$ for some quasi-tree $K$ of $G$, then the boundary of $K'$ is essentially the same as the boundary of $K$, 
    in particular it is connected.
    If we add one of $e_{u,uw}$ and $e_{w,uw}$, for any edge $uw\notin K$, to $K'$ then we do not increase the number of boundary components. Hence the so obtained spanning subgraph $L$ is a quasi-tree of $D_G$, moreover $\pi(L)=K$.

    We need to show that there is exactly one way to extend $K'$ to $L$ so that 
    $L$ is a Jaeger quasi-tree.

    Let $uw\notin K$ and suppose that $(u,uw) <_K (w,uw)$. Then, let $e_{u,uw}\notin L$ and $e_{w,uw}\in L$. If we make this choice for each $uw\notin K$, then clearly $L$ is naive Jaeger and in all other cases it is not.
    By Proposition \ref{prop:naive_Jaeger_in_DG_is_Jaeger}, in $D_G$ naive Jaeger 
    and Jaeger quasi-trees coincide, so the proof is complete.
\end{proof}

Using the above correspondence between Jaeger quasi-trees of $D_G$ and quasi-trees of $G$, we can re-write the abstract Homfly polynomial of $D_G$ (with an arbitrary root) 
as a sum whose terms are labeled by 
the quasi-trees of $G$. Luckily, we can also translate the activities of Jaeger quasi-trees of $D_G$ to embedding activities of the quasi-trees of $G$.

\begin{prop}\label{prop:activities_in_K_and_pi(K)}
    Let $K$ be an arbitrary Jaeger quasi-tree of $D_G$.
    Let $uw\in E(G)$ be such that $e_{w,uw} <_K e_{u,uw}$. Then $e_{u,uw}\in K$, and it is is internally active in $K$.
    Furthermore, $e_{w,uw}\in K$ if and only if $uw\in \pi(K)$, and $e_{w,uw}$ is internally/externally active/passive in $K$ if and only if $uw$ is internally/externally active/passive in $\pi(K)$.
\end{prop}

\begin{proof}
We start with considering $e_{u,uw}$.
If $uw\in\pi(K)$ then clearly $e_{u,uw}\in K$. Since $e_{w,uw} <_K e_{u,uw}$, in this case $e_{u,uw}$ is reached at its head $v_{uw}$, which makes $e_{u,uw}$ tight of type~$0$ and thus internally active in $K$ by Theorem \ref{thm:internal_activity_condition}.
If $uw\notin \pi(K)$, then exactly one of $e_{u,uw}$ and $e_{w,uw}$ is in $K$ and, by the argument in the proof of Proposition \ref{prop:quasi-trees_of_G_and_DG}, this needs to be $e_{u,uw}$.
Since $e_{w,uw}\notin K$, the edge $e_{u,uw}$ is a bridge in  $\AuxContr{D_G,K,e_{u,uw}}$, whence $e_{u,uw}$ is tight of type~$1$ and internally active in $K$.

Now let us turn to $e_{w,uw}$. From the fact that $e_{u,uw}\in K$ always holds, it follows that $e_{w,uw}\in K$ if and only if $uw\in \pi(K)$.

Consider the tour of $\pi(K)$ in $G$. We identify the edge $e_{u,uw}$ of $D_G$ with the half-edge $(u,uw)$ of $G$. One can check that with this identification, the ordering $<_K$ on $E(D_G)$ agrees with the ordering $<_{\pi(K)}$ on the half-edges of $G$.
Since $e_{w,uw} <_K e_{u,uw}$, the position of $e_{w,uw}$ (in $<_K$) determines the position of $uw$ (in $<_{\pi(K)}$). 

If $e_{w,uw}\notin K$ is a first difference between $K$ and another Jaeger quasi-tree $L$ of $D_G$, then $e_{w,uw}\in L$ and $e_{w,uw} <_L e_{u,uw}$, wherefore also $e_{u,uw}\in L$. Hence $uw\in \pi(L)$ and it is a first difference between the quasi-trees $\pi(K)$ and $\pi(L)$. Similarly, if $e_{w,uw}\in K$ is a first difference between $K$ and another Jaeger quasi-tree $L$, then $uw\in \pi(K)-\pi(L)$ is a first difference between $\pi(K)$ and $\pi(L)$.

Conversely, if $uw\notin \pi(K)$ is a first difference between $\pi(K)$ and another quasi-tree of $G$ and $(w,uw) <_{\pi(K)} (u,uw)$, then by Proposition \ref{prop:quasi-trees_of_G_and_DG} that other quasi-tree has the form $\pi(L)$ for a Jaeger quasi-tree $L$ of $D_G$. Then, $uw\in\pi(L)$ implies $e_{w,uw}\in L$. Hence $e_{w,uw}\in L-K$ is a first difference between $K$ and $L$. If $uw\in \pi(K)$ is a first difference between $\pi(K)$ and another quasi-tree of $G$, than that other quasi-tree again has the form $\pi(L)$ for a Jaeger quasi-tree $L$ of $G$. 
Since the tours of $\pi(K)$ and $\pi(L)$ agree until the half-edge $(w,uw)$, we also have $e_{w,uw} <_L e_{u,uw}$ and thus $e_{u,uw}\in L$. Hence $uw\notin \pi(L)$ means $e_{w,uw}\notin L$, and $e_{w,uw}$ is a first difference between $K$ and $L$. 
\end{proof}

\begin{proof}[Proof of Theorem \ref{thm:abstract_homfly_well_def_elosztott_graf}]
    Since each edge in $D_G$ is positive by definition, we have 
    $$P_{D_G,v_0,e_0}(v,z)=\sum_{J\in\jaeg(G,v_0,e_0)} (v^2)^{ip^+(J)}\cdot (vz)^{ep^+(J)}\cdot (vz+vz^{-1}-v^3z^{-1})^{ea^+(J)}.$$
    By Propositions \ref{prop:quasi-trees_of_G_and_DG} and \ref{prop:activities_in_K_and_pi(K)}, this is further equal to $$\sum_{K \in \quasi-trees(G)} (v^2)^{ip^{em}(K)}\cdot (vz)^{ep^{em}(K)}\cdot (vz+vz^{-1}-v^3z^{-1})^{ea^{em}(K)},$$
    where $ip^{em}(K)$ ($ep^{em}(K)$) denotes the number of edges of $K$ (of $G-K$) that are not internally (externally) embedding-active. 
    By rearranging the terms, $P_{D_G,v_0,e_0}(v,z)$ is further equal to  
    $$(vz)^{|E(G)|}\sum_{K \in \quasi-trees(G)} \left(\frac{v}{z}\right)^{|K|}\cdot (v^{-2})^{ia^{em}(K)}\cdot (1+z^{-2}-v^2z^{-2})^{ea^{em}(K)}.$$
    Now we can apply \eqref{Z_G_via_embedding_activities} with $a=v/z$ and $b=(v^{-1}-v)/z$ to conclude the statement of the Theorem.
\end{proof}

\subsection{The abstract Homfly polynomial can depend on the root}
\label{ss:root_dependence_of_abstract_homfly}

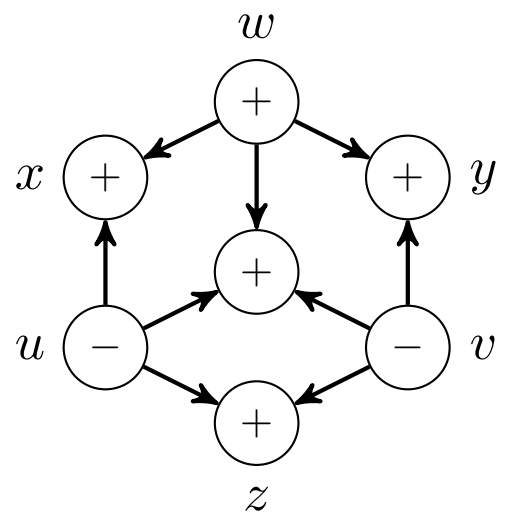
\begin{figure}[t]
	\begin{tikzpicture}[scale=.25,>=stealth']
        \tikzstyle{c}=[circle,draw,scale=0.8]
			\node[c] (a) at (14, 6.5) {$+$};
			\node[c] (b) at (22, 6.5) {$+$};
			\node[c] (c) at (18, 4) {$+$};
			\node[c] (d)  at (18, 0) {$+$};
			\node[c] (e) at (14, 2) {$-$};
			\node[c] (f) at (22, 2) {$-$};
			\node[c] (g) at (18, 8.5) {$+$};
			
			\node at (12, 6.5) {$x$};
			\node at (24, 6.5) {$y$};
			\node at (18, -2) {$z$};
			\node at (12, 2) {$u$};
			\node at (24, 2) {$v$};
			\node at (18, 10.5) {$w$};
            
            \draw [thick,->] (g) -- (a);
			\draw [thick,<-] (a) -- (e);
			\draw [thick,->] (g) -- (b);
			\draw [thick,->] (g) -- (c);
			\draw [thick,->] (e) -- (c);
			\draw [thick,->] (e) -- (d);
			\draw [thick,<-] (d) -- (f);
			\draw [thick,->] (f) -- (c);
			\draw [thick,->] (f) -- (b);
			
	\end{tikzpicture}
	\caption{A ribbon digraph where the abstract Homfly polynomial depends on the basepoint.}\label{fig:basepoint_dependence}
\end{figure}

Consider the directed graph of Figure \ref{fig:basepoint_dependence}. We endow it with the ribbon structure that is given by the positive orientation of the plane around the vertices with a $+$ sign, and by the negative orientation of the plane around the vertices with a $-$ sign. Moreover, we take each edge to be positive. (We note that at vertices of degree $2$ the signs $+$ and $-$ are equivalent because a two-element set has a unique cyclic order.)

This signed ribbon digraph is the Seifert graph of a link in the following thickened surface: We need to take the canonical surface for $G$ with the ribbon structure where the cyclic ordering is induced by the positive orientation of the plane around each vertex except for $w$, where it is induced by the negative orientation.

The abstract Homfly polynomial of this ribbon graph for the root $(w,wx)$ is 
$$
	\begin{array}{llll}
		v^3z^3 & +3v^5z^3 & +3v^7z^3  & \\  
		& +3v^5z &  +3v^7z&  -3v^9z  \\  
		&  +v^5z^{-1}&  -v^7z^{-1}&   
	\end{array},
$$
while the abstract Homfly polynomial for the basis $(v,vz)$ is 
$$
	\begin{array}{lllll}
		v^3z^3 & +3v^5z^3 & +3v^7z^3  & &\\  
		& +3v^5z &  +4v^7z&  -4v^9z & \\  
		&  &  +2^7z^{-1}& -3v^9z^{-1} & +v^{11}z^{-1}   
	\end{array}.
$$

One may also wonder if Corollary \ref{cor:conway} can be the basis for a notion of Conway polynomial of abstract ribbon digraphs. Namely, one may hope that the numbers of Jaeger quasi-trees of various genera (i.e., number of edges) in certain ribbon digraphs, where we take all edges to be positive, do not depend on the choice of basepoint. This is true for the example presented in this section (indeed, at $v=1$ both abstract Homfly polynomials become $7z^3+3z$, meaning that there are $7$ Jaeger trees and $3$ Jaeger quasi-trees of genus $1$, and this remains true for all other basepoints, too) but unfortunately, it fails in general. Counterexamples are provided by Nagashima \cite{nagashima}.

\section{Jaeger quasi-trees and naive Jaeger quasi-trees}
\label{s:Jaeger_and_naive_Jaeger}

In this section we highlight some cases where considering all naive Jaeger quasi-trees in place of just Jaeger quasi-trees is either sufficient, or definitely insufficient. In this we will only rely on the results of Sections \ref{s:intro}, \ref{s:prelim}, and \ref{s:Jaeger_trees_def_and_formula}.

Checking if a quasi-tree is naive Jaeger is more straightforward than checking whether it is Jaeger, which is why it would be desirable if one could use naive Jaeger quasi-trees instead of Jaeger quasi-trees.
However, in Section \ref{ssec:naive_Jaeger_trees_are_not_all_trivial} we show by an example that our formula from Theorem \ref{thm:Homfly_formula} would not work if we used naive Jaeger quasi-trees instead of Jaeger quasi-trees.

On the other hand, there are some natural cases where the naive Jaeger property  implies the Jaeger property, that is, we can rely on the simpler notion. 
In Section \ref{ssec:naive_J_trees_are_typically_J_trees} we show that if the Seifert graph is acyclic, then its set of Jaeger \emph{trees} (i.e., minimal genus Jaeger quasi-trees) agrees with that of naive Jaeger trees.
In Sections \ref{ssec:special_alternating_links} and \ref{ss:edge_subdivided_graphs_naive_Jaeger_is_Jaeger} we exhibit two classes of oriented link diagrams whose Seifert graphs are such that all naive Jaeger quasi-trees are actually Jaeger: this is the case for Seifert graphs of special alternating links, and for ribbon digraphs obtained from (undirected) ribbon graphs by subdividing each edge with a new point, and orienting the edges toward the subdividing point.

\subsection{Naive Jaeger quasi-trees do not always yield trivial links}\label{ssec:naive_Jaeger_trees_are_not_all_trivial}

In this section we show that Theorem \ref{thm:Homfly_formula} would not stay true if we substituted ``Jaeger quasi-trees'' with ``naive Jaeger quasi-trees'' in the statement. To explain the reason, first recall that
our strategy to compute the Homfly polynomial was the following: we performed a skein relation if we arrived at a crossing such that 
both the retention and the removal of
the corresponding edge of the Seifert graph 
was compatible with the goal of finding
a Jaeger quasi-tree. (The edge, then, was internally/externally passive in all Jaeger quasi-trees of the two new
graphs). 
Furthermore, if we arrived at a crossing whose corresponding edge was tight of type~$1$, then we contracted a bad separating cycle in the Seifert graph. 
This just meant
``pulling apart'' two parts of the link diagram that had already had a trivial relationship to each other. The final key point was Lemma \ref{l:Jaeger_quasi-trees_link_trivial}, claiming that each branch of the computation leads to a trivial knot. 
In other words, at the moment when the algorithm isolated and recognized a Jaeger quasi-tree, the corresponding diagram was always that of the unknot. It is this last point that we cannot guarantee for all naive Jaeger quasi-trees.

Contracting the nearest bad separating cycles for type~$1$ edges was a convenient technical tool that ensured that 
no later action was taken at type~$2$ edges, 
and it also helped us simplify some proofs, but all arguments could have been modified to omit 
the contracting steps. 
(We could, for example, mark type~$2$ crossings and prohibit applying the skein relation to them.)
Now if we wanted to recast our process as a search for all \emph{naive} Jaeger quasi-trees, then we should refrain from contracting the nearest bad separating cycles for the type~$1$ edges and allow type~$2$ edges to be removed. That way we would end up performing the skein relation at passive edges \emph{with respect to the set of naive Jaeger quasi-trees}. The Seifert graph corresponding to the branch of a Jaeger quasi-tree $K$ would be the quasi-tree 
obtained from $K$ by reversing the internally passive (with respect to naive Jaeger quasi-trees) edges. We show that such a Seifert graph does not necessarily correspond to a trivial knot if $K$ is only a naive Jaeger quasi-tree, in other words the analogue of Lemma \ref{l:Jaeger_quasi-trees_link_trivial} is not true for naive Jaeger quasi-trees.

\begin{prop}\label{prop:naive_Jaeger_can_be_nontrivial}
    There exists a link diagram such that its Seifert graph $G$ has a naive Jaeger quasi-tree $K$ such that if we take $K'$ where we reverse those edges of $K$ that are internally passive with respect to the set of naive Jaeger quasi-trees of $G$, 
    then $K'$ is the Seifert graph of a diagram of a non-trivial knot.
\end{prop}

\begin{figure}[h]
    \centering
    \begin{tikzpicture}[decoration={markings,
    mark= at position 0.6 with {\arrow{stealth}}}]
    \begin{scope}[scale=.3,shift={(0,0)},>=stealth]
    \tikzstyle{e}=[circle,scale=.7,draw,color=e]
    \tikzstyle{v}=[circle,scale=.7,draw,color=blue]
    \node [v] (u) at (10,0) {$u$};
    \node [v] (v) at (0,0) {$v$};
    \node [v] (w) at (2,8) {$w$};
    \node [e] (x) at (5,4) {$x$};
    \node [e] (y) at (5,0) {$y$};
    \node [e] (z) at (5,-4) {$z$};
    \node [e] (t) at (-3,4) {$t$};
    \draw [<-,thick] (u) to (x);
    \draw [->,thick] (u) to (y);
    \draw [->,thick] (u) to (z);
    \draw [->,thick] (v) to (x);
    \draw [<-,thick] (v) to (y);
    \draw [->,thick] (v) to (z);
    \draw [->,thick] (v) to (t);
    \draw [->,thick] (w) to (t);
    \draw [<-,thick] (w) to (x);
    \draw [thick,postaction={decorate},color=gray] (11.5,0) to [out=90,in=-10] (7.2,2.2) to [out=170,in=0] (5,2.5) to [out=180,in=10] (2.6,2.05) to [out=-170,in=-10] (-1.25,2.1);
    \draw [thick,postaction={decorate},color=gray] (-2,2.2) to [out=170,in=-90] (-4.2,4) to [out=90,in=180] (-0.5,6) to [out=0,in=180] (2,6) to [out=0,in=170] (4.2,5.65) to [out=-10,in=140] (6.3,4.8) to [out=-40,in=110] (7.2,2.6);
    \draw [thick,color=gray] (7.3,1.8) to [out=-80,in=80] (7.3,0.3);
    \draw [thick,postaction={decorate},color=gray] (7.2,-0.3) to [out=-110,in=-70] (2.85,-0.3);
    \draw [thick,color=gray] (2.7,0.3) to [out=110,in=-100] (2.55,1.8);
    \draw [thick,postaction={decorate},color=gray] (2.6,2.3) to [out=80,in=-100] (3.6,5.5);
    \draw [thick,postaction={decorate},color=gray] (3.7,6) to [out=80,in=0] (2,9.5) to [out=180,in=90] (-0.5,6.3);
    \draw [thick,postaction={decorate},color=gray] (-0.5,5.5) to [out=-90,in=70] (-1.6,2) to [out=-110,in=120] (-1.5,-0.3) to [out=-50,in=-180] (2.8,-2.6);
    \draw [thick,color=gray] (3.5,-2.55) to [out=10,in=170] (6.8,-2.4) to [out=-10,in=-90] (11.5,0);
    \draw [thick,color=gray] (7.05,-2.1) to [out=90,in=-60] (7.25,0) to [out=120,in=0] (5,1.2) to [out=180,in=70] (2.7,0) to [out=-110,in=100] (3.1,-2.4) to [out=-80,in=180] (5,-5.2) to [out=0,in=-90] (7.05,-2.75);
    \draw[->,thick,color=gray] (4.9, 1.2) to (5.1, 1.2);
    \node[circle,draw,fill,scale=.3] at (11.5,0) {};
	\draw[->,thick] (11.5, 0) to (11.5, 1);
    \end{scope}
    \begin{scope}[scale=.3,shift={(20,0)},>=stealth]
    \tikzstyle{e}=[circle,scale=.7,draw,color=e]
    \tikzstyle{v}=[circle,scale=.7,draw,color=blue]
    \node [v] (u) at (10,0) {$u$};
    \node [v] (v) at (0,0) {$v$};
    \node [v] (w) at (2,8) {$w$};
    \node [e] (x) at (5,4) {$x$};
    \node [e] (y) at (5,0) {$y$};
    \node [e] (z) at (5,-4) {$z$};
    \node [e] (t) at (-3,4) {$t$};
    \draw [<-,thick] (u) to (x);
    \draw [->,thick,dotted] (u) to (y);
    \draw [->,thick] (u) to (z);
    \draw [->,thick] (v) to (x);
    \draw [<-,thick] (v) to (y);
    \draw [->,thick] (v) to (z);
    \draw [->,thick] (v) to (t);
    \draw [->,thick] (w) to (t);
    \draw [<-,thick] (w) to (x);
    \draw [thick,postaction={decorate},color=gray] (11.5,0) to [out=90,in=-10] (7.2,2.2) to [out=170,in=0] (5,2.5) to [out=180,in=10] (2.6,2.05) to [out=-170,in=-10] (-1.25,2.1);
    \draw [thick,postaction={decorate},color=gray] (-2,2.2) to [out=170,in=-90] (-4.2,4) to [out=90,in=180] (-0.5,6) to [out=0,in=180] (2,6) to [out=0,in=170] (4.2,5.65) to [out=-10,in=140] (6.3,4.8) to [out=-40,in=110] (7.2,2.6);
    \draw [thick,color=gray] (7.3,1.8) to [out=-80,in=80] (7.9,0) to [out=-80,in=90] (7,-2.1);
    
    \draw [thick,color=gray] (7.1,-2.8) to [out=-90,in=0] (5,-5.2) to [out=180,in=-90] (3.1, -2.4) to [out=110,in=-110] (2.7,0) to [out=70,in=-180] (5,1.2) to [out=0,in=90] (7.1,0) to [out=-110,in=-70] (2.85,-0.3);
    \draw [thick,color=gray] (2.7,0.3) to [out=110,in=-100] (2.55,1.8);
    \draw [thick,postaction={decorate},color=gray] (2.6,2.3) to [out=80,in=-100] (3.6,5.5);
    \draw [thick,postaction={decorate},color=gray] (3.7,6) to [out=80,in=0] (2,9.5) to [out=180,in=90] (-0.5,6.3);
    \draw [thick,postaction={decorate},color=gray] (-0.5,5.5) to [out=-90,in=70] (-1.6,2) to [out=-110,in=120] (-1.5,-0.3) to [out=-50,in=-180] (2.8,-2.6);
    \draw [thick,color=gray] (3.5,-2.55) to [out=10,in=170] (6.8,-2.4) to [out=-10,in=-90] (11.5,0);
    \draw[->,thick,color=gray] (4.9, 1.2) to (5.1, 1.2);
    \node[circle,draw,fill,scale=.3] at (11.5,0) {};
	\draw[->,thick] (11.5, 0) to (11.5, 1);
    \end{scope}
    \end{tikzpicture}
    \caption{
    A naive Jaeger quasi-tree that is associated with a non-trivial knot.}
    \label{fig:non_Jaeger_example}
\end{figure}
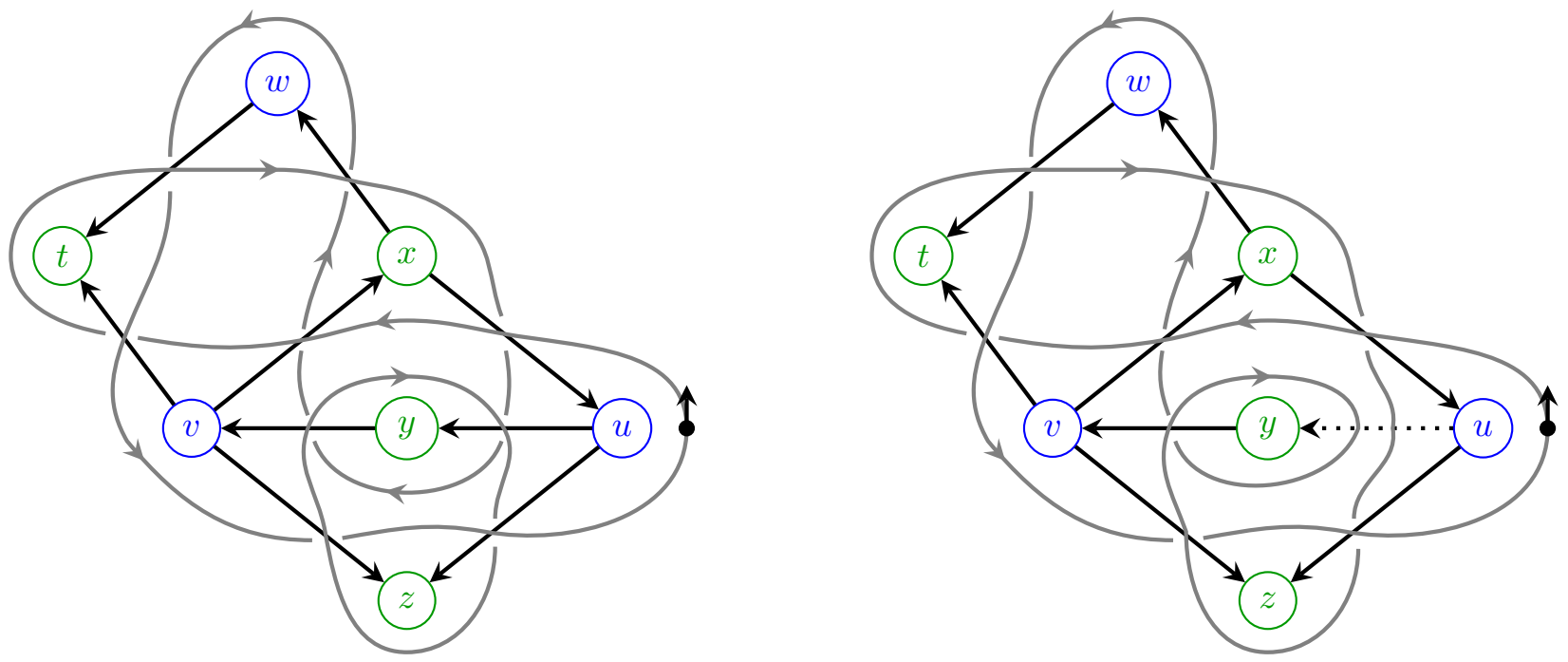

\begin{proof}
    Take the Seifert graph $G$ of the link diagram of Figure \ref{fig:non_Jaeger_example}. (This is in fact a complicated diagram of a Hopf link.) The ribbon structure is given by the positive orientation of the plane around the blue vertices ($u$, $v$, and $w$), and by the negative orientation of the plane around the green vertices ($x$, $y$, $z$, and $t$).
    Choose $(u,ux)$ as the root.

    Then, $K=G-uy$ is a naive Jaeger quasi-tree that is not a Jaeger quasi-tree.

    To see that $K-uy$ is a naive Jaeger quasi-tree, note that $<_K$ orders the half-edges as: $(u,ux)$, $(x,xv)$, $(v,vt)$, $(t,tw)$, $(w,wx)$, $(x,xu)$, $(u,uy)$, $(u,uz)$, $(z,zv)$, $(v,vy)$, $(y,yu)$, $(y,yv)$, $(v,vx)$, $(x,xw)$, $(w,wt)$, $(t,tv)$, $(v,vz)$, $(z,zu)$. For the only non-edge of $K$ we have $(u,uy) <_K (y,uy)$, where $u$ is the tail of $uy$, whence $K$ is indeed a naive Jaeger quasi-tree.

    To see that $K$ is not a Jaeger quasi-tree, notice that since $u$ is the head of $ux$ and $x$ is the head of $xv$, we have $\AuxContr{G,K,vt}=G$. However, for $(v,vt)$, the vertex $v$ is the tail, moreover, $G-vt$ has $\{ux,xv,vy,yu\}$ as its nearest bad separating cycle. (See also Figure \ref{fig:Jaeger_elimination_example} and Example \ref{ex:Jaeger_test}.)
    Since $uy\notin K$, the edge $uy$ fails the Jaeger test (with case \ref{e_notin_K_contracted}). 
    
    One can check that each edge of $K$ is internally active with respect to the set of naive Jaeger quasi-trees, implying $K'=K$, and that the link corresponding to $K'$, see the right panel of Figure \ref{fig:non_Jaeger_example}, is the negative trefoil knot.
\end{proof}

\subsection{Jaeger trees and naive Jaeger trees coincide for acyclic ribbon digraphs}\label{ssec:naive_J_trees_are_typically_J_trees}

Recall that we call a Jaeger quasi-tree $K$ a Jaeger tree if $K$ is a tree.
Similarly, we call a naive Jaeger quasi-tree $K$ a \emph{naive Jaeger tree}, if $K$ is a tree.

We denote the set of Jaeger trees of $G$ with root $(v_0,e_0)$ by $\jaegtree(G,v_0,e_0)$, and the set of naive Jaeger trees of $G$ with root $(v_0,e_0)$ by $\naivejaegtree(G,v_0,e_0)$.

\begin{prop}\label{prop:naive_Jaeger_tree_is_Jaeger_tree}
    If $G$ is an acyclic connected rooted ribbon digraph, then each naive Jaeger tree of $G$ is a Jaeger tree.
\end{prop}

\begin{remark}
    As pointed out in Lemma \ref{lem:link_has_diagram_with_connected_acyclic_seifert}, using Reidemeister II moves, each link diagram can be modified so that the Seifert graph becomes connected and acyclic. Thus, the proposition can be applied to appropriately chosen diagrams of any link.
\end{remark}

\begin{proof}
    Suppose that $G$ contains a naive Jaeger tree $T$ that is not a Jaeger tree. As the naive Jaeger property prevents $T$ from failing the Jaeger test via case \ref{e_notin_K_head}, it has to fail the test via case \ref{e_notin_K_contracted} for some edge $f$. 
    This requires that some edges be contracted during the Jaeger test of $T$, in other words, certain other edges have to pass the Jaeger test via case \ref{e_in_K_tail_active} in such a way that the contracted cycle is non-empty. Let $e$ be the first such edge according to $<_T$. 
    This implies that the auxiliary graph $\AuxContr{G,T,e}$ is formed from $G$ without contracting any edges. I.e., it is a subgraph of $G$ where the orientations of some edges are reversed, moreover, each edge $g$ that is oriented differently in $G$ and in $\AuxContr{G,T,e}$ satisfies $g <_T e$. 
    Let us use the shorthand $G'=\AuxContr{G,T,e}$.

    Let $C$ be the nearest bad separating cycle of $G'-e$, which is non-empty by assumption and directed by 
    definition. 
    As $G$ is acyclic, this is only possible if at least one edge $f\in C$ has been reversed where, then, $f<_T e$. 

    The subgraph $T-e$ has two connected components. Let $V_0$ be the vertex set of the connected component containing the root vertex $v_0$, and let $V_1=V(G)-V_0$. Note that $V_0$ is the set of vertices reached either before $(t(e),e)$ or after $(h(e),e)$ in 
    the tour of $T$,
    while $V_1$ is the set of vertices reached between $(t(e),e)$ and $(h(e),e)$.
    As we have some $f\in C$ with $f<_T e$, 
    the vertex set of $C$ needs to intersect $V_0$.
    
    From the fact that $C$ is a separating cycle in $G'-e$ but not in $G'$, we infer that the two half-edges of $e$ interlace with the two half-edges of $g$ for each edge $g\in C$. Moreover, because $C$ is also ``bad,'' the half-edge $(h(g),g)$ is either before the first, or after the second traversal of $e$. In other words, each $g\in C$ has either 
    \begin{equation}\label{eq:fund1}
    (h_{G'}(g),g)<_T (t_{G'}(e),e) <_T (t_{G'}(g),g)<_T (h_{G'}(e),e)
    \end{equation}
    or
    \begin{equation}\label{eq:fund2}
    (t_{G'}(e),e) <_T (t_{G'}(g),g)<_T (h_{G'}(e),e) <_T (h_{G'}(g),g).
    \end{equation}
    If $(x,xy)$ is a half-edge with $x\in V_0$ then, since $T$ is a tree, we either have $(x,xy) <_T (t_{G'}(e),e)$ or $(h_{G'}(e),e) <_T (x,xy)$. 
    Comparing this to \eqref{eq:fund1} and \eqref{eq:fund2}, we conclude that $C$ cannot have any edge 
    with both endpoints in $V_0$. 

    Finally, consider any edge $xy\in C$ with $x\in V_0$ and $y\in V_1$. (By the paragraph above, such edges exist.)
    We have $xy\notin T$ since $e\notin C$ is the only edge of $T$ connecting $V_0$ and $V_1$. Since $C\subseteq \AuxContr{G,T,e}$, it has to be that $xy >_T e$, for otherwise $xy$ would 
    be deleted during the formation of $\AuxContr{G,T,e}$. 
    This means that we are in the situation of \eqref{eq:fund2} with $x=h_{G'}(xy)$. In other words, any $xy\in C$ with $x\in V_0$ and $y\in V_1$ is directed toward $x$ in $G'$. 
    This contradicts the fact that $C$ is a directed cycle in $G'$.
\end{proof}

\subsection{Special alternating links}
\label{ssec:special_alternating_links}

A \emph{special alternating diagram} is an oriented link diagram 
that is alternating with
crossings of uniform sign, which we will assume to be positive. Such a diagram only has type~$1$ Seifert circles.
Let us also suppose that the diagram is connected and that its `point at infinity' is chosen so that none of the Seifert circles separate it from the rest of the diagram.

A link is \emph{special alternating} if it possesses a special alternating diagram. 
Trying to understand just the `top' (maximum $z$-degree) part of the Homfly polynomial of a special alternating link has already lead to advances in combinatorics by the authors and others \cite{hiperTutte,KP_Ehrhart,KM,semibalanced}.
Much of the current paper grew out of attempts to 
find graph-theoretic descriptions of further coefficients of these Homfly polynomials.

Let us color a vertex of the Seifert graph blue if the corresponding Seifert circuit has a positive (counterclockwise) orientation, and color a vertex green otherwise. Then in our cases of interest, we get a Seifert graph of the following structure.

\begin{fact}\label{fact:Seifert_graph_of_spec_alt_diagram}
    The Seifert graph of a special alternating diagram is a plane-embedded bipartite graph with blue and green vertices, such that each edge is positive, and is oriented from the green endpoint to the blue endpoint. Furthermore, the ribbon structure of the Seifert graph is induced by the positive orientation of the plane around blue vertices, and by the negative orientation around green vertices.
\end{fact}

We use $\blue(e)$ and $\green(e)$ to refer to the blue and green endpoints, respectively, of an edge $e$. Since each edge is positive, we initially have $\blue(e)=h(e)$ and $\green(e)=t(e)$, however in the auxiliary graphs $\AuxContr{G,K,e}$ this might change. 
This explains the need for extra notation.

The next proposition claims that the computation of the Homfly polynomial of a special alternating link is particularly easy to track combinatorially, because the distinction between `Jaeger' and `naive Jaeger' quasi-trees disappears.

\begin{prop}\label{prop:spec_alt_naive_Jaeger_is_Jaeger}     If $G$ is the Seifert graph of a special alternating diagram, 
    then each naive Jaeger quasi-tree of $G$ is a Jaeger quasi-tree.
\end{prop}

\begin{proof}
Let $(v_0,e_0)$ be an arbitrarily fixed root for $G$. 
Suppose for a contradiction that some naive Jaeger quasi-tree $K$ of $G$ fails the Jaeger test. Since $K$ is naive Jaeger, it cannot fail the test via case \ref{e_notin_K_head}, i.e., it fails the test via case \ref{e_notin_K_contracted}. 
Let $e$ be the first edge according to $<_K$ such that $e$ is tight of type~$1$, and the nearest bad separating cycle $C$ of $\AuxContr{G,K,e}-e$ has an edge $f\notin K$. 
Recall that $C$ is a directed cycle in $\AuxContr{G,K,e}$.

Let $G'$ be the graph obtained from $\AuxContr{G,K,e}$ by undoing the contraction steps. In other words, $G'$ is derived from $G$ by deleting edges $g$ such that $g\notin K$ and $g<_K e$, and reversing edges $g\in K$ with $g<_K e$ if they are internally passive in $K$.
Note that $G'$ is still a Seifert graph of a link diagram $D'$. In particular $G'$ is embedded into the plane, and its 
ribbon structure is given by the positive orientation of the plane around blue vertices, and by the negative orientation around green vertices. 
This is simply because $D'$ may be obtained from our original special alternating diagram $D$
by smoothing the crossings corresponding to $g<_K e$ with $g\notin K$ and switching the under- and overstrands at crossings that correspond to $g <_K e$ where the edge $g$ is internally passive in $K$. We can recognize the edges $r$ that have been reversed by the property that they have $\blue(r)=t_{G'}(r)$. 

In $G'-e$, some separating cycle $C'=C\cup S$ corresponds to $C$, where the edges of $S$ were contracted when we formed $\AuxContr{G,K,e}-e$. 
This claim is more easily verified if we examine the dual picture: the dual $(\AuxContr{G,K,e}-e)^*$ can be obtained from $(G'-e)^*$ by deleting some edges, in particular the vertex sets of the two dual graphs agree. We know that $C^*$ is an elementary cut in $(\AuxContr{G,K,e}-e)^*$. Let $B_0$ and $B_1$ be the vertex sets of $(\AuxContr{G,K,e}-e)^*$ such that $C^*$ is the set of edges connecting $B_0$ and $B_1$, with $B_0$ containing the base vertex. 
Then $(C')^*$ is the set of edges of $(G'-e)^*$ connecting a vertex of $B_0$ and a vertex of $B_1$. As $(G'-e)^*$ can be obtained from $(\AuxContr{G,K,e}-e)^*$ by adding some edges, $(C')^*$ is still an elementary cut in $(G'-e)^*$ (that is, $B_0$ and $B_1$ still span connected subgraphs in $(G'-e)^*$). Therefore $C'$ is indeed a separating cycle and $S=C'-C$ is a set of edges that we contracted when we constructed $\AuxContr{G,K,e}$. 

By smoothing the crossing of $e$ in $D'$ we obtain a link diagram $D''$ with Seifert graph $G'-e$. The components of the link depicted by $D''$ can be identified with the dual vertices $B_0\cup B_1$. We think of each component of $D'$ (and of $D''$) as a cyclically ordered set of half-edges of $G'$ (or of $G'-e$). Recall that these half-edges correspond to an under- or an overstrand at a crossing.

In our argument, it will be important that $G'$ is embedded not only in its `abstract' canonical surface $\Sigma_{G'}$, but also in the plane.
Take a look at $C'$. As a graph, it is a union of edge-disjoint cycles. The usual checkerboard coloring divides the plane into two (not necessarily connected) regions, such that we need to cross $C'$ in order to get from one region to the other. For simplicity, we will call the two regions the \emph{interior} of $C'$ and the \emph{exterior} of $C'$. See Figure \ref{fig:spec_alt}, where the shaded region is the interior.

The graph $C'$ also inherits a ribbon structure from $G'$ (or $G$), which is still given by the positive orientation around blue vertices and the opposite around green ones.
If we follow a boundary component of $C'$, its half-edges alternately cross from the outside to the inside of $C'$, and from the inside to the outside. See again Figure \ref{fig:spec_alt}.

\begin{figure}[h]
    \centering
    \begin{tikzpicture}[decoration={markings,
    mark= at position 0.6 with {\arrow{stealth}}}]
    \begin{scope}[scale=.3,shift={(0,0)},>=stealth]
    \tikzstyle{e}=[circle,scale=.7,draw,color=e,fill=e]
    \tikzstyle{v}=[circle,scale=.7,draw,color=blue,fill=blue]
    \draw[fill=light-gray] (10,0) -- (8,5) -- (2,8) -- (-3,4)--(-2,-2)--(5,-5)--(10,0);
    \draw[fill=white] (10,0)--(5,-1)--(1,1)--(5,3)--(10,0);
    \node [v] (u) at (10,0) {};
    \node [v] (v) at (-2,-2) {};
    \node [v] (v') at (1,1) {};
    \node [v] (w) at (2,8) {};
    \node [e] (x) at (5,3) {};
    \node [e] (x') at (8,5) {};
    \node [e] (y) at (5,-1) {};
    \node [e] (z) at (5,-5) {};
    \node [e] (t) at (-3,4) {};
    \draw [<-,thick] (u) to (x);
    \draw [<-,thick] (u) to (x');
    \draw [<-,thick] (u) to (y);
    \draw [<-,thick] (u) to (z);
    \draw [<-,thick] (v') to (x);
    \draw [<-,thick] (v') to (y);
    \draw [<-,thick] (v) to (z);
    \draw [<-,thick] (v) to (t);
    \draw [<-,thick] (w) to (t);
    \draw [<-,thick] (w) to (x');
    \draw[thick,color=gray,rounded corners=20pt,postaction={decorate}] (7.9, -2.5) -- (12.5, 0) -- (7, 4) -- (5.2,6.2);
    \draw[thick,color=gray,rounded corners=20pt,postaction={decorate}] (4.8, 6.8) -- (1.8, 10) -- (-0.5, 6) -- (-2.2,1.5);
    \draw[thick,color=gray,rounded corners=20pt,postaction={decorate}] (-2.8, 0.5) -- (-4, -4) -- (4.5, -3.5) -- (6.9,-2.8);
    \draw[thick,color=gray,postaction={decorate}] (0, 6) to[out=-10,in=-130] (5, 6.5) to[out=60,in=100] (9.5, 6) to [out=-80,in=30] (9.2,3);
    \draw[thick,color=gray,postaction={decorate}] (8.6, 2.5) to [out=-150,in=0] (5.5, 1.8)to [out=-180,in=-20] (3.5, 1.8);
    \draw[thick,color=gray,postaction={decorate}] (2.8, 2.2) to [out=160,in=90] (0, 1) to [out=-90,in=-150] (3,0) to [out=30,in=140] (6.8, -0.2);
    \draw[thick,color=gray,postaction={decorate}] (7.3, -1) to [out=-40,in=80] (7.4, -4) to [out=-100,in=0] (5,-6) to [out=180,in=-50] (1.9,-4);
    \draw[thick,color=gray,postaction={decorate}] (1.6, -3.5) to [out=130,in=-70] (-2.5, 1) to [out=110,in=-110] (-4,4.2) to [out=80,in=170] (-0.7,6.2);
    \draw[thick,color=gray,postaction={decorate}] (6.9, 1.3) to [out=-80,in=60] (7.1, -0.5) to [out=-120,in=0] (5,-1.8) to [out=180,in=-80] (3,-0.5);
    \draw[thick,color=gray,postaction={decorate}] (2.8, 0.3) to [out=110,in=-110] (3.5, 3) to [out=80,in=180] (5,4.2) to [out=0,in=100] (6.7,2.2);
    \end{scope}
    \end{tikzpicture}
    \caption{A cycle $C'$ and its interior, with edges carrying their original orientations in $G$.}
    \label{fig:spec_alt}
\end{figure}
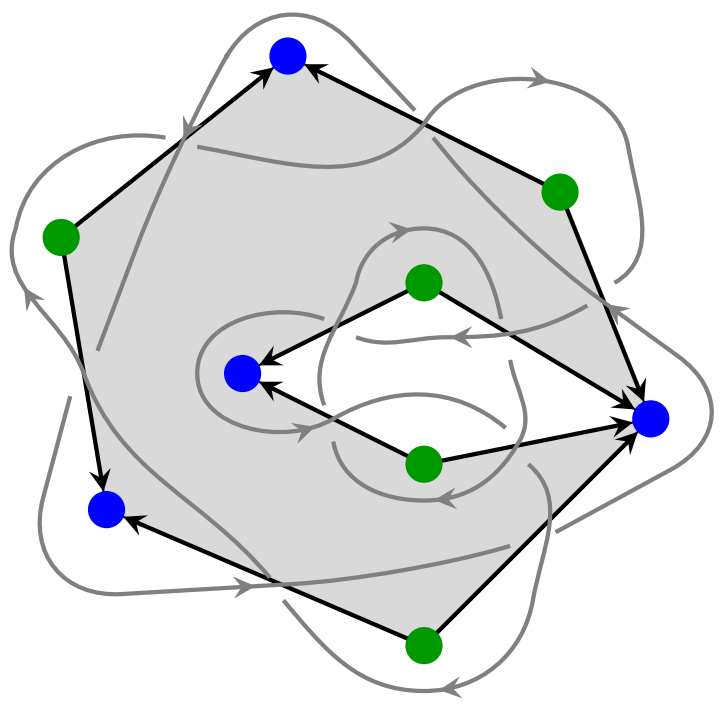

Since $C'$ is a separating cycle in $G'-e$ (for the Seifert ribbon structure), each edge along $C'$ has one of its half-edges as part of an element of $B_0$, while the other half-edge belongs to an element of $B_1$. Furthermore, 
regarding the half-edges `on the $B_0$ side' and their planar projections, by potentially switching the roles of interior and exterior we may suppose that those of the type $(\blue(g),g)$ cross from the exterior of $C'$ to the interior, while half-edges of the type $(\green(g),g)$ cross from the interior of $C'$ to the exterior. 
(Half-edges on the $B_1$ side behave in the opposite manner, cf.\ Figure \ref{fig:spec_alt}.)
Let us use the notation $g\in C'_{\blue}$ if $g\in C'$ with $(\blue(g),g)$ part of an element of $B_0$, and $g\in C'_{\green}$ if $g\in C'$ with $(\green(g),g)$ part of an element of $B_0$. 
This gives us a partition $C'=C'_{\blue}\sqcup C'_{\green}$ with an alternating pattern along $C'$, cf.\ Figure \ref{fig:spec_alt} again. 
Let us also use the notation $C_{\blue}=C'_{\blue} \cap C$ and $C_{\green}=C'_{\green}\cap C$.

As $C$ is a bad separating cycle in $\AuxContr{G,K,e}-e$, in $(\AuxContr{G,K,e}-e)^*$, each edge of $C^*$ points toward $B_0$. This stays true in $(G'-e)^*$ as well, because the vertex sets of the two dual graphs agree. In other words, for each edge of $C$, its half-edge in $B_0$ is the upper strand in $D'$. This means that the edges in $C_{\green}$ changed orientation in $G'$ with respect to $G$. In particular, for each $g\in C_{\green}$, we have $(\green(g),g) <_K (t(e),e)$.

We claim that this is also true for edges of $C'_{\green}\cap S$. That is, 
\begin{equation}\label{eq:C_g_condition}
\forall g\in C'_{\green},  \quad (\green(g),g) <_K (t(e),e).
\end{equation}
Let $g\in C'_{\green}\cap S$ and suppose that $g$ was contracted because of some edge $\tilde{e} <_K e$, which is such that $\AuxContr{G,K,\tilde{e}}-\tilde{e}$ has a nearest bad separating cycle $\tilde{C}$ with $g\in\tilde{C}$. If $(\green(g),g)$ is on the base side of $\tilde C$ with respect to $\AuxContr{G,K,\tilde{e}}-\tilde{e}$, then in order for $\tilde{C}$ to be a bad separating cycle, we need $h_{\AuxContr{G,K,\tilde{e}}}(g)=\green(g)$, that is, $g$ had to be reversed before reaching $\tilde{e}$, whence $(\green(g),g)<_K (t(\tilde{e}),\tilde{e}) <_K (t(e),e)$.
Else if the half-edge $(\blue(g),g)$ 
(which, recall, is part of an element of the non-base side $B_1$ with respect to $e$, in particular satisfies $(t(e),e) <_K (\blue(g),g) < (h(e),e)$) 
belongs to the base component of $(\AuxContr{G,K,\tilde{e}}-\tilde{e})/\tilde C$, then either $(\blue(g),g) <_K (t(\tilde{e}),\tilde{e})$ or $(h(\tilde{e}),\tilde{e}) <_K (\blue(g),g)$. 
Since $\tilde{e}$ is a bridge in $\AuxContr{G,K,e}$, cf.\ Proposition \ref{prop:type1_edge_becomes_bridge}, we know that $e$ and $\tilde{e}$ are not interlacing in the tour of $K$. Using also the fact that $\tilde{e} <_K e$, we conclude that $(t(\tilde{e}),\tilde{e})<_K (h(\tilde{e}),\tilde{e}) <_K (t(e),e) <_K (\blue(g),g) <_K (h(e),e)$. 
Since $g$ lies along the separating cycle $\tilde{C}$ of  $\AuxContr{G,K,\tilde{e}}$, we know that $g$ interlaces with $\tilde{e}$.
This leaves us with $(t(\tilde{e}),\tilde{e})<_K (\green(g),g) <_K (h(\tilde{e}),\tilde{e})$, which implies $(\green(g),g) <_K (t(e),e)$, i.e., that \eqref{eq:C_g_condition} holds in this case, too.

Let $f$ be the first edge according to $<_K$ such that $f\in C-K$. 
Let us follow the tour of $K$, which starts with a half-edge that belongs to an element of $B_0$ and which only traverses edges that are in $G'$. 
As $C'$ is a separating cycle in $G'-e$, the tour can only switch from half-edges along $B_0$ to half-edges along $B_1$ by either traversing $e$, or deleting an edge of $C'$.

Take the moment when $f$ is considered in the Jaeger test. As $f\notin K$, the edge $f$ needs to be reached at its tail, that is, on the half-edge $(\green(f),f)$. Since $g\in K$ for each $g\in C'_{\green}$, we have $f\in C'_{\blue}$, in particular $(\green(f),f)$ lies along $B_1$. Since $f$ is the first edge of $C'-K$, we could only get to $B_1$ through $e$. In other words, in order to reach $f$, one had to traverse $e$ beforehand.  

The half-edge $(\green(f),f)$ crosses from the exterior of $C'$ to the interior.  
Since $f\notin K$, after the half-edge $(\green(f),f)$, the tour of $K$ continues in the interior of $C'$, and with a half-edge of $B_0$. At this point, the tour has yet to visit the other half-edge $(\blue(f),f)$ of $f$. Since $(\green(f),f)$ lies along $ B_1$ and $f\in C'$, we have that $(\blue(f),f)$ is part of an element of $B_0$, as well as that $(\blue(f),f)$ crosses from the exterior of $C'$ to the interior of $C'$. 
Therefore between $(\green(f),f)$ and $(\blue(f),f)$ the tour has to return to the exterior of $C'$. However, this is not possible while staying in $B_0$, since all the half-edges of $B_0$ that cross from the interior of $C'$ to the exterior are in $C'_{\green}$, and as we argued above, all of these half-edges were traversed in the tour of $K$ before $(t(e),e)$, whence also before $(\green(f),f)$.

This means that in order to get to the exterior of $C'$, the tour needs to switch back to $B_1$ at some point between $(\green(f),f)$ and $(\blue(f),f)$. Notice that we cannot switch from $B_0$ to $B_1$ by deleting an edge of $C'$, because in $G'$, the edges of $C'$ all have their half-edge from $B_0$ as the overcrossing strand. 
(In other words, 
if the tour is to reach an edge of $C'$ after $(\green(f),f)$, then that edge will be reached at its head and will not be eligible for deletion in the tour of a naive Jaeger quasi-tree.)
Thus, to get to $B_1$, our only chance is to use $e$ again in the $B_0$ to $B_1$ direction. 
However, $e$ has already been traversed in that direction, before the tour even reached $f$.
Thus we are out of options, giving us a contradiction and finishing the proof.
\end{proof}

\subsection{Edge-subdivided graphs}\label{ss:edge_subdivided_graphs_naive_Jaeger_is_Jaeger}

We will now consider the class of graphs mentioned in Section \ref{ss:Homfly_well_def_on_subdivided_graphs}. That is, let $G$ be be an arbitrary ribbon graph, and let $D_G$ be the signed ribbon digraph obtained from $G$ by subdividing each edge with a new vertex, orienting the edges toward the subdividing points, and defining each edge to be positive. (We shall use some notation from Section \ref{ss:Homfly_well_def_on_subdivided_graphs}.)

\begin{prop}\label{prop:naive_Jaeger_in_DG_is_Jaeger}
    If $J$ is a naive Jaeger quasi-tree of $D_G$, then $J$ is also a Jaeger quasi-tree of $D_G$.
\end{prop}

\begin{proof}
    Suppose that $J$ is a naive Jaeger quasi-tree of $D_G$. To see that $J$ is a Jaeger quasi-tree, we need to show that it does not fail the Jaeger test. There are two ways to fail, cases \ref{e_notin_K_contracted} and \ref{e_notin_K_head}. Suppose for a contradiction that $J$ does fail the Jaeger test, and consider the edge $e$ where that happens.
    As before, case \ref{e_notin_K_head} for $e$ would mean that $J$ is not a naive Jaeger quasi-tree, a contradiction. We claim that $e$ cannot fail under case \ref{e_notin_K_contracted}, either. More precisely, we claim that for the Jaeger test of $J$ within $D_G$, if some edge $e$ gets contracted in the process of forming the auxiliary graph $\AuxContr{D_G,J,f}$, then we have $e<_J f$, whence neither case \ref{e_notin_K_contracted}, nor in fact case \ref{e_in_K_contracted}, can occur.

    Suppose that the statement holds for edges $e'$ that are contracted before $e$ gets contracted. Let us prove it for $e$, and the edges that get contracted in the same step as $e$.
    An edge $e$ can only get contracted because some other edge $f$ satisfies the condition in case \ref{e_in_K_tail_active} in such a way that $e\in C$ for the nearest bad separating cycle $C$ of $\AuxContr{D_G,J,f}-f$. We show that each edge of $C$ precedes $f$ in $<_J$.
    
    A bad separating cycle is necessarily a directed cycle. However, in $D_G$ each cycle $C\subseteq D_G$ consists of pairs of edges $e_{u,uw}$, $e_{w,uw}$ that are oppositely oriented. 
    Hence for $C$ to be a directed cycle in $\AuxContr{D_G,J,f}-f$, for each such pair of edges, one of them needs to be either reversed or contracted. 
    
    Now if $e_{u,uw}$ was contracted in $\AuxContr{D_G,J,f}$, then $e_{w,uw}$ would be contracted as well, because we always contract cycles. 
    But in that case neither edge could be in $C$, showing that this case in fact does not occur.
    If one member of the pair, say $e_{u,uw}$, is reversed, then it was traversed before $f$ in the tour of $J$. Since $v_{uw}$ has degree $2$, the other edge $e_{w,uw}$ was also traversed immediately before or after $e_{u,uw}$. Note that $f\neq e_{w,uw}$, because $e_{w.uw}\in C$ but $f\notin C$. Therefore both members of the pair are indeed before $f$ in the order of $J$, which finishes the proof.
\end{proof}

\subsection*{AI declaration.}
No artificial intelligence was used for obtaining the results of this paper, and no AI was used for typesetting.

\bibliographystyle{plain}
\bibliography{Bernardi}

\end{document}